\documentclass{amsart}
\usepackage{graphicx} 
\usepackage{hyperref}
\usepackage{amsmath, amsthm, amssymb, amscd}
\allowdisplaybreaks
\usepackage{mathtools}
\usepackage{enumerate}
\usepackage{hyperref}
\usepackage{verbatim}
\usepackage{cite}
\usepackage{color}
\usepackage{esint}
\usepackage{tikz-cd}
\usetikzlibrary{shapes.geometric}
\usepackage[shortlabels]{enumitem}
 
\usepackage{bm}

\title[CLP for symmetric self-similar metric spaces]{The Combinatorial Loewner Property and super-multiplicativity inequalities for symmetric self-similar metric spaces}
\author{Riku Anttila}
\email{riku.t.anttila@jyu.fi}
\address{Department of Math. and Stat.
P.O. Box 35 \\
FI-40014 University of Jyväskylä}

\author{Sylvester Eriksson-Bique}
\email{sylvester.d.eriksson-bique@jyu.fi}
\address{Department of Math. and Stat.
P.O. Box 35 \\
FI-40014 University of Jyväskylä}
\date{\today}


\usepackage[shortlabels]{enumitem}
 
\newtheorem{theorem}[equation]{Theorem}
\newtheorem{lemma}[equation]{Lemma}
\newtheorem{proposition}[equation]{Proposition}
\newtheorem{corollary}[equation]{Corollary}
\newtheorem{assumption}[equation]{Assumption}

\numberwithin{equation}{section}

\theoremstyle{definition}
\newtheorem{definition}[equation]{Definition}

\theoremstyle{remark}
\newtheorem{remark}[equation]{Remark}

\newtheorem{example}[equation]{Example}
    \newcommand*{\N}{\mathbb{N}}
    \newcommand*{\Z}{\mathbb{Z}}
    
    \newcommand*{\R}{\mathbb{R}}
    
        \DeclarePairedDelimiter\Span{\langle}{\rangle}
        \DeclareMathOperator{\id}{id}
        \DeclareMathOperator{\diam}{diam}
        \DeclareMathOperator{\sgn}{sgn}
        \DeclareMathOperator{\dist}{dist}
        
        \DeclareMathOperator{\fold}{fold}
        \DeclareMathOperator{\divr}{div}
        \DeclareMathOperator{\len}{len}

        \DeclareMathOperator{\cRes}{\rm Res}
        \DeclareMathOperator{\cE}{\mathcal{E}}
        \DeclareMathOperator{\cT}{\mathcal{T}}
        
        \DeclareMathOperator{\ft}{\mathfrak{t}}

        \DeclareMathOperator{\cL}{\mathcal{L}}

        \DeclareMathOperator{\cM}{\mathcal{M}}
        
        \DeclareMathOperator{\cG}{\mathcal{G}}
        \DeclareMathOperator{\cH}{\mathcal{H}}

        \DeclareMathOperator{\dims}{dim}
        
        \DeclareMathOperator{\Mod}{Mod}

        \DeclarePairedDelimiter\abs{\lvert}{\rvert}
        
        \DeclarePairedDelimiter{\ceil}{\lceil}{\rceil}
        
\begin{document}
\maketitle
\begin{abstract}
This paper introduces a general construction of self-similar metric spaces as limits of discrete graphs. Our framework produces many classical examples, such as the Sierpi\'nski carpet and the higher dimensional Menger sponges, but also a rich class of new examples.
The main result of the work roughly speaking states: If the construction is sufficiently symmetric then the limiting object supports useful moduli estimates, namely the Combinatorial Loewner property of Bourdon--Kleiner and the super-multiplicativity inequalities.
The latter are established on Menger sponges for which it had not been previously known.
The main new technique the work offers is a general framework of flows and resistance estimates.
\end{abstract}

\section{Introduction}

\subsection{Overview}
This paper introduces a general construction of symmetric self-similar metric spaces as a limit of discrete graphs. We analyze the quasiconformal geometry and non-linear potential theory of these limit spaces using the fairly common tool in the related literature, the so-called \emph{combinatorial}/\emph{discrete modulus}.

We are particularly interested, in the general context, in understanding when a metric space supports a certain family of scale-invariant moduli estimates, namely the \emph{Combinatorial Loewner property (CLP)} of Bourdon--Kleiner \cite{BourK}. This question is motivated to deepen our understanding on the conjecture posed by Kleiner \cite{KleinerICM}, stating that a self-similar metric space carrying the combinatorial Loewner property is quasisymmetrically equivalent to a Loewner space of Heinonen--Koskela \cite{HK}.

The conjecture was recently disproved by the authors \cite{anttila2024constructions} but Kleiner's question, nevertheless, remains open and relevant for many interesting examples. See for instance Bourdon--Kleiner \cite{BourK} and Clais \cite{clais} for a large class of boundaries of hyperbolic groups where the conjecture is open, Bourdon--Pajot \cite{BourdonPajot} for group boundaries where the answer is positive. If a group boundary is homeomorphic to the 2-sphere then, according to the results of Bonk--Kleiner \cite{BK01,BK05}, Kleiner's question is closely related to Cannon's conjecture in geometric group theory.
See also the recent work of Murugan and Shimizu \cite{murugan2023first} which connects Kleiner's question on the Sierpi\'nski carpet, which remains open, to potential theory and Sobolev spaces.
Lastly, the recent work of the second author and G. C. David \cite{david2024analytically} shows that positive answer to Kleiner's question on the pillow space (Example \ref{ex:DS}) would lead to a quite pathological differentiability structure of the Euclidean plane.

The main goal of the present work is to produce new self-similar metric spaces and verify their combinatorial Loewner property, with the aim to support future research on Kleiner's conjecture in the general context.
We introduce new resistance/flow methods and employ them to establish delicate moduli estimates, in particular the \emph{super-multiplicativity inequality}.
Our flow techniques can be understood as a generalization of Barlow--Bass \cite{BarlowBassResistanceOnSC} and Kwapisz \cite{Kwapisz}.

\subsection{Basic framework}
Let us begin by providing a rough description of our construction of the limit spaces through some examples;
the method is quite similar to the one recently introduced by the authors \cite{anttila2024constructions}, but therein construction is only able to produce topologically one dimensional objects.
It is also unable to produce the following model example of the present work, the Sierpi\'nski carpet.
The general construction is described in Section \ref{Section:IGSs} in terms of the language of \emph{(vertex)-iterated graph systems}.

\begin{example}[\textbf{Sierpi\'nski carpet}]
\label{intro:EX}
    Sierpi\'nski carpet is a self-similar fractal which can be obtained by removing smaller squares from the unit square $[0,1]^2$ in an iterative fashion. It is also an attractor of an iterated function system in $\R^2$ consisting of eight contraction maps.
    \begin{figure}[h!]
    \centering\includegraphics[width=260pt]{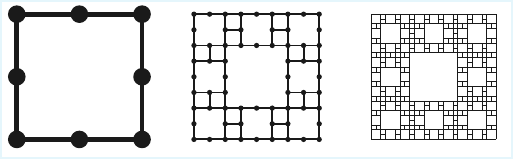} 
    \caption{Construction of the Sierpi\'nski carpet as a limit of graphs.
    }
    \label{fig:SC}
\end{figure}
    However, in the present work, we are interested in a different construction; the first steps are drawn in Figure \ref{fig:SC}. It goes as follows.
    First, take the left most graph in Figure \ref{fig:SC} which consists of eight vertices and eight edges, and name it $G_1$.
    The graph in the middle of the figure, say $G_2$, is constructed by first replacing each vertex in $G_1$ by $G_1$. Then each horizontal/vertical edge in $G_1$ is replaced by three horizontal/vertical edges. This produces $G_2$. The third graph $G_3$ is obtained by applying an analogous replacement procedure to $G_2$, namely we replace each vertex in $G_2$ by $G_1$ and edges in a simlar way. 
    By repeating, we obtain an infinite sequence $\{G_n = (V^n,E_n)\}_{n \in \N}$ of self-similar graphs.
    The limit space is obtained by taking a suitable limit of the rescaled path metrics,
    \[
        (X,d) = \lim_{n \to \infty} (V^n,3^{-n}d_{G_n})
    \]
    where the limit can be understood in a (subsequential) Gromov--Hausdorff sense.
    The limit space is biLipschitz equivalent to the Sierpi\'nski carpet.
\end{example}

Similar replacement construction can be used to produce many other classical examples. For instance, by adding dimensions to the graphs, we can produce higher dimensional carpets, namely the $d$-dimensional Menger sponges for $d \geq 3$, and also the broader class of generalized Sierpi\'nski carpets; see \cite{UniquenessofBrownianMotionOnGSC} for the definition.
We can also consider other shapes as well. The present work introduces a new pentagonal self-similar metric space.

\begin{example}[\textbf{Pentagonal carpet}]
    The replacement rule in Figure \ref{fig:PC} produces a pentagonal fractal space which we call the \emph{pentagonal carpet}.
    \begin{figure}[h!]
        \centering\includegraphics[width=300pt]{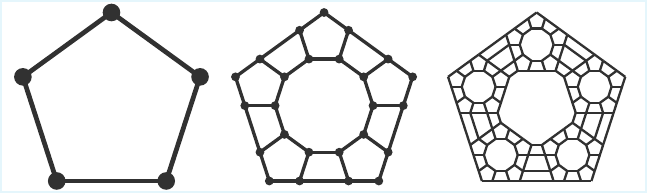} 
        \caption{First three steps in the construction of the pentagonal carpet
        }
        \label{fig:PC}
    \end{figure}
    It is homeomorphic to the Sierpi\'nski carpet, thus the name pentagonal carpet, but its geometry and analysis seems to be rather different from the standard carpet.
    A driving reason is the difference of symmetry groups.
    The limit can be obtained by taking a suitable rescaled limit,
    \[
      (X,d) = \lim_{n \to \infty} (V^n,2^{-n}d_{G_n}).
    \]
\end{example}

We are also interested in the following fractal space which can be understood as a two-dimensional variant of the Laakso diamond space; see \cite{Laakso,LangPlaut,AndreaSchioppa2015,David_Schul_2017} for some studies on the Laakso diamond.

\begin{example}[\textbf{Pillow space}]
\label{ex:DS}
    The \emph{Pillow space} is self-similar metric space that can be constructed in a similar manner as the Sierpi\'nski carpet with the difference that the middle squares are doubled instead of being removed. Of course, such a procedure is impossible to perform inside the Euclidean plane. Instead, the construction is done in an abstract setting.
    \begin{figure}[h!]
        \centering\includegraphics[width=240pt]{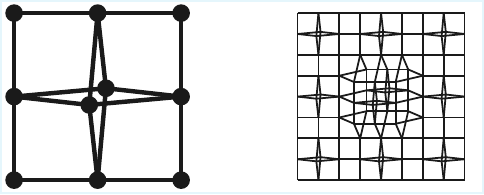} 
        \caption{First two steps in the construction of the pillow space.
        }
        \label{fig:DS}
    \end{figure}
    This idea was inspired by the Laakso diamond space, which can be constructed via similar doubling procedure of intervals.
    In Figure \ref{fig:DS} we construct the pillow space as a limit of graphs. The limit is, again, defined as a rescaled limit of graphs,
    \[
      (X,d) = \lim_{n \to \infty} (V^n,3^{-n}d_{G_n}).
    \]
\end{example}

Section \ref{Section:IGSs} provides the precise definitions and detailed descriptions of the constructions, and the rest of the main examples are introduced in Section \ref{Section: Main examples}.

\subsection{Main results}
First, we need some terminology; we remark that some of the notation and definitions below are slightly different from what we use in later sections for expository reasons.
Consider a general compact path-connected metric space $(X,d)$.
For a given positive integer $n \in \N$ we fix a $2^{-n}$-net $N_n \subseteq X$, meaning the balls $B(x,2^{-n})$ for $x \in N_n$ cover $X$ and $d(x,y) \geq 2^{-n}$ for all distinct $x,y \in N_n$. Let us denote the collections of open balls $G_n := \{ B(x,2^{-n}) \}_{x \in N_n}$. Note that $G_n$ is always finite because $(X,d)$ is compact.

Then take any family of (continuous and compact) curves $\Gamma$ in $X$.
By the term \emph{curve} we refer to a subset $\gamma \subseteq X$ which is an image of a continuous map $[0,1] \to X$. Now, the \emph{discrete modulus} of $\Gamma$ is defined as follows.
For a given function $\rho : G_n \to [0,\infty)$, we say that $\rho$ is \emph{$\Gamma$-admissible} if for all $\gamma \in \Gamma$,
\[
    \sum_{ \substack{v \in G_n \\ \gamma \cap v \neq \emptyset} } \rho(v) \geq 1.
\]
Then, for a given $p \geq 1$, the \emph{$G_n$-discrete $p$-modulus} is
\[
    \Mod_p(\Gamma,G_n) := \inf_{\rho} \sum_{v \in G_n} \rho(v)^p
\]
where the infimum is taken over all $\Gamma$-admissible functions $\rho$.
The value $\Mod_p(\Gamma,G_n)$ in general depends on the parameters $p,\,\Gamma$ and $G_n$, but the dependence on $G_n$ is very mild; 
see \cite[Proposition 2.2]{BourK} which is also stated in Proposition \ref{prop: Modulus change levels} of the present work.

The most important discrete moduli for the present work are
\begin{equation}\label{Intro:CombMod}
    \cM_{p,n} := \Mod_p(\Gamma_\delta,G_n) \text{ where } \Gamma_\delta := \{ \gamma : \diam(\gamma) \geq \delta \},
\end{equation}
where $\delta > 0$ is some fixed sufficiently small constant.
We are interested in the asymptotic behavior of $\cM_{p,n}$ as $n \to \infty$.
By definition, $\Gamma_\delta$ contains the large scale geometric data of the ambient space $(X,d)$, but fails to see the scales smaller than $\delta$. However, in a setting where small scale geometry is identical to large scale, or in other words when the ambient space is self-similar, $\Gamma_\delta$ is particularly useful. For instance, Carrasco Piaggio and Keith--Kleiner showed that, by analyzing the asymptotic behavior of $\cM_{p,n}$ for different $p \geq 1$, one can determine the \emph{conformal dimension} of the ambient space; see \cite[Corollary 1.4]{Carrasco}, and this is stated in Proposition \ref{prop:confdimchar} of the present work as-well. We note that the result of Carrasco Piaggio is far more general than we indicate here.
See also \cite{Shaconf,Kigamiweighted,EBConf,MurCA} for related results.

We now state the first main result of the work, which regards an important regularity estimate of the discrete moduli, the \emph{super-multiplicativity inequality}; see Theorem \ref{thm:supermult} for the general result and the proof.
\begin{theorem}\label{introTHM:supermult}
    Let $(X,d)$ be a $d$-dimensional Menger sponge for $d \geq 3$, the pentagonal carpet, the pillow space or more generally any metric space satisfying the assumptions in Theorem \ref{thm:supermult}. For $p \geq 1$ let $\cM_{p,n}$ be as we defined above. Then the following super-multiplicativity inequality holds. There is a constant $C \geq 1$ such that for all $n,m \in \N$,
    \begin{equation}\label{intro:Supermult}
    C^{-1} \cM_{p,n}\cdot \cM_{p,m} \leq \cM_{p,n+m} \leq C \cM_{p,n}\cdot \cM_{p,m}.
    \end{equation}
    The constant $C$ depends on the ambient space and continuously on $p \in [1,\infty)$.
\end{theorem}

Frankly, super-multiplicativity only refers to the lower inequality, and the upper inequality is called \emph{sub-multiplicativity inequality}. We nevertheless wrote both for expository reasons, but also because sub-multiplicativity is virtually always true in self-similar setting; see \cite[Proposition 1.5]{BourK}.

Super-multiplicativity has numerous applications.
Let us review a few notable ones; we provide further details in Section \ref{sec: Some further results}.
Bourdon--Kleiner proved the super-multiplicativity inequality for the Sierpi\'nski carpet, and used it to establish the Combinatorial Loewner property.
They also verified the Combinatorial Loewner property for the higher dimensional Menger sponges, but using a weaker variant of the super-multiplicativity; see \cite[Section 4]{BourK}. The present work resolves super-multiplicativity of Menger sponges for all $p \geq 1$.
We note that, for exclusively $p = 2$, this was already known by \cite{Mcgillivray}.

In the recent work of Murugan and Shimizu \cite{murugan2023first}, where they constructed Sobolev spaces on the Sierpi\'nski carpet using the Kusuoka--Zhou apporach \cite{KusuokaZhou}, the main reason why their methods did not generalize to Menger sponges was the lack of super-multiplicativity; see \cite[Remark 8.24]{murugan2023first}. By combining the results of the present work with the techniques of Murugan--Shimizu, it should be possible to construct Sobolev spaces also on Menger sponges. 
However, we note that some of their techniques have a certain dimension restriction, which is a problem when $p \in (1,\infty)$ is very small; see \cite[Problem 10.1]{murugan2023first}. 
In $\R^n$ this restriction can be understood as $p > n -1$ and is related to the fact that in this range curves have positive capacity.
The present work does not resolve this deep issue.
We also note that, with even a more restrictive dimensional condition, namely $p$ being strictly larger than the conformal dimension, the framework of Kigami \cite{kigami} and the earlier work of Shimizu \cite{shimizu} already provides good Sobolev spaces also on the Menger sponges. This dimensional condition corresponds to $p > n$ in $\R^n$, which is supported by the validity of Morrey-type inequality \cite[Theorem 3.21]{kigami}.

The last application we mention is a special case of what we discussed in the previous paragraph. When $p = 2$ and the dimensional restriction from Murugan--Shimizu \cite[Problem 10.1]{murugan2023first} is satisfied, one can perform an analytic construction of diffusion processes and Dirichlet forms. In particular, the method works on the 3-dimensional Menger sponge. 
This construction would coincide with the Brownian motion constructed by Barlow--Bass using a probabilistic method \cite{barlow1989construction,barlow1999brownian}, because a self-similar diffusion process on the Menger sponge is unique according to \cite{UniquenessofBrownianMotionOnGSC}. We note that the probabilistic method of Barlow--Bass does not suffer from the dimension restriction of \cite{murugan2023first}.

Next, we state our main result regarding the combinatorial Loewner property. The precise definition is postponed to Section \ref{Section: Preliminary}, but roughly speaking the combinatorial Loewner property is a quantitative formulation of the general principle that discrete modulus is scale/conformally invariant. More specifically, a metric space $(X,d)$ is \emph{combinatorially $Q$-Loewner} for $Q>1$ if for any disjoint pair of connected and compact subsets $E,F \subseteq X$ the discrete moduli $\Mod_Q(\Gamma(E,F),G_n)$ is controlled by the \emph{relative distance}
\[
    \Delta(E,F) := \frac{\dist(E,F)}{\diam(E) \land
     \diam(F)}.
\]
Here $\Gamma(E,F)$ denotes the continuous curves connecting $E$ and $F$.
This definition was introduced in Bourdon--Kleiner \cite{BourK} and it was inspired by the Loewner spaces of Heinonen--Koskela \cite{HK}. Unlike the Loewner condition, the combinatorial Loewner proeprty is a quasisymmetry invariant; see \cite[Theorem 2.6]{BourK}.

\begin{theorem}
    Let $(X,d)$ be the pentagonal carpet, the pillow space or more generally any metric space satisfying the assumption of Theorem \ref{thm: CLP}.
    Then there is a unique $Q > 1$ such that
    \[
        C^{-1} \leq \liminf_{n \to \infty} \cM_{Q,n} \leq \limsup_{n \to \infty} \cM_{Q,n} \leq C
    \]
    for some constant $C \geq 1$. This $Q$ is equal to the conformal dimension of $(X,d)$ and $(X,d)$ satisfies the combinatorial $Q$-Loewner property.
\end{theorem}

\subsection{Main techniques}
The main new technique is the \emph{replacement flow method}.
The general idea and the terminology ``replacement flow'' was inspired by the work of Kwapisz \cite{Kwapisz} where he derived impressive estimates on the conformal dimension of the Sierpi\'nski carpet by studying flows on the graph approximations. See also the related work \cite{malo2015discrete}.
The key idea is that flows provide lower bounds to discrete moduli through the standard duality; see \cite[Theorem 5.1]{NakamuraYamasakiDuality} and \cite{ACFPC}. We also state a variant of it in Proposition \ref{prop: Duality}.

Flow methods on Sierpi\'nski carpets go back to Barlow--Bass \cite{BarlowBassResistanceOnSC} which was later generalized to higher dimensional carpets by McGillivray \cite{Mcgillivray}; see also \cite{Resistance4Ncarpets}.
These works exclusively regard the case $p = 2$ and some of their methods use the classical Harmonic analysis of the Euclidean space. It seems unclear how well their techniques generalize to general exponents $p \in (1,\infty)$, let alone to cases where the ambient space is not a subset of the Euclidean space.
We note that, because we are interested in the Combinatorial Loewner property, we in particular need these techniques for $p$ equal to the conformal dimension.
Kwapisz applied the flow method for general exponents $p$ but exclusively in the case of the planar Sierpi\'nski carpet.

The replacement flow method is by far the most technical part of the paper. But the basic idea is rather simple. Consider three graphs $G_n,\, G_m$ and $G_{n+m}$ described in Example \ref{intro:EX}, for $n,m \in \N$. Assume we are given a flow $\mathcal{J}_n$ from $S_n$ to $T_n$ where $S_n,\, T_n \subseteq V(G_n)$ is a pair of disjoint subsets. The idea of the replacement flow is to produce a flow on $\mathcal{J}_{n+m}$ on $G_{n+m}$ by suitably ``lifting'' $\mathcal{J}_n$.
This is performed by first constructing a family of flows $\mathcal{B}$ on $G_m$ so that for every pair of the four line segments (as in Figure \ref{fig:SC}) in the boundary, there is a flow between them.
In later sections we call $\mathcal{B}$ a \emph{flow basis}.
Then the flow $\mathcal{J}_{n+m}$ is obtained by suitably embedding the flows on $G_{m}$ into $G_{n+m}$ by using the self-similarity of the graphs.
We take suitably linear combinations of them and the coefficients are determined by $\mathcal{J}_n$. The construction of a flow basis is a bit delicate, and involves substantial use of symmetries to bend and reflect flows.

\subsection*{Organization}
In Section \ref{Section: Preliminary} we review some key concepts and results regarding conformal dimension.

Section \ref{Section:IGSs} introduces the definitions and main techniques used for producing the self-similar metric space, namely the general framework of \emph{(vertex)-iterated graph systems} and \emph{replacement graphs}.
In Section \ref{Section: Main examples} we review the main examples of IGSs and their limit spaces this work is interested in.

Section \ref{Section: Replacement flow} regards flows on replacement graphs. We establish the most important ingredient of the work, namely the replacement flow technique.

Section \ref{Section: Bounded geometry} studies asymptotic geometry of the replacement graphs that seems to be crucial for constructing a good limit space.
Section \ref{sec:LS} introduces the limit spaces and studies their basic properties under certain axioms.

Finally in Section \ref{Section: CLP}, we prove our main results, Theorems \ref{thm:supermult} and \ref{thm: CLP}. Some further results related to Sobolev spaces are discussed in Section \ref{sec: Some further results}.

\subsection*{Acknowledgements} We thank the referee for a very careful reading, that encouraged us to substantially improve the paper from its original draft. 
We thank Pietro Poggi-Corradini, Hrant Hakobyan and Ryosuke Shimizu for helpful discussions, and give a special thanks to Shimizu for his comments on an earlier version of the paper. The first author was funded by Jenny and Arttu Wihuri Foundation, Eemil Aaltonen Foundation, University of Oulu and University of Jyväskylä. The second author
was supported by the Research Council of Finland grants 354241, 345005 and 356861.

\section{Preliminaries}\label{Section: Preliminary}
\subsection{Conformal dimension}
Consider complete metric spaces $(X,d)$. All metric balls $B(x,r):=\{y\in X : d(x,y)<r\}$ are open for $x\in X, r>0$. 
We usually refer to the metric space $(X,d)$ simply by $X$ and omit the metric $d$.

\begin{definition}
    Let $Q \in (0,\infty)$. We say that a metric space $X$ is $Q$-\emph{Ahlfors regular}, if there exists a Radon measure $\mu$ of $(X,d)$ and constant $C \geq 1$ such that for all $r\in (0,2\diam(X))$ and all $x\in X$
\[
    C^{-1} \cdot r^Q \leq \mu(B(x,r)) \leq C \cdot r^Q.
\]
If a metric space is $Q$-Ahlfors regular for some $Q \in (0,\infty)$, we simply say that it is \emph{Ahlfors regular}.
A metric space $(X,d)$ is \emph{metric doubling}, if there is a constant $N\in \N$ such that for every ball $B(x,r)\subset X$ there are $x_1,\dots, x_N \in X$ for which $B(x,r)\subset \bigcup_{i=1}^N B(x,r/2)$. A simple volume counting argument shows that every Ahlfors regular metric space is metric doubling. 
It is also well-known that if $X$ is Ahlfors regular, then any measure $\mu$ as above is comparable to the $Q$-\emph{Hausdorff measure} $\cH^Q$ given by
\[
    \cH^Q(A) := \lim_{\delta \to 0} \cH^Q_\delta(A),
\]
where $\cH^Q_\delta$ is the \emph{Hausdorff content}
\[
    \cH^Q_\delta(A):=\inf \left\{\sum_{i=1}^\infty \diam(A_i)^Q : A \subset \bigcup_{i=1}^\infty A_i, \diam(A_i)\leq \delta \right\}.
\]
In particular, the \emph{Hausdorff dimension} of $(X,d)$ then is
\[
    Q = \dim_H(X)  := \inf \left\{ d > 0 : \cH^d(X) = 0 \right \}.
\]
\end{definition}

\begin{definition}
    A metric space $X$ is \emph{quasiconvex} if there is a constant $C \geq 1$ such that for every pair $x,y \in X$ there is a continuous curve $\gamma : [0,1] \to X$ with $\gamma(0) = x, \gamma(1) = 1$ and $\len(\gamma) \leq C \cdot d(x,y)$. Here the \emph{length} of $\gamma$ is given by
    \[
       \len(\gamma) := \sup \left\{ \sum_{i = 1}^{N - 1} d(\gamma(x_i),\gamma(x_{i+1})) : 0 = x_1 < x_2 < \dots < x_N = 1   \right\}.
    \]
\end{definition}

Next, we introduce the notion of self-similarity relevant to our work. See \cite{Carthesis} for more background.
\begin{definition}
If $L \geq 1$ a map $f:(X,d_X)\to (Y,d_Y)$ is an $L$-\emph{biLipschitz map} if for all $x,y\in X$
\[
\frac{1}{L} d_X(x,y) \leq d_Y(f(x),f(y)) \leq L d_X(x,y).
\]
A metric space $X$ is \emph{approximately self-similar}, if there is a constant $L\geq 1$ such that for every $x\in X$ and every $r\in (0,\diam(X)]$, there exists an open set $U_{x,r}\subset X$ and a $L$-biLipschitz map $f_{x,r}:(B(x,r),d/r)\to U_{x,r}$. Here the ball $B(x,r)$ is equipped with the scaled metric $d/r$. 
\end{definition}
If $L=1$, such a biLipschitz scaling map is called a homothetic homeomorphism with scaling factor $1/r$.

\begin{definition}
Given a homeomorphism $\eta:[0,\infty)\to [0,\infty)$, we say that homeomorphism $f:(X,d_X)\to (Y,d_Y)$ is an $\eta$-\emph{quasisymmetry} if for every $x,y,z\in Z$ with $x\neq z$, we have
\[
\frac{d_Y(f(x),f(y))}{d_Y(f(x),f(z))} \leq \eta\left(\frac{d_X(x,y)}{d_X(x,z)}\right).
\]
We say that $f:X\to Y$ is a \emph{quasisymmetry}, if it is an $\eta$-quasisymmetry for some $\eta$. In these cases, we say that $X$ is quasisymmetric to $Y$. For more background on these mappings and the notions below, see \cite{He}. We will usually omit the subscript of a metric $d$, where the space is clear from context.
The \emph{conformal gauge} of a metric space $X$ is given by
\[
\cG(X,d):= \left\{d' : d' \text{ is a metric on } X \text{ and } \id:(X,d)\to (X,d') \text{ is a quasisymmetry} \right\}.
\]
We usually drop the metric from $\cG(X,d)$ and simply write $\cG(X)$.
\end{definition}

Understanding the structure of the conformal gauge allows one to study different types of \emph{quasisymmetric uniformization problems}, i.e., finding a metric in the conformal gauge more suitable for the problem in question. An important example is  to minimize the Hausdorff dimension among Ahlfors regular metrics in the conformal gauge. From this problem, naturally arises the notion of \emph{(Ahlfors regular) conformal dimension}, which is given by
\[
\dims_{\rm AR}(X):=\inf \left\{ Q' : d' \in \cG(X,d) \text{ and } (X,d') \text{ is $Q'$-Ahlfors regular} \right\}.
\]
Determining whether the above infimum is a minimum is called \emph{the attainment problem of (Ahlfors regular) conformal dimension}.

The Ahlfors regular metrics in $\cG(X)$ can be described by quasi-isometries on a certain hyperbolic space, and the conformal dimension is characterized as the critical dimension of suitable discrete moduli problems; see e.g. \cite{Kigamiweighted,Shaconf,Carrasco}. Despite this, the attainment problem is notoriously hard, and the answer is known only in some exceptional cases. Finding the exact value of the conformal dimension also tends to be challenging. See \cite{Kwapisz} for numerical estimates of the conformal dimension of the Sierpi\'nski carpet, whose exact value of the conformal dimension remains unknown.
See \cite{anttila2024constructions} for examples where the conformal dimension is computable.

Other variants of the conformal dimensions have also been studied. The conformal Hausdorff dimension was first considered by Pansu \cite{P89}, the conformal Assouad dimension was studied in \cite{MurCA}, the conformal walk-dimension was introduced in \cite{MN}, and the present one, the Ahlfors regular conformal dimension, was introduced in \cite{BourdonPajot}. For further background on the conformal dimension, see \cite{MT}. The various definitions (except the conformal walk-dimension) are equivalent in self-similar settings; see \cite{EBConf}. Since all the spaces considered in this paper are self-similar, we refer to the Ahlfors regular conformal dimension simply as the conformal dimension.

\subsection{Terminology of graphs}\label{subsec: Graph Theory}
A \emph{graph} is a pair $G = (V,E)$ where $V(G) := V$ is a non-empty finite set of \emph{vertices} and $E(G) := E \subseteq V \times V$ is the set of edges.
The set of edges are always assumed to satisfy: If $(x,y) \in E$ then $(y,x) \notin E$. In particular, $(x,x) \notin E$ for all $x \in V$.

Usually, the above definition of a graph is known as oriented simple graph. In our work the given orientations are more or less ``artificial'' in the sense that they are needed only for describing certain constructions in Subsection \ref{subsec:RR} and are mostly disregarded afterwards. Whenever we only want to emphasize that a pair of vertices is connected by an edge and the orientation is not important, we write the edge as $\{x,y\}$. If the choice of the orientation does not matter, we sometimes leave it unspecified which orientation we have chosen.

For a vertex $x \in V$ we define its \emph{degree} as the number of its neighbors, meaning
\[
    \deg(x) := \abs{ \{ y \in V(G) : \{x,y\} \in E(G) \} }.
\]
The \emph{degree} of $G$, denoted $\deg(G)$, is then the maximum of the degrees of its vertices.

Next we discuss concepts related to discrete modulus; see e.g. \cite{ShimizuParabolicIndex,ACFPC} for background.
First, a \emph{path} in a graph $G$ is a sequence $\theta = [x_1,\dots,x_k]$ and $\{ x_i,x_{i + 1} \} \in E$ for all $i = 1,\dots,k - 1$ when $k > 1$. If $k = 1$ then $\theta$ consists of only one vertex. The \emph{length} of a path $\theta$ then is $\len(\theta) := k - 1$.
Given non-empty subsets $A,B \subseteq V(G)$ we denote
\[
    \Theta(A,B,G) := \{ \theta : \theta = [x_1,\dots,x_n] \text{ is a path in } G \text{ and } x_1 \in A, x_n \in B   \}.  
\]
If $\Theta(A,B,G) \neq \emptyset$ whenever $A,B \neq \emptyset$, we say that $G$ is \emph{connected}, and denote the \emph{(shortest) path metric} by $d_G$.
We also say that a subset $U \subseteq V(G)$ is \emph{connected} if the induced subgraph $G[U]$ is connected. Furthermore, given a path $[x_1,\dots,x_k]$, without an explicit mention, we sometimes regard it as the connected subset $\{ x_1,\dots,x_k \} \subseteq G$.

\begin{definition}
    Let $G$ be a graph and $\rho : V(G) \to [0,\infty)$. If $\theta$ is a path we define its \emph{$\rho$-length}
    \[
        L_\rho(\theta) := \sum_{x \in \theta} \rho(x),
    \]
    namely we sum over all vertices appearing in the path $\theta$.
    Note that each vertex is included at most once in the sum even if it appears multiple times in $\theta$.
    Given a family of paths $\Theta$ in $G$ we say that $\rho : V(G) \to \R_{\geq 0}$ is \emph{$\Theta$-admissible} if $L_\rho(\theta) \geq 1$ for all $\theta \in \Theta$.
    Given a density $\rho : V \to \R_{\geq 0}$ and $p \geq 1$, the \emph{$p$-mass} of $\rho$ is given by
    \[
        \cM_p(\rho) := \sum_{x \in V} \rho(x)^p.
    \]
    The \emph{discrete $p$-modulus} of $\Theta$ is then given by
    \[
        \Mod_p(\Theta,G) := \inf_{\rho} \cM_p(\rho),
    \]
    where the infimum is over all $\Theta$-admissible densities. For simplicity we write $\Mod_p(A,B,G) := \Mod_p(\Theta(A,B,G),G)$.
\end{definition}

Discrete modulus is a convex minimization problem. When $p > 1$ then the minimum always exists and is unique. 
Because it is a minimization problem, obtaining upper bound estimates is straightforward by constructing admissible densities with suitable mass.
Lower bounds are more involved. Our method for deriving these is to analyze the dual problem of resistance and flows; see \cite{ShimizuParabolicIndex,ACFPC,NakamuraYamasakiDuality} for background.

\begin{definition}
    A function $\mathcal{J} : V \times V \to \R$ is \emph{antisymmetric} if it satisfies the following conditions.
    \begin{enumerate}
        \item $\mathcal{J}(x,y) = -\mathcal{J}(y,x)$ for all $x,y \in V$.
        \item $\mathcal{J}(x,y) = 0$ unless $\{ x,y \} \in E$.
    \end{enumerate}
    The \emph{divergence} of $\mathcal{J}$ at $x \in V$ is given by
    \[
        \divr(\mathcal{J})(x) := \sum_{\{ x,y \} \in E} \mathcal{J}(x,y).
    \]
    Given non-empty disjoint subsets $A,B \subseteq V$, an antisymmetric function $\mathcal{J}$ is a \emph{flow} from $A$ to $B$ if
    \[
        \divr(\mathcal{J})(x) = 0 \text{ for all } x \notin A \cup B.
    \]
    The \emph{total flow} is given by
    \[
        I(\mathcal{J}) := \sum_{x \in A} \divr(F)(x).
    \]
    If $I(\mathcal{J}) = 1,$ we say that $\mathcal{J}$ is a \emph{unit flow} from $A$ to $B$. 
    For $p \in (1,\infty)$, the \emph{$p$-energy} of a flow $\mathcal{J}$ is
    \[
        \cE_q(\mathcal{J}) := \sum_{\{ x,y \} \in E } \abs{\mathcal{J}(x,y)}^q,
    \]
    where $q = p/(p-1)$.
    The \emph{$p$-resistance} between $A$ and $B$ is
    \[
        \cRes_p(A,B,G) := \inf_{\mathcal{J}} \cE_q(\mathcal{J}),
    \]
    where the infimum is over all unit flows from $A$ to $B$.
\end{definition}

There is a well-known duality result between resistances and modulus; see \cite[Theorem 5.1]{NakamuraYamasakiDuality} and also \cite{ACFPC}. In our setting, it reads as follows.

\begin{proposition}\label{prop: Duality}
Let $G$ be a connected graph, $A,B \subseteq V(G)$ be disjoint non-empty subsets and $p \in (1,\infty)$.
Write $C = 2^p\deg(G)$. Then
\begin{equation}\label{eq: Duality}
    C^{-1} \cRes_p(A,B,G)^{1-p} \leq \Mod_p(\Theta(A,B),G) \leq C \cRes(A,B,G)^{1-p}
\end{equation}
\end{proposition}

\begin{remark}
    The standard version of duality \eqref{eq: Duality} concerns the edge modulus, instead of the vertex counterpart we use. The above version of the duality follows from the fact that edge-modulus and vertex-modulus are comparable;
    see \cite[Lemma 4.3]{anttila2024constructions} for details.
\end{remark}

In our analysis, it will be convenient to expand the graph at sources and targets of the flow. The following lemma provides control over this procedure.

\begin{lemma}\label{lemma: Resistance after joining}
    Let $G$ be a connected graph, $S,T \subseteq V(G)$ be disjoint non-empty subsets and $\mathcal{J}$ be a flow from $S$ to $T$ in $G$.
    We define an expanded graph $G^*$, which is obtained by adding the vertices
    \[
        V(G^*) :=  V(G) \cup \{ x^* : x \in A \cup B \},
    \]
    and the edges
    \[
        E(G^*) := E(G) \cup \{ \{ x^*,x \} : x \in A \cup B \}.
    \]
    We denote $S^* := \{ x^* : x \in S \} \text{ and } T^* := \{ x^* : x \in T \}$.
    Now, let $J^*$ be a flow from $S^*$ to $T^*$ be a flow obtained by extending $\mathcal{J}$ by
    \[
       \mathcal{J}^*(x^*,x) := \divr(\mathcal{J})(x) = \sum_{\{ x,y \}\in E} \mathcal{J}(x,y).
    \]
    Then there is a constant $C \geq 1$ depending only on $p \in (1,\infty)$ and $\deg(G)$ such that 
    \begin{equation}
        \cE_q(\mathcal{J}^*) \leq  3\deg(G)^{\frac{1}{p-1}}\cE_q(\mathcal{J})
    \end{equation}
    where $q = p/(p-1)$.
\end{lemma}

\begin{proof}
    The claim follows from Hölder's inequality,
    \begin{align*}
        \cE_q(\mathcal{J}^*) & = \sum_{\{ x,y \} \in E} \abs{\mathcal{J}^*(x,y)}^q + \sum_{x \in S \cup T} \left|\widehat{\mathcal{J}}(x^*,x) \right|^q\\
        & = \sum_{\{ x,y \} \in E} \abs{\mathcal{J}(x,y)}^q + \sum_{x \in S \cup T} \left|\sum_{ \{ x,y \} \in E} \mathcal{J}(x,y) \right|^q\\
        & \leq \cE_q(\mathcal{J}) + \deg(G)^{\frac{1}{p-1}}\sum_{x \in S \cup T} \sum_{ \{ x,y \} \in E } \abs{\mathcal{J}(x,y)}^q\\
        & \leq (1 + 2\deg(G)^{\frac{1}{p-1}}) \cE_q(\mathcal{J})\\
        & \leq 3\deg(G)^{\frac{1}{p-1}} \cE_q(\mathcal{J}).
    \end{align*} 
\end{proof}

\subsection{Graph approximations}\label{Subsec: Discrete modulus and conformal dimension}

Let $X$ be metric space.  We say that a graph $G = (V,E)$ is an \emph{incidence graph} of $X$ if $V$ is a covering of $X$ and $\{ u,v \} \in E$ if and only if $u \cap v \neq \emptyset$ and $u \neq v$. For a given subset $F \subseteq X$ we write
\[
    V[F] := \{ v \in V : v \cap F \neq \emptyset \}.
\]
If $\Gamma$ is a family of subsets of $X$ we define
\[
    V[\Gamma] := \{ V[\gamma] : \gamma \in \Gamma \}.
\]
In our analysis $\Gamma$ is always a collection of (continuous and compact) \emph{curves}, meaning a collection of subset $F = \gamma([0,1]) \subseteq X$ where $\gamma$ is a continuous map $\gamma : [0,1] \to X$.
When $F_1,F_2 \subseteq X$ are two non-empty sets we denote by $\Gamma(F_1,F_2)$ the set of continuous curves $\gamma : [0,1] \to X$ such that $\gamma(0) \in A$ and $\gamma(1) \in B$.
We often identify a continuous curve $\gamma$ by their image.

For $\alpha,L_*>1$, an \emph{$\alpha$-approximation} to a compact metric space $X$ is a sequence of incidence graphs $\{G_n=(V^n,E_n)\}_{n \in \N}$, such that each graph $G_n$ satisfies the following conditions.
\begin{enumerate}
    \vspace{2pt}
    \item For every $v\in V^n$ there exists a $z_v\in v$ for which \[
    B(z_v, \alpha^{-1}L_*^{-n}) \subset v \subset B(z_v, \alpha L_*^{-n}).
    \]
    \item For every pair of distinct $v,w\in V^n$ we have 
    \[B(z_v,\alpha^{-1}L_*^{-n})\cap B(z_w, \alpha^{-1}L_*^{-n})=\emptyset.\] 
\end{enumerate}

In \cite{BourK}, the value $L_*=2$ is used exclusively, but by adjusting the levels $n\in \N$ (by repeating or skipping indices), and the constant $\alpha$, any $\alpha$ approximation with $L_*>1$ also yields an $\alpha'$-approximation with $L_*=2$ for some $\alpha'>1$. In the following, it will be however convenient to allow for general $L_*$. 

Fix $n\in \N, p \geq 1$ and let $\rho:V^n\to \R_{\geq 0}$. Given a family of continuous curves $\Gamma$ in $X$, we say that $\rho$ is \emph{$\Gamma$-admissible} if
\[
    L_\rho(\gamma):=\sum_{v \in V^n[\gamma]} \rho(v) \geq 1
\]
for all $\gamma \in \Gamma$. Here, the curves are identified by their image.
We define the \emph{$G_n$-discrete $p$-modulus} of $\Gamma$ by 
\[
\Mod_p(\Gamma, G_n):=\inf_{\rho} \cM_p(\rho),
\]
where the infimum is over all $\Gamma$-admissible densities $\rho : V^n \to \R_{\geq 0}$.
If $\Gamma = \Gamma(F_1,F_2)$ for $F_1,F_2 \subseteq X$ we denote $\Mod_p(F_1,F_2,G_n) := \Mod_p(\Gamma(F_1,F_2), G_n)$.
We note that this notion of modulus is very similar to the one discussed the previous subsection but with the differences that the curves here live in a general metric space rather than on a graph.

The same metric space can have many approximations but a different choice only affects the discrete moduli by a multiplicative constant.

\begin{proposition}[Proposition 2.2 \cite{BourK}]\label{prop: Modulus change levels}
    Let $\alpha,\alpha',L_* > 1$ and $G_n$ and $G_n'$ be $\alpha$- and $\alpha'$-approximations of a compact metric space $X$ respectively, and assume $X$ satisfies the metric doubling property. Then for all $p \geq 1$ there is a constant $C \geq 1$ depending only on $\alpha,\alpha'$, continuously on $p$ and the doubling constant of $X$, such that for all family of paths $\Gamma$
    \begin{equation}\label{eq: Modulus change levels}
        C^{-1} \cdot \Mod_p(\Gamma,G_n') \leq \Mod_p(\Gamma,G_n) \leq  C \cdot \Mod_p(\Gamma,G_n').
    \end{equation}
\end{proposition}

Keith and Kleiner described the conformal dimension of self-similar spaces as a critical exponent for discrete modulus. Carrasco Piaggio gave an independent proof of this fact in \cite{Carrasco}. This idea  appeared already in Pansu's original paper, where conformal dimension is actually defined as a critical dimension for a certain discrete modulus \cite{P89}. For $\delta>0$ we write $\Gamma_\delta =\{\gamma : \diam(\gamma) \geq \delta\}$  and 
\begin{equation}\label{eq:moddeltadef}
\cM_{p,\delta}^{(n)} := \Mod_p(\Gamma_{\delta}, G_n).
\end{equation}

\begin{proposition}[Corollary 1.4 \cite{Carrasco}]\label{prop:confdimchar} Let $X$ be an approximately self-similar path-connected metric space. There exists $\delta_0>0$ such that for all $\delta\in (0,\delta_0)$ we have
\[
\dims_{\rm AR}(X)=\inf\left\{p \geq 1 : \lim_{n\to\infty} \cM_{p,\delta}^{(n)} =0 \right\}.
\]
\end{proposition}

This characterization yields a numerical way to estimate conformal dimension. However, since the definition involves a limit, we need some further estimates to give rigorous bounds. Upper bounds for conformal dimension can be obtained through an argument using the sub-multiplicative inequality \eqref{eq: Sub-mult}. Kwapisz used this strategy in \cite{Kwapisz} and obtained good numerical estimates on the conformal dimension of the standard Sierpi\'nski carpet.

\begin{proposition}[Proposition 3.6 \cite{BourK}]\label{prop: Sub-mult}
    Let $X$ be compact, path-connected and approximately self-similar metric space. There exists $\delta_0 > 0$ such that for all $p \in [1,\infty)$ and $\delta \in (0,\delta_0)$ there is a constant $C$ depending on $\delta$ and continuously on $p \in [1,\infty)$ such that for all $n,m \in \N$,
    \begin{equation}\label{eq: Sub-mult}
        \cM_{p,\delta}^{(n+m)} \leq C \cM_{p,\delta}^{(n)}\cM_{p,\delta}^{(m)}.
    \end{equation}
    for all $n,m \in \N$.
\end{proposition}

In Section \ref{Section: CLP} we consider the reverse inequality, the \emph{super-multiplicative inequality}.  In contrast to sub-multiplicativity, which is virtually always true in self-similar settings, establishing super-multiplicativity will be much more involved and requires additional assumptions.

\begin{definition}\label{def:CLP}
    Let $Q, L_*,\alpha > 1$ and $\{ G_n \}_{n \in \N}$ be an $\alpha$-approximation of a compact approximately self-similar metric space.
    We say that $X$ satisfies the \emph{combinatorial $Q$-Loewner property} if there exists two positive increasing functions $\phi,\psi$ on $(0,\infty)$ with $\lim_{t \to 0} \psi(t) = 0$, which satisfies the following two conditions.
    \begin{enumerate}[label={\color{blue}{\textup{(CLP\arabic*})}}, widest=a, leftmargin=*]
        \vspace{3pt}
        \item \label{CLP1} If $F_1,F_2 \subseteq X$ are two disjoint non-degenerate continua such that $L_*^{-n} \leq \diam(F_1) \land \diam(F_2)$, then for all $m \in \N$
        \begin{equation*}
            \phi(\Delta(F_1,F_2)^{-1}) \leq \Mod_{Q}(F_1,F_2, G_{n + m}),
        \end{equation*}
        where $\Delta(F_1,F_2)$ is the \emph{relative distance}
        \[
            \Delta(F_1,F_2) := \frac{\dist(F_1,F_2)}{\diam(F_1) \land \diam(F_2)}.
        \]
        \item \label{CLP2} If $x \in X$, $r \geq L_*^{-n}$ and $C > 0$ then for all $m \in \N$
        \begin{equation*}
            \Mod_{Q}\left(B(x,r),X \setminus B(x,(C+1)r)\right), G_{n + m}) \leq \psi(C^{-1}).
        \end{equation*}
    \end{enumerate}
\end{definition}

The definition of the combinatorial Loewener property given above slightly differs from the one in \cite{BourK}. This version appeared first in \cite{clais}, although \ref{CLP2} was written there in a slightly different, but equivalent, form.

\section{Iterated graph systems}\label{Section:IGSs}
In this section, we introduce the core concepts of the paper, the iterated graph systems (IGSs) and replacement graphs. Similar notions have been considered in many works with different names, among them hierarchical graph systems, hierarchical graph products and self-similar graphs \cite{neroli2024fractal, Leeslash,XI2017ScaleFree,kron2002green,Qi2017ConsensusIS,godsil,chilakamarri2013self}.
The term ``iterated graph systems'' originates from the work of Neroli Z. \cite{neroli2024fractal}.

\subsection{Definitions and notations}

\begin{definition}\label{def: Iterated graph system}
An \emph{iterated graph system (IGS)} is a collection of different classes of data.
It consists of a connected graph $G$, which we refer to as the \emph{generator}, a finite non-empty set $\mathcal{T}$ of \emph{types}, a surjective \emph{typing function} $\mathfrak{t} : E \to \cT$ and, finally, each type $t \in \mathcal{T}$ is associated a \emph{gluing rule}.
For each $t \in \mathcal{T}$ the gluing rule consists of a pair of non-empty disjoint sets $I_{t,+},I_{t,-} \subseteq V(G)$ and bijective functions $\alpha_t : I_{t,+} \cup I_{t,-} \to I_{t,+} \cup I_{t,-}$ which satisfies $\alpha_t(I_{t,+}) = I_{t,-}$ and $\alpha^2_t = \id_{V(G)}$.
\end{definition}

When we specify the the data associated to a given IGS, we always write $\Xi := (G,\mathcal{T},\ft,\{ \alpha_t \}_{t \in \mathcal{T}})$. The sets $I_{t,\pm}$ are omitted for brevity but are always understood to be part of the data.
Whenever there is no danger of confusion, we will always implicitly associate the notation $G,\mathcal{T},\ft,\alpha_t$ to the given IGS without further mention.
We also implicitly always denote the vertex set $V := V(G)$ and the edge set $E := E(G)$.

\begin{remark}
    In \cite{anttila2024constructions}, we also discuss another type of IGSs. To distinguish the two, one can call the present ones \emph{vertex-IGSs} and the ones there \emph{edge-IGS}, according to whether the replacement procedure replaces vertices or edges.
    In this paper we only discuss vertex-IGSs, and thus simplify the terminology.
\end{remark}

We introduce some symbolic operations. Fix a finite non-empty set $V$; in our setting $V$ is always the vertex set of the generator $G = (V,E)$ of an IGS.

First symbolic operations are the \emph{projections}. For $n \geq m$ we denote $\pi_{n,m} : V^n \to V^{m}, \, v_1\ldots v_n \mapsto v_1\ldots v_m$. 
We abbreviate $\pi_{n} = \pi_{n,n-1}$.

The second are the \emph{shift functions}.
For $n \geq m$ we set $\sigma_{n,m} : V^n \to V^{n-m}$, $v_1\ldots v_n \mapsto v_{m+1}\ldots v_n$. If $m = n-1$ we abbreviate $\sigma_{n} = \sigma_{n,n-1}$.

Lastly, for $x = v_1\ldots v_n, \, y = w_1\ldots w_m$ and $n,m \in \N$ and $v_i,w_j \in V$, we define the \emph{product} $x \cdot y := v_1\ldots v_m \mapsto w_1\ldots w_n v_1\ldots v_m \in V^{n+m}$. In particular, $x = \pi_{n,m}(x) \cdot \sigma_{n,m}(x)$ for all $x \in V^n$ and $n \geq m$. Given a subset $L \subseteq V^m$ and $x \in V^n$ we denote $x \cdot L \subseteq V^{n+m}$ the image of $L$ under the map $y \mapsto x \cdot y,\, V^{m} \to V^{n+m}$. Note that $x \cdot V^m = \pi_{n+m,n}^{-1}(x)$.

\subsection{Replacement rule}
\label{subsec:RR}
For a given IGS and its data $\Xi := (G,\cT,\ft,\{ \alpha_t \}_{t \in \mathcal{T}})$, where $G = (V,E)$, we consider a recursive replacement rule which produces a sequence of graphs $\{G_n := (V^n,E_n)\}_{n \in \N}$.
We refer to them as \emph{replacement graphs}.
At each step in the recursion, we need a graph $G_n$ and and a typing function $\ft_n$. In the first step we simply set $G_1 := G$ and $\mathfrak{t}_1 := \mathfrak{t}$.
Next assume $G_{n-1}$ and $\mathfrak{t}_{n-1}$ have been defined, and let $x, y\in V^{n}$.
Then $(x,y) \in E_{n}$ if and only if one of the following two conditions holds.
\begin{enumerate}
    \vspace{2pt}
    \item $\pi_n(x) = \pi_n(y)$ and $(\sigma_{n}(x),\sigma_{n}(y)) \in E$.
    \vspace{2pt}
    \item $(\pi_{n}(x),\pi_{n}(y)) \in E_{n-1}$, $\sigma_n(x) \in I_{t,+}$ and $\alpha_t(\sigma_n(x)) = \sigma_n(y) \in I_{t,-}$ where $t = \ft_{n-1}(\pi_n(x),\pi_n(y))$.
\end{enumerate}
Finally, the new typing function is
\[
    \mathfrak{t}_{n}(x,y) := 
    \begin{cases}
        \mathfrak{t}(\sigma_{n}(x),\sigma_{n}(y)) & \text{ if } \pi_{n}(x) = \pi_{n}(y) \\
        \mathfrak{t}_{n-1}(\pi_{n}(x),\pi_{n}(y)) & \text{ if } (\pi_{n}(x),\pi_{n}(y)) \in E_{n-1}.
    \end{cases}
\]

Next, we review our model examples; we strongly encourage the reader to draw some pictures for themselves.
For simplicity, we abbreviate the notation regarding the gluing rules. We identify $\alpha_t$ as a subset of $V \times V$ such that $(v,w) \in \alpha_t$ if and only if $v \in I_{t,+}$ and $\alpha_t(v) = w \in I_{t,-}$.
Such identification is well-defined by the properties of the gluing rules in Definition \ref{def: Iterated graph system}.
We use this convention only to simplify the description of explicit examples.

\begin{example}[\textbf{Sierpi\'nski carpet}]
Let $G = (V,E)$ be defined with vertices $V=\{0,\dots, 7\}$ and edges 
\begin{align*}
    E=\{&(0,1),(1,2),(5,6),(6,7),\\
    & (0,3),(2,4),(3,5),(4,7)\}.
\end{align*}
The first four edges which lie on the first line have type $h$ (horizontal), and the four remaining edges on the second row have type $v$ (vertical). The gluing rules are given by the functions 
\begin{align*}
    \alpha_h&=\{(2,0),(4,3),(7,5)\}\\
    \alpha_v&=\{(5,0),(6,1),(7,2)\}.
\end{align*}
These graphs produce a sequence of (discrete) metric spaces where the metric is the path metric.
By taking a (subsequential) Gromov-Hausdorff limit of $(V^n,3^{-n}d_{G_n})$, we obtain a biLipschitz equivalent version of the standard Sierpi\'snki carpet.
The first three replacement graphs can be found in Figure \ref{fig:SC}.
\end{example}

\begin{example}[\textbf{Menger sponge}] Let $G=(V,E)$ be defined with vertices $V=\{0,\dots, 19\}$ and edges 
\begin{align*}
    E=\{&(0,1),(1,2),(5,6),(6,7),(12,13),(13,14),(17,18),(18,19)\\
    & (0,3),(2,4),(3,5),(4,7), (12,15),(14,16),(15,17),(16,19) \\
    & (0,8),(2,9),(5,10),(7,11),(8,12),(9,14),(10,17),(11,19)\}.
\end{align*}
The first eight edges which lie on the first line have type $x$ ($x$-axis), second eight edges on the second row have type $y$ ($y$-axis) and the final eight on the last row have type $z$ ($z$-axis). The gluing rules are given by
\begin{align*}
    \alpha_x &=\{(2,0),(4,3),(7,5),(9,8),(11,10),(14,12),
    (16,15),(19,17)\}\\
    \alpha_y &=\{(5,0),(6,1),(7,2),(10,8),(11,9),(17,12),(18,13),(19,14)\}\\
    \alpha_z &=\{(12,0),(13,1),(14,2),(15,3),(16,4),(17,5),(18,6),(19,7)\}.
\end{align*}
This construction produces a limit space which is biLipschitz equivalent to the usual Menger sponge.
\end{example}

By mimicking the construction in the previous two examples, we can construct an IGS producing higher dimensional Menger sponges and also the  \emph{generalized Sierpi\'nski carpet}; see \cite[Subsection 2.2]{UniquenessofBrownianMotionOnGSC} for the definition.

\begin{example}
[\textbf{Pillow space}]
    Let $G = (V,E)$ be defined  with the vertices
    \[
        V := \{ 0,1,2,3,4a,4b,5,6,7,8 \}
    \]
    and edges
    \begin{align*}
        E=\{&(0,1),(1,2),(3,4a),(3,4b),(4a,5),(4b,5),(6,7),(7,8)\\
    & (0,3),(3,6),(1,4a),(1,4b),(4a,7),(4b,7),(2,5),(5,8)\}.
    \end{align*}
    So $G$ is a $3 \times 3$ square graph where the central vertex is doubled. The edges in the first row have the type $h$ (horizontal) and on the second are of the type $v$ (vertical). The gluing rules are given by
    \begin{align*}
    \alpha_h&=\{(2,0),(5,3),(8,6)\}\\
    \alpha_v&=\{(6,0),(7,1),(8,2)\}.
\end{align*}
First two replacement graphs can be found in Figure \ref{fig:DS}.
\end{example}

Next, we review few simple lemmas. We begin by addressing the obvious self-similarity of our construction. We recall some standard terminology first.
Given two general graphs $G $ and $ G'$ we say that $f : V(G) \to V(G)$ is \emph{simplicial}, and write $f : G \to G'$, if for all $\{x,y\} \in E(G)$ we have $\{f(x),f(y)\} \in E(G')$ or $f(x) = f(y)$. If $f$ is also a bijective function and $f^{-1} : G' \to G$ is simplicial, we say that $f : G \to G'$ is a \emph{graph isomorphism}.
Given a subset $U \subseteq V(G)$ the corresponding \emph{induced subgraph} is the graph $\Span{U}$ satisfying $V(\Span{U}) = U$ and $E(\Span{U}) = E(G) \cap (U \times U)$.

\begin{lemma}\label{lemma:CharEdges}
    Let $n \in \N, $ and $ x = v_1\ldots v_n, \, y = w_1\ldots w_n \in V^n$ be distinct. Fix $1 \leq l \leq n$
    such that $\pi_{n,l}(x) \neq \pi_{n,l}(y)$. Then the following are equivalent.
    \begin{enumerate}
        \vspace{2pt}
        \item $(x,y) \in E_n$ edge and is of type $t \in \mathcal{T}$.
        \vspace{2pt}
        \item $(\pi_{n,l}(x),\pi_{n,l}(y))$ is an edge of type $t$, $w_k \in I_{t,+}$ and $\alpha_t(w_{k}) = v_k$ for all $l < k \leq n$.
    \end{enumerate}
    Moreover, the replacement graphs are self-similar in the following sense. For all $n,m \in \N$ and $x \in V^m$ the mapping $V^n \to V^{n+m}$ given by the product $y \mapsto x \cdot y$ is a graph isomorphism onto the subgraph induced by $\pi_{n+m,m}^{-1}(x)$, and the inverse is given by the shift function $\sigma_{n+m,n}$.
\end{lemma}

\begin{proof}
    These follow from an induction on the steps of the replacement rule.
\end{proof}

The following lemma is needed to ensure that the replacement graphs do not grow uncontrollably.

\begin{lemma}\label{lemma: characterize bounded degree}
    Assume that the following conditions hold:
    \begin{enumerate}
    \vspace{2pt}
        \item For all $t \in \mathcal{T}$, if $v \in I_{t,+}$ then there is no edge $(v,w) \in E(G)$ of type $t$.
        \vspace{2pt}
        \item For all $t \in \mathcal{T}$, if $w \in I_{t,-}$ then there is no edge $(v,w) \in E(G)$ of type $t$.
    \end{enumerate}
    Then the replacement graphs have uniformly bounded degree, meaning
    \[
        \sup_{n \in \N} \deg(G_n)  < \infty.
    \]
\end{lemma}

\begin{proof}
    For any $n \in \N$ and $v \in V^n$ define
    \[
        \deg_{t,+}(x) := \abs{\{ y \in V^{n} : (x,y) \in E_{n} \text{ and } \ft(x,y) = t \}}.
    \]
    We also define $\deg_{t,-}(x)$ in an analogous way.
        
    By the iterative construction of $G_n$, for $n \geq 2$ we have the bound
    \begin{equation*}
        \deg_{t,+}(x) \leq \deg_{t,+}(\sigma_n(x)) + \deg_{t,+}(\pi_n(x))
    \end{equation*}
    
    Fix a vertex $x \in V^{n}$ for $n \geq 2$. If $\sigma_{n}(x) \in I_{t,+}$, it follows from the previous inequality and the first assumption in the lemma,
    \[
        \deg_{t,+}(x) \leq \deg_{t,+}(\sigma_n(x)) + \deg_{t,+}(\pi_n(x)) = \deg_{t,+}(\pi_n(x)).
    \]
    If $\sigma_{n}(x) \notin I_{t,+}$, then it follows from the recursive construction of $G_n$,
    \begin{align*}
        \deg_{t,+}(x) & = \deg_{t,+}(\sigma_{n}(x)) \leq \deg(G).
    \end{align*}
    So in conclusion, $\deg_{t,+}(x) \leq \max\{ \deg_{t,+}(\pi_n(x)),\deg(G) \}$. After applying this inequality $n-1$ times, we get $\deg_{t,+}(x) \leq \deg(G)$.

    By an identical argument and using the properties of the gluing rules $\alpha_{t}$ in Definition \ref{def: Iterated graph system}, we see that $\deg_{t,-}(x) \leq \deg(G)$.
    Since this holds for all $t \in \mathcal{T}$ and $n \in \N$, we can conclude the proof.
\end{proof}

\begin{remark}
    The converse to Lemma \ref{lemma: characterize bounded degree} is also valid but we leave the details to an interested reader.
\end{remark}

\subsection{Symmetric IGSs}\label{subsec:Symmetries}
Our methods in later sections heavily rely on suitable symmetries.
Here we describe what kind of conditions we are looking for.

\begin{definition}\label{def:symmetry}
Let $\Xi$ be an IGS and $G$ be its generator. We say that $\alpha : V \to V$ is a \emph{symmetry} of $\Xi$ if the map $v_1\ldots v_n \mapsto \alpha(v_1),\ldots,\alpha(v_n)$ defines a graph isomorphism $G_n \to G_n$ for all $n \in \N$. For simplicity, we usually denote higher level maps by $\alpha$ as well.
\end{definition}

Symmetries of IGSs have an easily verifiable characterization, but it is given by a somewhat technical symbolic condition. We nevertheless state it for the convenience of the reader.

\begin{lemma}\label{lemma:CharSymmetry}
    Let $\Xi$ be an IGS and $G$ be its generator. Then a graph isomorphism $\alpha : G \to G$ is a symmetry of $\Xi$ if and only if the following conditions hold.
    \begin{enumerate}
        \vspace{2pt}
        \item If $(v,w) \in E(G)$ is an edge of type $t$ and $(\alpha(v),\alpha(w)) \in E(G)$ is an edge of type $t'$ then $\alpha(I_{t,+}) = I_{t',+}$, $\alpha(I_{t,-}) = I_{t',-}$ and
        \[
        (\alpha_{t'} \circ \alpha)(v) = (\alpha \circ \alpha_{t})(v) \text{ for all } v \in I_{t,+} \cup I_{t,-}.
        \]
        \item If $(v,w) \in E(G)$ is an edge of type $t$ and $(\alpha(w),\alpha(v)) \in E(G)$ is an edge of type $t'$ then $\alpha(I_{t,+}) = I_{t',-}$, $\alpha(I_{t,-}) = I_{t',+}$
        \[
        (\alpha_{t'} \circ \alpha)(v) = (\alpha \circ \alpha_{t})(v) \text{ for all } v \in I_{t,+} \cup I_{t,-}.
        \]
    \end{enumerate}
\end{lemma}

\begin{proof}
    The ``if'' part follows from Lemma \ref{lemma:CharEdges}, and the ``only if'' part follows from the fact that $G_2 \to G_2$, $v_1v_2 \mapsto \alpha(v_1)\alpha(v_2)$ is a graph isomorphism.
\end{proof}

\begin{definition}\label{def:SymmetricIGS}
    We say that the gluing rules of an IGS $\Xi$ are \emph{symmetric} if for every $t \in \mathcal{T}$ there is a symmetry $\beta_t : V \to V$ of $\Xi$ with
    \[
        \text{$\beta_t|_{I_{t,+} \cup I_{t,-}} = \alpha_t$} \text{ and } \beta_t^2 = \id_{V}.
    \]
    Whenever the gluing rules are symmetric, we implicitly assume that the gluing rules $\alpha_t = \beta_t$ are restrictions of symmetries of $\Xi$ to the set $I_{t,+} \cup I_{t,-}$.
    For simplicity we also denote the higher level map $v_1\ldots v_n \mapsto \beta_t(v_1)\ldots \beta_t(v_n)$ by $\alpha_t$.
\end{definition}

Symmetric gluing rules allow us to perform \emph{foldings} which simplify moduli computations in later sections. These ideas were inspired by Bourdon--Kleiner \cite{BourK} but similar techniques also appear in e.g. \cite{murugan2023first,kigami}.
Below, we describe an algorithm which takes $n,m \in \N$ and a path $\theta=[x_1,\dots, x_k]$ in $G_{n+m}$ as input and produces as output the \emph{folding} $\fold_m(\theta)$ of $\theta$, which is also a path in $G_{n+m}$. The folded path has the additional property, that it lies entirely in $y\cdot V^m = \pi_{n+m,n}^{-1}(y)$ for some $y\in V^n$.

Consider the sequence $y_1,\ldots y_k \in V^n$ where $y_i = \pi_{n+m,n}(x_i)$. If $y_{i_1},\ldots y_{i_l}$ is the maximal subsequence of $y_1,\ldots y_k$ with no consecutive tuples, then our algorithm has $l$ iteration steps in total.

First, if $y_i$ is a constant sequence, then we return $\theta = \fold_m(\theta)$.
If not, let $j$ be the least index for which $y_{j} \neq y_{j+1}$. 
Then, by definition of the projections, $x_i = y_j \cdot \tilde{x}_i$ for some $\tilde{x}_i \in V^m$, for all $1 \leq i \leq j$.
Because the shift maps are isomorphisms (onto their image) according to Lemma \ref{lemma:CharEdges}, the sequence $[z_1,\ldots z_j]$ where $z_i = y_{j+1}\cdot \alpha_t(\tilde{x}_i)$ is a path in $G_{n+m}$.
Here $t$ is the type of $\{y_j,y_{j+1}\} \in E_n$.
By Lemma \ref{lemma:CharEdges}, $z_j = x_{j+1}$.
We can define the new path
\[
    \tilde{\theta} := [z_1,\ldots,z_j,x_{j+2},\ldots x_k].
\]
Next, we replace $\theta$ with $\tilde{\theta}$ and repeat the previous procedure. 
Note that $\tilde{\theta}$ is shorter than the path $\theta$ that we started with so this algorithm will always terminate. Moreover, $\tilde{\theta} \subseteq \pi_{n+m,n}^{-1}(\pi_{n+m,n}(x_k))$ where $x_k$ is the last vertex in the original path $\theta$.

The following somewhat technical lemma will be very useful later.
Roughly speaking it states: If the folded path is small, then the original path is also small.
In the statement, we use the same notation as above and also define $\mathcal{L}_m$ to be the collection consisting of the sets $I_{t,+}^m,I_{t,-}^m \subseteq V^m$.
Recall also the notation $x \cdot L$ from Subsection \ref{subsec:RR}.

\begin{lemma}\label{lemma:foldingFar}
    Define the collection
    \[
        \mathcal{C} := \{ \pi_{n+m,n}(x_k) \cdot L \subseteq V^{n+m} : L \in \mathcal{L}_m \text{ and } \fold_m(\theta) \cap (\pi_{n+m,n}(x_k) \cdot L) \neq \emptyset \}.
    \]
    If all sets in $\mathcal{C}$ have a common vertex, then there is a path $[h_1,\ldots, h_l]$ in $G_{n+m}$ such that $\pi_{n+m}(h_1) = y_1,\, \pi_{n+m}(h_k) = y_k$ and $\pi_{n+m,n}(h_i) \neq \pi_{n+m,n}(h_{i+1})$ for all $i = 1,\ldots,l-1$.
\end{lemma}

\begin{proof}
    Let $z_*$ be a common vertex of all members in $\mathcal{C}$.
    We will prove the claim by backtracking the construction of $\fold_m(\theta)$.
    
    Suppose that $y_{i_1},\ldots y_{i_l}$ is the maximal subsequence of $y_1,\ldots y_k$ with no consecutive tuples. By construction, the folding procedure had $l$ iteration steps in total.
    Also, $\{ y_{i_j},y_{i_{j+1}} \} \in E_n$ for all $j$, and
    let $t_j$ be its type.
    We claim that $h_j$ can be given inductively as $h_k = z_*$ and $h_{j-1} = y_{i_{j-1}} \cdot \alpha_{t_j}(\sigma_{n+m,n}(h_j))$.

    To see why, first note that $\{h_k,h_{k-1}\} \in E_{n+m}$ by construction and the definition of symmetric gluing rules.
    Let $\theta'$ be the path in the beginning of the last round in the folding. Then $\theta'$ contains a vertex in $y_{i_{l-1}} \cdot I_{t_{l-2},\pm}^m$ where the sign depends on the orientation of $\{ y_{i_{l-2}},y_{i_{l-1}} \}$.
    For simplicity, let the sign be $+$.
    Note that $\alpha_{t_{l}}(I_{t_{l-2},+}^m) \in \mathcal{L}_m$ by Lemma \ref{lemma:CharSymmetry}. After the last folding step, $x_k \cdot \alpha_{t_{l}}(I_{t_{l-2},+}^m) \in \mathcal{C}$, which means that $z_* \in \alpha_{t_{l}}(I_{t_{l-2},+}^m)$. Because $\alpha_{t_j}^2$ is the identity by Definition \ref{def:SymmetricIGS}, we must have $\alpha_{t_j}(z_*) \in I_{t_{l-2},+}^m$.
    So we conclude that $\{ h_{k-2},h_{k-1} \} \in E_{n+m}$. Then we repeat this backtracking process until we reach the beginning of the path.
\end{proof}

Symmetries are also used to ``bend'' flows, and the following definition captures when this is possible. It is applicable for general graphs $G$.

\begin{definition}\label{def:Reflection}
    Let $G$ be a general graph and $\alpha$ be a graph isomorphism $G \to G$. We say that $\alpha$ is a \emph{reflectional symmetry of $G$} if there is a partition $V(G) = C \cup D \cup \Delta$, where $\Delta := \{ v \in V(G) : \alpha(v) = v \}$, so that $\alpha(C) = D$, $\alpha^2 = \id_{V}$ and there is no edge between vertices in $C$ and $D$.
    Here, partition refers to a pair-wise disjoint covering.
\end{definition}

The model example of reflectional symmetries are the diagonal reflections of hypercubes, e.g. the reflection over the line $y = x$ in $\R^2$. Later, we will also apply reflections to the pentagonal carpet. 
Their main application is to bend flows, which we state in following proposition; similar arguments were used by Barlow--Bass \cite{BarlowBassResistanceOnSC} and Kwapisz \cite{Kwapisz}.

\begin{proposition}\label{prop:TwistFlow}
    Let $G$ be a graph and $S,T \subseteq V(G)$ be disjoint non-empty subsets. Let $\alpha$ be a reflectional symmetry of $G$ and with the associated partition $V(G) = C \cup D \cup \Delta$. Lastly, we assume that $S \subseteq C, \,  T \subseteq D$ and that $S \cap \alpha(T) = \emptyset = \alpha(S) \cap T$. Then given a unit flow $\mathcal{J}$ from $S$ to $T$,
    \[
        \widetilde{\mathcal{J}}(x,y) :=
        \begin{cases}
            \mathcal{J}(x,y) & \text{ if } x,y \in \Delta\\
            \mathcal{J}(x,y) + \mathcal{J}(\alpha(x),\alpha(y)) & \text{ if } \{x,y\} \cap C \neq \emptyset\\
            0 & \text{ otherwise}
        \end{cases}
    \]
    is a unit flow from $S$ to $\alpha(T)$.
\end{proposition}

\begin{proof}
    We begin by checking the divergences at vertices in $x \notin S \cup \alpha(T)$. Fix such a vertex $x$.
    First, assume $x \in C$. We note that $x,\alpha(x) \notin S \cup T$. This is because $x \notin S$ is already assumed and $x \notin T$ follows from $T \subseteq D$.
    Also, $\alpha(x) \notin S$ follows from $S \subseteq C$, and $\alpha(x) \notin T$ from $x \notin \alpha(T)$ and $\alpha^2 = \id_V$.
    Now, since there is no edge between $C$ and $D$,
    \[
        \divr(\widetilde{\mathcal{J}})(x) = \divr(\mathcal{J})(x) + \divr(\mathcal{J})(\alpha(x)) = 0 + 0 = 0.
    \]
    Next, assume $x \in \Delta$. By the definition of $\widetilde{\mathcal{J}}$,
    \begin{align*}
    \divr(\widetilde{\mathcal{J}})(x) & = \sum_{ \substack{ \{x,y\} \in E\\ y \in \Delta }} \widetilde{\mathcal{J}}(x,y) + \sum_{ \substack{ \{x,y\} \in E\\ y \in C }} \widetilde{\mathcal{J}}(x,y) + \sum_{ \substack{ \{x,y\} \in E\\ y \in D }} \widetilde{\mathcal{J}}(x,y)\\
    & = \sum_{ \substack{ \{x,y\} \in E\\ y \in \Delta }} \mathcal{J}(x,y) + \sum_{ \substack{ \{x,y\} \in E\\ y \in C }} (\mathcal{J}(x,y) + \mathcal{J}(x,\alpha(y)))\\
    & = \sum_{ \substack{ \{x,y\} \in E\\ y \in \Delta }} \mathcal{J}(x,y) + \sum_{ \substack{ \{x,y\} \in E\\ y \in C }} \mathcal{J}(x,y) + \sum_{ \substack{ \{x,y\} \in E\\ y \in D }} \mathcal{J}(x,y)\\
    & = \divr(\mathcal{J})(v) = 0.
    \end{align*}
    In the second to last row, we used the equality
    \[
        \sum_{ \substack{ \{x,y\} \in E\\ y \in C }} \mathcal{J}(x,\alpha(y)) = \sum_{ \substack{ \{x,y\} \in E\\ y \in D}} \mathcal{J}(x,y)
    \]
    which follows from the properties of $\alpha$.
    Lastly, assume $x \in D$. Since there is an edge between $x$ and $C$, it trivially follows from the definition of $\mathcal{J}$ that $\divr(\mathcal{J})(x) = 0$.

    Next, we verify that $\widetilde{\mathcal{J}}$ is a unit flow. By the assumptions, $\alpha(x) \notin (S \cup T)$ if $x \in S$. Therefore
    \[
        \divr(\widetilde{\mathcal{J}})(x) = \divr(\mathcal{J})(x) + \divr(\mathcal{J})(\alpha(x)) = \divr(\mathcal{J})(x),
    \]
    and we get
    \[
    I(\widetilde{\mathcal{J}}) = \sum_{x \in S} \divr(\mathcal{J})(v) = 1.
    \]
\end{proof}

In some cases of interest, such as the Sierpi\'nski gasket, we need to consider symmetries that do not satisfy the assumptions above. Such symmetries enable us to reflect flows.

\begin{proposition}\label{prop:TwistFlow2}
    Let $G$ be a graph and $\alpha$ a graph isomorphism $G \to G$ such that there is a partition $V(G) = C \cup D \cup \Delta$ where $\Delta = \{ v \in V(G) : \alpha(v) = v \}$, $\alpha(C) = D$ and $\alpha^2 = \id_{V(G)}$. We also assume that there is no vertex between two vertices in $\Delta$ and that whenever $\{ x,y \} \in E(G)$ for $x\in C$ and $y \in D$ then $\alpha(x) = y$. Now let $\mathcal{J}$ be a unit flow from $S$ to $T$ such that $S \subseteq C$ and $T \subseteq D$. Then
    \[
        \widetilde{\mathcal{J}}(x,y) = \begin{cases}
            \mathcal{J}(x,y) & \text{ if } \{x,y\} \cap C \neq \emptyset\\
            -\mathcal{J}(\alpha(x),\alpha(y)) & \text{ otherwise.}
        \end{cases}
    \]
    Then $\widetilde{\mathcal{J}}$ is a unit flow from $S$ to $\alpha(S)$.
\end{proposition}

\begin{proof}
    We check the divergences at $x \notin S \cup \alpha(S)$.
    If $x \in C$, then $\divr(\widetilde{\mathcal{J}})(x) = \divr(\mathcal{J})(x) = 0$ by the definition of $\widetilde{\mathcal{J}}$.
    Next, let $x \in \Delta$. Because there are no edges between vertices in $\Delta$, we have
    \begin{align*}
        \divr(\widetilde{\mathcal{J}})(x) & = \sum_{\substack{ \{x,y\} \in E(G) \\ y \in C }} \widetilde{\mathcal{J}}(x,y) + \sum_{\substack{ \{x,y\} \in E(G) \\ y \in D }}\widetilde{\mathcal{J}}(x,y)
        \\
        & = \sum_{\substack{ \{x,y\} \in E(G) \\ y \in C }} \mathcal{J}(x,y) - \sum_{\substack{ \{x,y\} \in E(G) \\ y \in D }}\mathcal{J}(x,\alpha(y))\\
        & = \sum_{\substack{ \{x,y\} \in E(G) \\ y \in C }} \mathcal{J}(x,y)-\sum_{\substack{ \{x,y\} \in E(G) \\ y \in C }} \mathcal{J}(x,y) = 0.
    \end{align*}
    The penultimate equality holds since $\alpha$ is a graph isomorphism with  $\alpha(C) = D$ and $\alpha^2 = \id_{V}$.
    Lastly, assume $x \in D$.
    Note that $\alpha(x) \notin S$ because we assumed $x \notin \alpha(S)$ and $\alpha^2 = \id_{V(G)}$.
    Also $\alpha(x) \notin T$ because $\alpha(x) \in C$ and $T \subseteq D$.
    Now, if the neighbors of $x$ are contained in $D \cup \Delta$, $\divr(\mathcal{J})(x) = -\divr(\mathcal{J})(\alpha(x)) = 0$ by the definition of $\widetilde{\mathcal{J}}$.
    If this is not the case, by the assumptions on $\alpha$, the only neighbor of $x$ that is not in $D \cup \Delta$ is $\alpha(x)$.
    We then obtain
    \begin{align*}
        \divr(\widetilde{\mathcal{J}})(x) &  = \mathcal{J}(x,\alpha(x)) - \sum_{ \substack{ \{ x,y \} \in E \\ y \in D \cup \Delta } } \mathcal{J}( \alpha(x),\alpha(y) )\\
        & = -\mathcal{J}(\alpha(x),x) - \sum_{ \substack{ \{ \alpha(x),y \} \in E \\ y \in C \cup \Delta } } \mathcal{J}( \alpha(x),y) = -\divr(\mathcal{J})(\alpha(x)) = 0.
    \end{align*}
    Finally, to see that $\widetilde{\mathcal{J}}$ is unit, this follows from $S \subseteq C$ and the definition of $\widetilde{\mathcal{J}}$.
\end{proof}

Finally, we end with a simple property for symmetries. 

\begin{proposition}\label{prop:ReflFold}
    Let $G, \alpha$ and $V(G) = C \cup D \cup \Delta$ be as in Proposition \ref{prop:TwistFlow2}. Then the function $f : V(G) \to V(G)$ given by
    \[
        f(x) :=
        \begin{cases}
            x & \text{ if } x \in C \cup \Delta\\
            \alpha(x) & \text{ if } x \in D
        \end{cases}
    \]
    is simplicial.
\end{proposition}

\begin{proof}
    This follows from the properties of $\alpha$.
\end{proof}

\section{Main Examples}\label{Section: Main examples}
This section introduces the main examples of IGSs that this work considers.
These are divided into two different classes; hypercube based constructions and regular polygon based constructions. The former generalizes the Sierpi\'nski carpet and higher dimensional Menger sponges, and also includes the pillow space.
For the latter class we only consider two specific examples, the pentagonal carpet and Sierpi\'nski gasket. We included a short discussion on further generalization in the end of the section.

\subsection{Hypercube based construction}\label{subsection: Cubical iterated graph systems}

Fix integer parameters $d \geq 2, L_* \geq 3$ and $K \geq 1$. 
Based on these choices, we first define a single large IGS, and the hypercube based IGSs are obtained as suitable ``sub-objects'' of these larger ones.
First, the set of types is $\mathcal{T} := \{ 1,\ldots, d \}$.
To define the underlying graph $G(d,L_*,K)$, we define the vertex set
\[
    V(d,L_*,K) := V(G(d,L_*,K)) := \{1,2,\dots, L_*\}^{d} \times \{1,\ldots K\}.
\]
For the sake of clearer notation, given $v = \left((z_i)_{i\in \mathcal{T}},s\right) \in V$, we write $c_i(v)= z_i$ for $i \in \mathcal{T}$ and $\xi(v) = s$.
The set of edges are defined so that
$(w,v) \in E(d,L_*,K) := E(G(d,L_*,K))$ is an edge of type $i \in \mathcal{T}$ if and only if
    \[
        c_j(v) =
        \begin{cases}
            c_j(w) & \text{ if } i \neq j\\
            c_j(w) + 1 & \text{ if } i = j.
        \end{cases}
    \]
    The typing function is denoted $\tilde{\ft}$.
    To define the gluing rules $\tilde{\alpha}_{i} : \tilde{I}_{i,+} \cup \tilde{I}_{i,-} \to \tilde{I}_{i,+} \cup \tilde{I}_{i,-}$ for $i \in \mathcal{T}$, we first fix the sets. We define $w \in I_{i,+}$ if and only if $c_i(w) = L_*$, and $w \in I_{i,-}$ if and only if $c_i(w) = 1$. Finally, we define $\tilde{\alpha}_{i}(w) = v$ if and only if
    \[
    \begin{cases}
        c_j(w) = c_j(v) & \text{ for all } i \neq j\\
        (c_i(w),c_i(v)) = (L_*,1)  & \text{ and }\\
        \xi(w) = \xi(v).
    \end{cases}
    \]
We denote the IGS $\Xi(d,L_*,K) := (G(d,L_*,K),\mathcal{T}(d),\tilde{\ft},\{\tilde{\alpha}_t\}_{\mathcal{T}(d)})$.

\begin{definition}
    We say that $\Xi = (G,\mathcal{T},\ft,\{\alpha_t\}_{t \in \mathcal{T}})$ is a \emph{hypercubic} IGS if there are $d \geq 2$, $L_* \geq 3$ and $K \geq 1$ such that $\Xi$ is a sub-ohject of the IGS described above in the following sense.
    \begin{enumerate}
    \vspace{2pt}
        \item $G$ is a subgraph of $G(d,L_*,K)$, meaning $V(G) \subseteq V(d,L_*,K)$ and the inclusion is simplicial.
    \vspace{2pt}
        \item $\mathcal{T} = \mathcal{T}(d)$ and $\ft = \tilde{\ft}|_{E(G)}$.
    \vspace{2pt}
        \item The gluing rules $\alpha_{i} : I_{i,+} \cup I_{i,-} \to  I_{i,+} \cup I_{i,-}$ satisfy $I_{i,\pm} \subseteq \tilde{I}_{i,\pm}$ and $\alpha_{i} = (\tilde{\alpha}_{i})|_{I_{i,+}}$ for all $i \in \mathcal{T}$.
    \end{enumerate}
    Unless specified otherwise, we always associate the notation $d,\,L_*,\,K$ to the parameters of a hypercubic IGS.
\end{definition}

\begin{remark}
    Replacement graphs of hypercubic IGSs with parameters $d,L_*$ and $K = 1$ have simplicial embeddings to the integer lattice graphs $\Z^d$.
\end{remark}

We are not yet done with the definition because we need suitable (dihedral) symmetries to ensure good analytic properties. We define them on $\Xi(d,L_*,K)$ first.
We begin by extending the gluing rules of $\Xi(d,L_*,K)$, namely $\tilde{\alpha}_{i}$ for $i \in \mathcal{T}$, to a symmetry of $\Xi(d,L_*,K)$ in the sense of Definition \ref{def:symmetry}.
For $i \in \mathcal{T}(d)$ we set $\tilde{\alpha}_{i}(w) = v$ where $\xi(w) = \xi(v)$ and
\[
    c_j(v) = 
    \begin{cases}
        c_j(w) & \text{ if } i \neq j\\
        L_* - c_i(w) + 1 & \text{ if } i = j.
    \end{cases}
\]
We call $\tilde{\alpha}_{i}$ a \emph{horizontal symmetry}.

Next for $1 \leq i < j \leq d$ we define the \emph{diagonal reflection} $\alpha_{i,j}^+ = \alpha_{j,i}^+ = \alpha$ such that $\alpha(w) := v$ where $\xi(v) = \xi(w)$ and
\[
    c_k(v) =  
    \begin{cases}
        c_k(w) & \text{ if } k \notin \{ i,j \}\\
        c_j(w) & \text{ if } k = i\\
        c_i(w) & \text{ if } k = j.
    \end{cases}
\]
for all $k \in \mathcal{T}$. We also define the \emph{anti-diagonal reflection} $\alpha_{i,j}^- = \alpha_{j,i}^- = \alpha$ such that $\alpha(w) = v$ where
    $\xi(w) = \xi(v)$ and
    \[
        c_k(v) =
        \begin{cases}
            c_k(w) & \text{ if } k \notin \{ i,j \}\\
            L_* - c_j(w_1) + 1 & \text{ if } k = i\\
            L_* - c_i(w_1) + 1 & \text{ if } k = j.
        \end{cases}
    \]
    
Finally, the set of the graph isomorphisms of $G(d,L_*,K)$ we defined above, namely $\tilde{\alpha}_i$ and $\alpha_{i,j}^{\pm}$, is denoted $\mathcal{G}(d,L_*,K)$.

\begin{definition}\label{Def:SymmHypercubic}
    We say that $\Xi = (G,\mathcal{T},\ft,\{\alpha_t\}_{t \in \mathcal{T}})$ is a \emph{symmetric hypercubic} IGS if $\Xi$ is hypercubic IGS for some parameters $d,L_*,K$ and the following additional conditions are satisfied.
    \begin{enumerate}[label={\textup{(\textcolor{blue}{C\arabic*})}}, widest=a, leftmargin=*]
    \vspace{2pt}
        \item \label{Def:SymmHypercubic(1)} For all $j \in \mathcal{T}$ the vertex set $V(G)$ contains all vertices $w$ with
        \begin{equation*}
            \begin{cases}
                c_i(w) \in \{ 1, L_*\} \text{ for all } i \neq j\\
                \xi(w) = 1,
            \end{cases}
        \end{equation*}
        \vspace{2pt}
        \item \label{Def:SymmHypercubic(2)} For all $j \in \mathcal{T}$ the edge set $E(G)$ contains all edges $(v,w)$ with
        \[
            \begin{cases}
                c_i(v) = c_i(w) \in \{ 1,L_* \} \text{ for all } i \neq j \\
                c_j(v) = c_j(w) + 1\\
                \xi(v) = \xi(w) = 1
            \end{cases}
        \]
        \vspace{2pt}
        \item \label{Def:SymmHypercubic(3)} For all $j \in \mathcal{T}$ we have $w \in I_{j,+}$ and $v \in I_{j,-}$ where
        \[
            \begin{cases}
                c_i(w) = c_i(v) & \text{ for all } i \neq j\\
                (c_j(w_2),c_j(v_2)) = (L_*, 1)\\
                \xi(w) = \xi(v) = 1.
            \end{cases}
        \]
        \vspace{2pt}
        \item \label{Def:SymmHypercubic(4)} For every $\alpha \in \mathcal{G}(d,L_*,K)$ the restriction $\alpha|_{V(G)}$ is a symmetry of $\Xi$ in the sense of Definition \ref{def:symmetry}.
    \end{enumerate}
\end{definition}

\begin{remark}\label{rem:SymmHyper}
    Symmetrical hypercubic IGSs are symmetric also in the sense of Definition \ref{def:SymmetricIGS}. This follows by noting that the gluing rules can be extended to the horizontal symmetries defined above.
\end{remark}

\begin{remark}
    Symmetrical hypercubic IGSs generalize the three examples in Subsection \ref{subsec:RR}.
    For $K = 1$ we get the Sierpi\'nski carpet, Menger sponge and also the broader class of generalized Sierpi\'nski carpets (see \cite{UniquenessofBrownianMotionOnGSC}. The pillow space on the other hand, can be understood as a symmetric hypercubic IGSs for $K = 2$.
\end{remark}

\begin{figure}[h!]
    \centering\includegraphics[width=260pt]{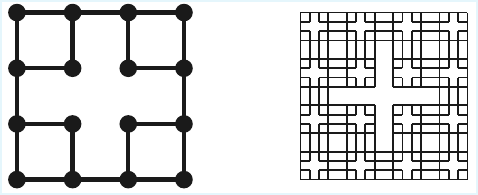} 
    \caption{Figure of a symmetric hypercubic IGS. The limit space can also be obtained by adding horizontal- and vertical ``slices'' to a unit square in a self-similar fashion.
    }
\end{figure}

\begin{figure}[!ht]
    \centering\includegraphics[width=260pt]{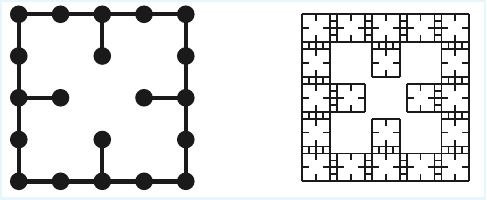} 
    \caption{In this figure we point out that our definition of hypercubic IGS does not allow direct diagonal connections. The construction in the figure is very similar to an Unconstrained Sierpi\'nski carpet in \cite{cao2024whether}. But it is nevertheless very different because the fractal space in the figure does not contain local cut-points.}
    \label{fig: non-gsc}
\end{figure}

\begin{figure}[!ht]
    \centering\includegraphics[width=230pt]{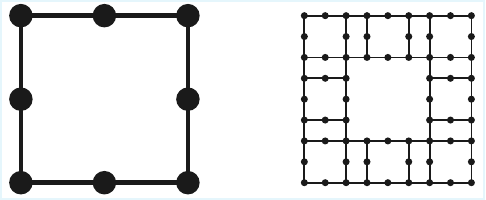} 
    \caption{A symmetric hypercubic IGS where the gluing takes place along a totally disconnected set.}
    \label{fig: Gluing along a disconnected set}
\end{figure}

\begin{figure}[!ht]
    \centering\includegraphics[width=230pt]{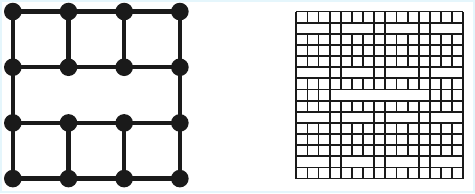} 
    \caption{A (non-symmetric) hypercubic IGS whose limit space is a variant of the slit carpet of Merenkov \cite{merenkov}.
    Because its conformal dimension is 2, by showing that the discrete 2-modulus from top to bottom decays to 0 with respect to the level of the graph, we could deduce that this space does not satisfy the combinatorial Loewner property.}
    \label{fig: Slit carpet}
\end{figure}

\subsection{Constructions by Regular Polygons}
We introduce the two other main examples of the work,
the \emph{Sierpi\'nski gasket} and the \emph{pentagonal carpet}.
The former is one of the model examples of fractal spaces in the mathematical literature, and it is also one of the most researched examples in analysis on fractals; see \cite{AnalOnFractals,barlow} for further references.
The latter seems to be a new construction, although see \ref{rmk:pentagonal tiling}.

We again identify the gluing rules $\alpha_t$ as a subset of $V \times V$ such that $(v,w) \in \alpha_t$ if and only if $v \in I_{t,+}$ and $\alpha_t(v) = w \in I_{t,-}$. We do not use this convention in later section.

\begin{example}[\textbf{Sierpi\'nski gasket}]\label{ex: SG}
    Let $G = (V,E)$ be the graph where $V = \{0,1,2\}$
    and $E = \{(0,1), (1,2), (0,2)\}$ with different types $a,b,c$, respectively.
    The gluing rules are given by
    \[
        \alpha_a = \{(1,0)\} \, , \alpha_b = \{ (2,1) \} \text{ and } \alpha_c = \{ (2,0) \}.
    \]
    We refer to this IGS as the \emph{Sierpi\'nski gasket} because the limit it produces is biLipschitz equivalent to the usual Sierpi\'nski gasket with the Euclidean metric.
    \begin{figure}[!ht]
    \centering\includegraphics[width=350pt]{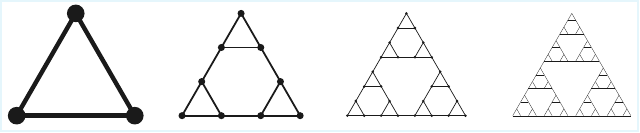} 
    \caption{First four replacement graphs of the Sierpi\'nski gasket.}
\end{figure}
\end{example}

\begin{remark}\label{rem:gasketSYmm}
    The Sierpi\'nski gasket is a symmetric IGS in the sense of Definition \ref{def:SymmetricIGS}. This is because the gluing rules $\alpha_a,\, \alpha_b, \, \alpha_c$ can be extended to symmetries of the equilateral triangle that yields graph isomorphisms $G_n \to G_n$ on the replacement graphs, according to Lemma \ref{lemma:CharSymmetry}.
    For instance, the non-trivial symmetry that fixes the vertex $i \in V$ can be written as
    \[
    \alpha(j) :=
        \begin{cases}
            i &\text{ if } j = i\\
            i + 1 &\text{ if } j = i - 1\\
            i - 1 &\text{ if } j = i + 1.\\
        \end{cases}
    \]
    These are symmetries of the IGS Sierpi\'nski gasket.
\end{remark}

\begin{example}[\textbf{Pentagonal carpet}]\label{ex: SP}
    Let $G = (V,E)$ be the graph $V = \{0,1,2,3,4\}$
    and $E = \{(0,1), (1,2), (2,3), (3,4), (4,0)\}$ with types $a,b,c,d,e$, respectively.
    The gluing rules are given by
    \begin{align*}
        \alpha_a & := \{ (1,0), (2,4) \},\, 
        \alpha_b := \{ (2,1), (3,0) \},\, 
        \alpha_c := \{ (3,2),(4,1) \}, \\
        \alpha_d & := \{ (4,3), (0,2) \}, \, 
        \alpha_e := \{ (0,4),(1,3) \}.
    \end{align*}
    We call this IGS the \emph{pentagonal carpet}.
    The first three replacement graphs can be found in Figure \ref{fig:PC}.
    We further introduce the notations
    \begin{align*}
        L_0 := \{ 2,3 \} = I_{b,+} = I_{d,-},\\
        L_1 := \{3,4\} = I_{c,+} = I_{e,-}, \\
        L_2 := \{ 0,4 \} = I_{a,-} = I_{d,+}, \\
        L_3 := \{ 0,1 \} = I_{b,-} = I_{e,+}, \\
        L_4 := \{ 1,2 \} = I_{a,+} = I_{c,-}.
    \end{align*}
    Notice that $L_i$ is the line segment opposite to the vertex $i \in V$. We heavily employ these notations in our arguments.
\end{example}

\begin{remark}
    The pentagonal carpet is also a symmetric IGS. To see this, it follows from Lemma \ref{lemma:CharSymmetry} that the symmetries of the regular pentagon are symmetries also in the sense of Definition \ref{def:SymmetricIGS}. For instance, the non-trivial symmetry of the regular pentagon that fixes the vertex $i \in V$ can be written as
    \[
        \alpha(j) :=
        \begin{cases}
            i &\text{ if } j = i\\
            i + 1 &\text{ if } j = i-1\\
            i - 1 &\text{ if } j = i+1\\
            i + 2 &\text{ if } j = i-2\\
            i - 2 &\text{ if } j = i+2.\\
        \end{cases}
    \]
    The arithmetic above is taken as modulo 5.
\end{remark}

The construction of a pentagonal Sierpi\'nski carpet is new, but it has a curious connection with a subdivision rule already studied in the literature.

\begin{remark}\label{rmk:pentagonal tiling} There exists a planar pentagonal subdivision rule \cite{bowers} which was first studied in \cite{CFP, CFKP}. This subdivision rule divides a pentagon into five pentagons recursively. This can be seen also from Figure \ref{fig:PC}, if one imagines a reflected pentagon being placed at the center of the construction. In \cite{CFKP} it is shown, that this subdivision rule is realized by a rational map, see also \cite{meyerquasi}. In \cite{CFP} and \cite{bowers} it is shown that the pentagonal subdivision rule can be realized as a conformal subdivision rule of the complex plane. See \cite{bowers} for visualizations of these conformal subdivision rules that are computed using circle packings. Such subdivision rules produce sequences of tilings of a pentagon. If one omits the central pentagon in the resulting tilings at every stage, one obtains a carpet. The remaining tiles of this carpet have the same combinatorial structure as in our graphs. It thus seems, that our pentagonal carpet is a subset of these subdivision rules equipped with an appropriately chosen visual metric in the sense of Bonk and Meyer \cite{BM}. Consequently, the carpets that we study may be interesting from a dynamical perspective, and this suggests that it may be useful to use the IGS framework in conjunction with the subdivision framework. This connection suggests many further problems, such as the extension of the rigidity work of Bonk and Merenkov; see \cite{BonkMerenkov}.
\end{remark}

\subsection{Further examples and discussions}
While our class of regular polygon based construction only consists of two explicit examples, we of course could consider many other examples. In the two figures below we included two IGSs that also fall within this class of constructions.

\begin{figure}[!ht]
    \centering\includegraphics[width=300pt]{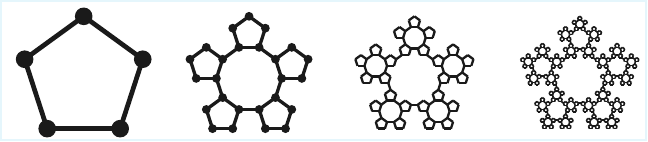} 
    \caption{Second pentagonal fractal space; see \cite{AnalOnFractals}.}
\end{figure}

\begin{figure}[!ht]
    \centering\includegraphics[width=300pt]{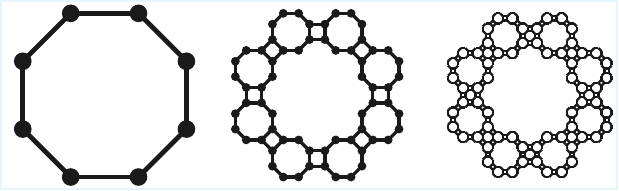} 
    \caption{An ocatagonal fractal space; see  \cite{Resistance4Ncarpets}.}
    \label{fig: octakun}
\end{figure}

The main reason for our decision to focus on only the two regular polygon based examples is that the analysis on these kinds of fractals seems to subtly depend on the particularities of each shape used. The cause of this seems to be the differences in the symmetry groups, and this complicates finding unified proofs for a larger class of examples.
Nevertheless, extensions to a broader class should be possible by adapting our methods. Later, we will also state a set of axioms on the replacement graphs, which seems to be quite general; see Section \ref{sec:LS}. If one is interested in an example that does not fall directly into one of the examples of this section, one can attempt to verify these axioms instead. These verifications should be quite similar to the various special cases that we discuss.

\section{Replacement flow}\label{Section: Replacement flow}

This section studies flows on replacement graphs.
Subsection \ref{subsec:ReplacementFLow} is devoted to prove the main result of the section, which regards the general replacement flow technique (Theorem \ref{thm: Replacement flow}). We state it for general IGSs with symmetric gluing rules and whose replacement graphs have uniformly bounded degrees.
This is the most important technique of the paper and it is a crucial ingredient in our proof of super-multiplicativity inequality and the combinatorial Loewner property in Section \ref{Section: CLP}.
In Subsection \ref{subsec:FB} we construct good flow bases for our main examples.
 
\subsection{Flow bases and replacement flow}\label{subsec:ReplacementFLow}

Throughout the subsection, we consider a fixed IGS $\Xi = (G,\mathcal{T},\ft,\{ \alpha_t \}_{t \in \mathcal{T}})$ such that its replacement graphs $G_n = (V^n,E_n)$ have uniformly bounded degree, meaning
\[
    C_{\deg} := \sup_{n \in \N} \deg(G_n) < \infty.
\]
We assume $\Xi$ to be symmetric in the sense of Definition \ref{def:SymmetricIGS}. The notation $\alpha_t$ for $t \in \mathcal{T}$ denotes both the gluing rules and the graph isomorphisms $G_n \to G_n$. We also need to be careful with the orientations of the edges. If $G$ is any oriented graph and $\{ x,y \} \in E$ is any edge we write
\[
    \sgn(x,y) :=
    \begin{cases}
        + & \text{ if } (x,y) \in E\\
        - & \text{ if } (y,x) \in E.
    \end{cases}
\]

For $m \in \N$ we define $\mathcal{L}_m$ to be the collection consisting of the subsets $I_{t,+}^m, I_{t,-}^m \subseteq V^m$ for $t \in \mathcal{T}$.
For our flow method it will be useful to introduce the following \emph{extended replacement graphs}. 
For $m \in \N$, we set $G^e_m := (V_m^e, E_m^e)$ to be the graph that is obtained  by extending the set of vertices
\[
    V_m^e := V^m \cup \{ x_L : L \in \cL_m \text{ and } x \in L \}
\]
and by adding the edges
\[
    E_m^e := E_m \cup \left( \bigcup_{L \in \cL_m}  \{ \{ x_L, x \} : x \in L \} \right).
\]
Finally, for $L \in \cL_m$, we denote $L^e := \{ x_L : x \in L \}$. Note that $L^e_1 \cap L_2^e = \emptyset$ for all distinct pairs $L_1,L_2 \in \mathcal{L}_m$.
Also $\alpha_t(I_{t,\pm}^m) = I_{t,\mp}^{m}$ for all $t \in \mathcal{T}$.
In particular if $L = I_{t,-}^m$, then the notation $\alpha_t(x)_L$ is sensible. We use this in \ref{(FB4)} below.

\begin{definition}\label{def:FB}
A collection $\mathcal{B}$ of flows in $G_m^e$ is a \emph{flow basis} if it satisfies the following conditions.
    \begin{enumerate}[label={\textup{(\textcolor{blue}{FB\arabic*})}}, widest=a, leftmargin=*]
    \vskip.3cm
        \item \label{(FB1)} Each flow in $\mathcal{B}$ is a unit flow from $L_1^e$ to $L_2^e$ for some distinct pair $L_1, L_2 \in \mathcal{L}_m$.
        Furthermore, for each distinct pair $L_1,L_2 \in \mathcal{L}_m$ there is exactly one unit flow from $L_1^e$ to $L_2^e$ in
        $\mathcal{J}$ denoted as $\mathcal{J}_{[L_1 \to L_2]}$.
        \vskip.3cm
        \item \label{(FB3)} If $L \neq L_1,L_2$ and $x \in L$,
        \[
            \mathcal{J}_{[L \to L_1]}(x_{L},x) = \mathcal{J}_{[L \to L_2]}(x_{L},x) \text{ and } \mathcal{J}_{[L_1 \to L]}(x,x_{L}) = \mathcal{J}_{[L_2 \to L]}(x,x_{L}).
        \]
        \vskip.3cm
        \item \label{(FB4)} It holds for all $t \in \mathcal{T}, \, L_1 = I_{t,+}^{m},\, L_2 = I_{t,-}^{m}$ and $x \in L_1$,
        \[
            \mathcal{J}_{[L_1 \to L_2]}(x_{L_1},x) = \mathcal{J}_{[L_1 \to L_2]}(\alpha_t(x),\alpha_t(x)_{L_2})
        \]
        and
        \[
            \mathcal{J}_{[L_1 \to L_2]}(x_{L_1},x) = \mathcal{J}_{[L_2 \to L_1]}(x,x_{L_1}).
        \]
        Here $\alpha_t$ is the gluing rule/symmetry of type $t$.
    \end{enumerate}
    The \emph{$p$-energy} of $\mathcal{B}$ is $\cE_q(\mathcal{B}) := \max \left\{ \cE_q(\mathcal{J}) : \mathcal{J} \in \mathcal{B} \right\}$ where $q = p/(p-1)$.
\end{definition}

For $n \in \N$ we fix a source and a target $S_n,T_n \subseteq V^n$, which is a disjoint pair of non-empty subsets. For every $u \in S_n \cup T_n$ we assign any choice of a pair $(\ft(u),\sgn(u)) \in \mathcal{T} \times \{+,-\}$. Then we define the source and the target in $G_{n+m}$,
\[
    S_{n+m} := \{ u \cdot x \in V^{n+m} : u \in S_n \text{ and } x \in I_{\ft(u),\sgn(u)}^m \} \subseteq V^{n+m}.
\]
Similarily we define
\[
    T_{n+m} := \{ u \cdot x \in V^{n+m} : u \in T_n \text{ and } x \in I_{\ft(u),\sgn(u)}^m \} \subseteq V^{n+m}.
\]

We can now state the replacement flow theorem.

\begin{theorem}\label{thm: Replacement flow}
Fix $n,m \in \N$ and $p \in (1,\infty)$.
Let $\Xi$ be an IGS with symmetric gluing rules and such that $\deg(G_n) \leq C_{\deg} < \infty$ for all $n \in \N$.
Let $\mathcal{B}$ be a flow basis of $G_m^e$, and $S_n,T_n,S_{n+m},T_{n+m}$ be as above. Then there is a constant $C \geq 1$ depending only on $C_{\deg}$ and $p \in (1,\infty)$ s.t. the following holds. For every unit flow $\mathcal{J}_n$ from $S_n \to T_n$ in $G_n$ there exists a unit flow $\mathcal{J}_{n+m}$ from $S_{n + m}$ to $T_{n+m}$ in $G_{n + m}$ satisfying
\begin{equation*}
    \cE_q(\mathcal{J}_{n + m}) \leq C \cdot \cE_q(\mathcal{B})\cE_q(\mathcal{J}_n).
\end{equation*}
The constant $C$ can be chosen so that $C^{1-p}$ is continuous in $p \in (1,\infty)$ and has a limit when $p \to 1^+$.
\end{theorem}

Theorem \ref{thm: Replacement flow} is quite powerful when we have good estimates for $\mathcal{E}_q(\mathcal{B})$. An important application is the \emph{super-multiplicativity inequality}. For its proof in Section \ref{Section: CLP}, we need the following theorem.

\begin{theorem}\label{thm:ReplacementflowSuper}
    Let $\Xi$ be an IGS with the same assumptions as in Theorem \ref{thm: Replacement flow} and $p \in (1,\infty)$.
    Assume that there is a constant $C_1 = C_1(p)$ such that $C_1^{1-p}$ is continuous in $p \in (1,\infty)$, has a limit when $p \to 1^+$ and for all $m \in \N$ there is a flow basis $\mathcal{B}$ on $G_m^e$ satisfying
    \[
        \mathcal{E}_q(\mathcal{B}) \leq C_1 \cRes_p(I_{t,+}^m,I_{t,-}^m,G_n)
    \]
    for all $t \in \mathcal{T}$.
    Then there is a constant $C_2$ depending on $C_1,p$ and $C_{\deg}$ such that for all $n,m \in \N$ and $t \in \mathcal{T}$,
    \[
        \cRes_p(I_{t,+}^{n+m},I_{t,-}^{n+m},G_{n+m}) \leq C_2\cRes_p(I_{t,+}^{n},I_{t,-}^{n},G_n)\cdot \cRes_p(I_{t,+}^{m},I_{t,-}^{m},G_m)
    \]
    The constant $C_2$ can be chosen so that $C_2^{1-p}$ is continuous in $p \in (1,\infty)$ and has a limit when $p \to 1^+$.
\end{theorem}

\begin{proof}
    Let $\mathcal{J}_n$ be the minimizer for $\cRes_p(I_{t,+}^n,I_{t,-}^n,G_n)$. In particular, its source is $S_n := I_{t,+}^n$ and target is $T_n := I_{t,-}^n$. We choose the sets $S_{n+m}$ and $T_{n+m}$ as follows. For each $u \in S_n$ we choose $(\ft(u),\sgn(u)) = (t,+)$, and for $y \in T_n$ we choose $(\ft(u),\sgn(u)) = (t,-)$. Then $S_{n+m} = I_{t,+}^{n+m}$ and $T_{n+m} = I_{t,-}^{n+m}$.
    Then the super-multiplicative inequality, as in the claim, follows from the energy estimate on the flow basis and the replacement flow theorem ( Theorem \ref{thm: Replacement flow}).
\end{proof}

The proof of Theorem \ref{thm: Replacement flow} is quite technical due to the delicate nature of flows. We will introduce some additional notation to ease the computation.
We begin by performing an extension of $G_n$ and $G_{n+m}$ as well. Set $G_n^* := (V_n^*,E_n^*)$ where 
\[
    V_n^* := V^n \cup \{ u^* : u \in S_n \cup T_n \}
\]
and
\[
    E_n^* := E_n \cup \{ \{ u^*,u \} : u \in S_n \cup T_n \}.
\]
To define the extension $G_{n+m}^* := (V_{n+m}^*,E_{n+m}^*)$, for $u \in S_n \cup T_n$ we first set, $\ft(u^*,u) := \ft(u)$ and $\sgn(u^*,u) := \sgn(u)$. We also choose $\sgn(u,u^*)$ so that
\[
    \{ \sgn(u^*,u),\sgn(u,u^*) \} = \{+,-\}.
\]
Lastly we introduce new vertices $u^* \cdot \alpha_{\ft(u,u^*)}(x)$. This is only a formal notation; see the remark below.
We now define the extended graph so that the vertex set is
\[
    V_{n+m}^* := V^{n+m} \cup \left\{ u^* \cdot \alpha_{\ft(u,u^*)}(x) : u \in S_n \cup T_n \text{ and } x \in I_{\ft(u),\sgn(u)}^{m}  \right\}
\]
and the edge set
\[
    E^*_{n+m} := E_{n+m} \cup \left\{ \{ u \cdot x, u^* \cdot \alpha_{\ft(u,u^*)}(x) \} : u \in S_n \cup T_n \text{ and } x \in I_{\ft(u),\sgn(u)}^{m} \right\}.
\]

The orientation of these added edges will play no role, but $\ft(u^* \cdot \alpha_{\ft(u,u^*)}(x),u \cdot x) := \ft(u)$ and $\sgn(u^* \cdot \alpha_{\ft(u,u^*)}(x),u \cdot x) := \sgn(u)$.
\begin{remark}\label{remark: Formal notation}
    For the sake of clarity, we discuss the meaning of the formal notation involving $G_{n+m}^*$. Because our gluing rules are symmetric, it follows from Lemma \ref{lemma:CharEdges} that the edge set $E_{n+m}^{*}$ contains precisely the following edges. 
\begin{enumerate}
    \vspace{2pt}
    \item[(i)] $\{u \cdot x_1, u \cdot x_2\} \in E_{n+m}$ where $u \in V^n$ and $\{x_1,x_2\} \in E_m$.
    \vspace{2pt}
    \item[(ii)] $\{ u_1 \cdot x, u_2 \cdot \alpha_{\ft(u_1,u_2)}(x) \}$ where $u_1 \in V^n$, $\{ u_1,u_2 \} \in E_n^*$ and $x \in I_{\ft(u_1,u_2),\sgn(u_1,u_2)}^{m}$.
\end{enumerate}
Recall the notation $\sgn(x,y)$ from the beginning of the section.
\end{remark}

Let us now begin the construction of $\mathcal{J}_{n+m}$ in Theorem \ref{thm: Replacement flow}.
We begin by extending the unit flow $\mathcal{J}_n$ from $S_n$ to $T_n$ in $G_n$ into a unit flow $\mathcal{J}_n^*$ from $S_n^*$ to $T_n^*$ in $G_n^*$ as in Lemma \ref{lemma: Resistance after joining}. Namely we set
\[
    \mathcal{J}^*_n(u^*,u) := \divr(\mathcal{J}_n)(u).
\]
We first describe $\mathcal{J}_{n+m}$ locally. Fix $u \in V^n$ and write
\[
    N_{u}^+ := \{ y \in V_n^* : \{u,y\} \in E_n^* \text{ and } \mathcal{J}^*_n(u,y) > 0 \}
\]
and also
\[
    N_{u}^- :=  \{ y \in V_n^* : \{u,y\} \in E_n^* \text{ and }  \mathcal{J}^*_n(u,y) \leq 0\}.
\]
For each $y \in N_{u}^- \cup N_{u}^+$ we choose the members of $\mathcal{L}_m$ so that
\[
    L_{u,y} := I_{\ft(u,y), \sgn(u,y)}^{m} \text{ and } \widetilde{L}_{u,y} := I_{\ft(u,y), \sgn(y,u)}^{m}.
\]

Now, fix $y \in N_u^-$ and $z \in N_u^+$.
If $L_{u,y} \neq L_{u,z}$, we define the flow $\mathcal{J}_{y,z}$ by suitably embedding $\mathcal{B}_{[L_{u,y} \to L_{u,z}]}$ into $G^*_{n+m}$.
To be precise, first if an edge in $G_{n+m}^*$ does not contain an endpoint of the form $u \cdot x$, the flow on the edge is zero. Now choose an edge $\{z_1,z_2\} \in E_{n+m}^*$. According to Remark \ref{remark: Formal notation}, we may assume that $z_1 = u \cdot x_1 \in V^{n+m}$ and $z_2 = h \cdot x_2 \in V_{n+m}^*$ for $h \in V_n^*$ and $x_2 \in V^{m}$. If $u = h$ then
\[
    \mathcal{J}_{y,z}(z_1,z_2) := \mathcal{B}_{[L_{u,y} \to L_{u,z}]}(x_1,x_2).
\]
If $u \neq h$, we define $\mathcal{J}_{y,z}(z_1,z_2) = 0$ unless
\begin{equation}\label{eq: Replacement flow vertices between rooms}
	h \in \{ y,z \}, x_1 \in L_{u,y} \text{ and } x_2 = \alpha_{\ft(u,h)}(x_1).
\end{equation}
If \eqref{eq: Replacement flow vertices between rooms} holds, we set
\[
\mathcal{J}_{y,z}(z_1,z_2) := \mathcal{B}_{[L_{u,y} \to L_{u,z}]}(x_1,(x_1)_{L_{u,h}}).
\]

Then we consider the case $L_{u,y} = L_{u,z}$. We again let  $\mathcal{J}_{y,z}$ vanish on edges that do not contain an endpoint $(z_1,z_2)$. We use the same notation for $z_1,z_2,x_1,x_2,h$ as above.
If \eqref{eq: Replacement flow vertices between rooms} fails, then $\mathcal{J}_{y,z}(z_1,z_2) = 0$. If \eqref{eq: Replacement flow vertices between rooms} holds for $h = y$ we define
\[
    \mathcal{J}_{y,z}(z_1,z_2) = \mathcal{B}_{[L_{u,y} \to \widetilde{L}_{u,y}]}(x,x_{L_{u,y}}).
\]
For $h = z$,
\[
    \mathcal{J}_{y,z}(z_1,z_2) = \mathcal{B}_{[\widetilde{L}_{u,z} \to L_{u,z}]}(x,x_{L_{u,z}}).
\]
\begin{lemma}\label{lemma:ReplacementFlow(1)}
    For each pair $y \in N_u^-\, z \in N_u^+$ the function $\mathcal{J}_{y,z}$ as defined above is a unit flow from $\{ y \cdot \alpha_{\ft(u,y)}(x) : x \in L_{u,y} \}$ to $\{ z \cdot \alpha_{\ft(u,z)}(x) : x \in L_{u,z} \}$ in $G_{n+m}^*$ and it satisfies
    \begin{equation}
        \cE_q(\mathcal{J}_{y,z}) \leq \cE_q(\mathcal{B}).
    \end{equation}
    where $q = p/(p-1)$ is the dual exponent of $p$.
\end{lemma}

\begin{proof}
    First assume $L_{u,y} \neq L_{u,z}$. For vertices of the form $u \cdot x$ for $x \in V^m$, we have by construction
    \[
        \divr(\mathcal{J}_{y,z})(u \cdot x) = \divr(\mathcal{B}_{[L_{u,y} \to L_{u,z}]})(x) = 0.
    \]
    To see that $\divr(\mathcal{J}_{y,z})$ vanishes at the rest of the claimed vertices, simply observe that $\mathcal{J}_{y,z}$ vanishes at the edges containing them. The flow is unit because
    \[
        I(\mathcal{J}_{y,z}) = \sum_{x \in L_{u,y}} \mathcal{J}_{y,z}(y \cdot \alpha_{\ft(u,y)}(x), u \cdot x) = I(\mathcal{B}_{[L_{u,y} \to L_{u,z}]}) = 1,
    \]
    where the last equality follows from \ref{(FB1)}.
    The energy estimate follows from
    \[
      \cE_q(\mathcal{J}_{y,z}) = \cE_q(\mathcal{B}_{[L_{u,y} \to L_{u,z}]}) \leq \cE_q(\mathcal{B}).  
    \]
    
    Next assume $L_{u,y} = L_{u,z}$. By definition of $\mathcal{J}_{y,z}$,
    \[
    		I(\mathcal{J}_{y,z}) = \sum_{x
    \in L_{u,y}} \mathcal{J}_{y,z}(y \cdot \alpha_{\ft(u,y)}(x), u \cdot x) = I(\mathcal{B}_{[L_{u,y} \to \widetilde{L}_{u,y}]}) = 1.
    \]
    For $x \in L_{u,y}$ we compute
\begin{align*}
	\mathcal{J}_{y,z}(y \cdot \alpha_{\ft(u,y)}(x),u \cdot x) & = \mathcal{B}_{[L_{u,y} \to \widetilde{L}_{u,y}]}(x_{L_{u,y}},x) =  \mathcal{B}_{[\widetilde{L}_{u,y} \to L_{u,y}]}(x,x_{L_{u,y}})\\
	& = \mathcal{B}_{[\widetilde{L}_{u,z} \to L_{u,z}]}(x,x_{L_{u,z}}) = \mathcal{J}_{y,z}(u \cdot x,z \cdot \alpha_{\ft(u,z)}(x)),
\end{align*}
where the second equality follows from \ref{(FB4)} and the third from \ref{(FB3)}. So we get
    \[
        \divr(\mathcal{J}_{y,z})(u \cdot x) = \mathcal{J}_{y,z}(u \cdot x, y \cdot \alpha_{\ft(u,y)}(x)) + \mathcal{J}_{y,z}(u \cdot x,z \cdot \alpha_{\ft(u,z)}(v)) = 0.
    \]
    The divergences vanish at the other claimed vertices because, again, the flow vanishes at the edges containing them, and $\mathcal{J}_{y,z}$ is a unit flow.
    To obtain the energy estimate, we use \ref{(FB4)} to compute
    \[
        \cE_q(\mathcal{J}_{y,z}) = 2\sum_{x \in L_{u,y}} \left|\mathcal{B}_{[L_{u,y} \to \widetilde{L}_{u,y}]}(x_{L_{u,y}},x)\right|^q \leq \cE_q(\mathcal{B}_{[L_{u,y} \to \widetilde{L}_{u,y}]}) \leq \cE_q(\mathcal{B}).
    \]
\end{proof}

Next we introduce a flow $\mathcal{J}_{u}$ for the chosen $u \in V^n$ by taking a suitable linear combination $\mathcal{J}_{y,z}$. To this end, we first set
\[
    I_u := \sum_{z \in N_u^+} \mathcal{J}_n^* (u,z) = -\sum_{y \in N_u^-} \mathcal{J}_n^* (u,y).
\]
Note that the latter equality follows from $\divr(\mathcal{J}_n^*)(u) = 0$. Also, we clearly have $I_u = 0$ if and only if $\mathcal{J}^*_n(u,y) = 0$ for all $y \in V_n^*$.
In such case, we set $\mathcal{J}_u$ to be the zero flow. Otherwise, we must have $N_u^+,N_u^- \neq \emptyset$ and $J_u > 0$. The flow $\mathcal{J}_u$ is then given by
\[
    \mathcal{J}_u := \sum_{z \in N_u^+}\left(\frac{\mathcal{J}_n^*(u,z)}{I_u}\sum_{y \in N_u^-} \mathcal{J}_n^*(y,u)\cdot \mathcal{J}_{y,z} \right).
\]
\begin{proposition}\label{prop: Replacement flow on connecting edges}
    For all $h \in V_n^*$ so that $\{ u,h \} \in E_n^*$ and $x \in L_{u,h}$,
    \begin{equation*}
        \mathcal{J}_u(h \cdot \alpha_{\ft(u,h)}(x),u \cdot x) = \mathcal{J}_n^*(h,u) \cdot \mathcal{B}_{[L_{u,h} \to \widetilde{L}_{u,h}]}(x_{L_{u,h}},x).
    \end{equation*}
\end{proposition}

\begin{proof}
    There is nothing to prove if $I_u = 0$. Thus, we may assume $I_u>0$. First we consider the case $h \in N_u^{-}$. By the construction of $\mathcal{J}_{y,z}$, whenever $h \neq y \in N_u^-$ and $z \in N_u^+$, 
    \begin{equation*}
        \mathcal{J}_{y,z}(h \cdot \alpha_{\ft(u,h)}(x),u \cdot x) = 0.
    \end{equation*}
    By using \ref{(FB3)} we see from the construction of $\mathcal{J}_{h,z}$ for $z \in N_u^+$,
    \begin{equation*}
        \mathcal{J}_{h,z}(h \cdot \alpha_{\ft(u,h)}(x),u \cdot x) = \mathcal{B}_{[L_{u,h} \to \widetilde{L}_{u,h}]}(x_{L_{u,h}},x).
    \end{equation*}
    By combining the previous two equalities with the definitions of $\mathcal{J}_u$ and $I_u$,
    \begin{align*}
        \mathcal{J}_u(h \cdot \alpha_{\ft(u,h)}(x), u \cdot x)
        & = \sum_{z \in N_u^+}\left( \frac{\mathcal{J}_n^*(u,z)}{I_u}\sum_{y \in N_u^-}\mathcal{J}_n^* (y,u)\cdot \mathcal{J}_{y,z}(h \cdot \alpha_{\ft(u,h)}(x), u \cdot x)\right)\\
        & = \sum_{z \in N_u^+}\left( \frac{\mathcal{J}_n^*(u,z)}{I_u} \mathcal{J}_n^* (h,u)\cdot \mathcal{J}_{h,z}(h \cdot \alpha_{\ft(u,h)}(x), u \cdot x) \right)\\
        & = \left(\sum_{z \in N_u^+} \frac{\mathcal{J}_n^*(u,z)}{I_u}\right)
        \mathcal{J}_n^*(h,u) \mathcal{B}_{[L_{u,h} \to \widetilde{L}_{u,h}]}(x_{L_{u,h}},x)\\
        & = \mathcal{J}_n^*(h,u)\mathcal{B}_{[L_{u,h} \to L_{u,x}^{*}]}(x_{L_{u,h}},x).
    \end{align*}
    The case $h \in N_x^+$ follows similarly 
    \begin{align*}
        \mathcal{J}_u(h \cdot \alpha_{\ft(u,h)}(x), u \cdot x)
        & = \sum_{z \in N_u^+}\left( \frac{\mathcal{J}_n^*(u,z)}{I_u}\sum_{y \in N_u^-}\mathcal{J}_n^* (y,u)\cdot \mathcal{J}_{y,z}(h \cdot \alpha_{\ft(u,h)}(x), u \cdot x)\right)\\
        & = \sum_{y \in N_u^-}\left( \frac{\mathcal{J}_n^*(y,u)}{I_u} \mathcal{J}_n^* (u,h)\cdot \mathcal{J}_{y,h}(h \cdot \alpha_{\ft(u,h)}(x), u \cdot x) \right)\\
        & = \mathcal{J}_n^*(u,h)\cdot \mathcal{B}_{[L_{u,h} \to \widetilde{L}_{u,h}]}(x_{L_{u,h}},x)\\
        & = \mathcal{J}_n^*(h,u)\cdot \mathcal{B}_{[L_{u,h} \to \widetilde{L}_{u,h}]}(x,x_{L_{u,h}})
    \end{align*}
    where the second equality follows from \ref{(FB4)}.
\end{proof}

\begin{proposition}\label{prop: Local energy of replacement flow}
    The flow $\mathcal{J}_u$ satisfies the energy estimate
    \begin{equation}\label{eq: Local energy of replacement flow}
        \cE_q(\mathcal{J}_u) \leq (C_{\deg} + 1)^{\frac{1}{p - 1}}\cE_q(\mathcal{B}) \sum_{y \in N_u^-}\abs{\mathcal{J}_n^*(u,y)}^q.
    \end{equation}
\end{proposition}

\begin{proof}
    The claim is clear when $I_u = 0$, and we may assume $I_u > 0$. We recall that $\cE_q(\mathcal{J})^{\frac{1}{q}}$ is just the $\ell^q$-norm of $\mathcal{J} \in \R^{E_m^*}$.
    By applying Minkowski's and Hölder's inequality we obtain
    \begin{align*}
    \cE_q(\mathcal{J}_{u})^{\frac{1}{q}} & \leq
    \sum_{z \in N_u^+}\left( \frac{\mathcal{J}_n^*(u,z)}{J_u}\sum_{y \in N_u^-}\mathcal{J}_n^*(y,u)\cdot\cE_q(\mathcal{J}_{y,z})^{\frac{1}{q}}\right) && (\text{Minkowski's ineq}) \\
    & \leq \cE_q(\mathcal{B})^{\frac{1}{q}} \left(\sum_{z \in N_u^+} \frac{\mathcal{J}_n^*(u,z)}{I_u}\right) \sum_{y \in N_u^-}\mathcal{J}_n^*(y,u) && (\text{Lemma \ref{lemma:ReplacementFlow(1)}}) \\
    & = \cE_q(\mathcal{B})^{\frac{1}{q}} \sum_{y \in N_u^-}\mathcal{J}_n^*(y,u)\\
    & \leq \cE_q(\mathcal{B})^{\frac{1}{q}} \abs{N_u^-}^{\frac{1}{p}} \left(\sum_{y \in N_u^-}\abs{\mathcal{J}_n^*(y,u)}^q \right)^{\frac{1}{q}} && (\text{Hölder's ineq}) \\
    & \leq \cE_q(\mathcal{B})^{\frac{1}{q}} (C_{\deg} + 1)^{\frac{1}{p}} \left(\sum_{y \in N_u^-}\abs{\mathcal{J}_n^*(y,u)}^q \right)^{\frac{1}{q}}.
    \end{align*}
    In the last inequality, we used $\abs{N_u^-} \leq \deg(G_n^*) \leq C_{\deg} + 1$. The conclusion now follows by raising the previous inequality to the power $q$.
\end{proof}

We are now in the position of defining the flow $\mathcal{J}_{n+m}$.
Now, for $u \in V^n$, $z_1 = u \cdot v \in V^{n+m}$ and $z_2 \in V_{n+m}^*$, we set
\begin{equation}\label{eq: Def replacement flow locally}
    \mathcal{J}_{n+m}^*(z_1,z_2) := \mathcal{J}_u(z_1,z_2),
\end{equation}
and $\mathcal{J}_{n+m} := (\mathcal{J}_{n+m}^*)|_{V^{n+m} \times V^{n+m}}$. We will first verify that the expression \eqref{eq: Def replacement flow locally} defines a well-defined antisymmetric function on edges.
What needs to be done is to take $u_1,u_2 \in V^n$ such that $(u_1,u_2) \in E_n$ is an edge of type $t \in \cT$ and $x \in I_{t,+}^{m}$, and to verify $\mathcal{J}_{u_1}(z_1,z_2) = \mathcal{J}_{u_2}(z_1,z_2)$ for $z_1 = u \cdot x$ and $z_2 = u_2 \cdot \alpha_{t}(x)$. Note that here
\[
    L_{u_1,u_2} = \widetilde{L}_{u_2,u_1} = I_{t,+}^{m} \text{ and } L_{u_2,u_1} = \widetilde{L}_{u_1,u_2} = I_{t,-}^{m}.    
\]
We have
\begin{align*}
    \mathcal{J}_{u_1}(z_1,z_2) & = \mathcal{J}_n^*(u_1,u_2) \cdot \mathcal{B}_{[L_{u_1,u_2} \to \widetilde{L}_{u_1,u_2}]}(x_{L_{u_1,u_2}},x) && (\text{Proposition \ref{prop: Replacement flow on connecting edges}}) \\
    & = \mathcal{J}_n^*(u_1,u_2) \cdot \mathcal{B}_{[\widetilde{L}_{u_1,u_2} \to L_{u_1,u_2}]}(\alpha_t(x)_{\widetilde{L}_{u_1,u_2}},\alpha_{t}(x)) && \text{\ref{(FB4)}} \\
    & = \mathcal{J}_{u_2}(z_1,z_2) && (\text{Proposition \ref{prop: Replacement flow on connecting edges}}).
\end{align*}
This shows that $\mathcal{J}_{u_1}(z_1,z_2)=-\mathcal{J}_{u_2}(z_2,z_1)$, and that the function $\mathcal{J}_{n+m}$ defines an antisymmetric function.

\begin{proof}[Proof of theorem \ref{thm: Replacement flow}]  
We first show that $\widehat{F}_{n + m}$ is a unit flow from $S_{n+m}^*$ to $T_{n+m}^*$. Fix $z \in V^{n+m} \setminus (S_{n+m}^* \cup T_{n+m}^*) = V^{n+m}$. Namely $z = u \cdot x$ for some $u \in V^n$ and $x \in V^m$. Then, by Lemma \ref{lemma:ReplacementFlow(1)} and the definition of $\mathcal{J}_{n+m}^*$,
\[
    \divr(\mathcal{J}_{n+m}^*)(v) = \divr(\mathcal{J}_u)(z) = 0.
\]
The flow is unit by the computation
\begin{align*}
    I(\mathcal{J}_{n+m}^*) & = \ \ 
    \sum_{u \in S_{n}} \sum_{x \in L_{u,u^*}}
    \mathcal{J}_{n + m}^*(u^* \cdot \alpha_{\ft(u^*,u)}(x),u \cdot x)\\
    & = \sum_{u \in S_{n}} \sum_{v \in L_{u,u^*}} \mathcal{J}_u(u^* \cdot \alpha_{\ft(u^*,u)}(x),u \cdot x) && \text{(Def. \eqref{eq: Def replacement flow locally})} \\
    &  = \sum_{u \in S_n} \mathcal{J}_u(u^*,u)  \cdot  \sum_{v \in L_{u,u^*}} \mathcal{B}_{[L_{u,u^*} \to \widetilde{L}_{u,u^*}]}(x_{L_{u,u^*}},x) && (\text{Proposition \ref{prop: Replacement flow on connecting edges}}) \\
    & = \sum_{u \in S_n} \mathcal{J}_u(u^*,u) \cdot I(\mathcal{B}_{[L_{u,u^*} \to \widetilde{L}_{u,u^*}]}) = I(\mathcal{J}_n^*) = 1. && \ref{(FB1)}
\end{align*}
By the construction of $G_{n+m}^*$, $\mathcal{J}_{n+m}$ is then a unit flow from $S_{n+m}$ to $T_{n+m}$.

We finish the proof in establishing the desired energy estimate.
\begin{align*}
    \cE_q(\mathcal{J}_{n+m})^{\frac{1}{q}} & \leq  
    \cE_q(\mathcal{J}_{n+m}^*)^{\frac{1}{q}} \leq \sum_{u \in W_n} \cE_q(\mathcal{J}_{u})^{\frac{1}{q}} && \text{(Minkowski's ineq)}  \\
    & \leq C_1 \cE_q(\mathcal{B})^{\frac{1}{q}} \left( \sum_{u \in V^n} \sum_{y \in N_u^-} \abs{\mathcal{J}_n^*(u,y)}^{q} \right)^{\frac{1}{q}} && \text{(Proposition \ref{prop: Local energy of replacement flow})} \\
    & = C_1 \cE_q(\mathcal{B})^{\frac{1}{q}} \cE_q(\mathcal{J}_n^*)^{\frac{1}{q}}\\
    & \leq C_2 \cE_q(\mathcal{B})^{\frac{1}{q}} \cE_q(\mathcal{J}_n)^{\frac{1}{q}} && (\text{Lemma \ref{lemma: Resistance after joining}}) .
\end{align*}
The aimed for inequality now follows by raising the previous inequality to the power $q$.
The claimed continuity properties of the constant $C := C_2$ follows from Proposition \ref{prop: Local energy of replacement flow} and Lemma \ref{lemma: Resistance after joining}.
\end{proof}

\subsection{Construction of flow bases}\label{subsec:FB}

Next we construct flow bases for our main examples. We recall symmetric hypercubic IGSs from Definition \ref{Def:SymmHypercubic}.

\begin{theorem}\label{thm: Flow bases for cub}
    Let $\Xi$ be a symmetric hypercubic IGS for parameters $d,L_*,K$ and $p \in (1,\infty)$.
    Then there is a constant $C \geq 1$ depending only on $p,d,L_*,K$ such that for all $m \geq 1$ there exists a flow basis $\mathcal{B}$ of $G_m^e$ satisfying
    \[
        \mathcal{E}_q(\mathcal{B}) \leq C \cRes_p\left( I_{i,+}^m, I_{i,-}^m, G_m\right)
    \]
    where $i \in \mathcal{T}$ is any type.
    The constant $C$ can be chosen so that $C^{1-p}$ is continuous in $p \in (1,\infty)$ and has a limit as $p \to 1^+$.
\end{theorem}

\begin{proof}
    Fix $m \in \N$ and write for simplicity $L_{i,\pm} := I_{i,\pm}^m$. We first define the flows $\mathcal{B}_{[L_{i,\pm} \to L_{i,\mp}]}$ for all $1 \leq i \leq d$. Let $\mathcal{J}_i$ be the unique solution to $\cRes_p(L_{i,+},L_{i,-},G_m)$. Extend $\mathcal{J}_i$ to a flow $\mathcal{J}_i^*$ on $G_m^e$ in the natural way. More precisely, we set for $L = L_{i,\pm}$,
    \[
        \mathcal{J}_i^*(x_L,x) := \divr(\mathcal{J}_i)(x),
    \]
    and $\mathcal{J}_i^*(x_L,x) = 0$ for $L \neq L_{i,\pm}$. Now, set $\mathcal{B}_{[L_{i,+} \to L_{i,-}]} = \mathcal{J}_i^*$, $\mathcal{B}_{[L_{i,-} \to L_{i,+}]} = -\mathcal{J}_i^*$.

    Before constructing the other flows, let us first show that $\mathcal{J}_{i}^*$ and $-\mathcal{J}_{i}^*$ satisfy \ref{(FB4)} where the symmetry $\alpha_i$ is the horizontal symmetry. Because $\alpha_i$ is a graph isomorphism on $G_m$ and $\alpha_i(L_{i,\pm}) = L_{i,\mp}$ (see Remark \ref{rem:SymmHyper}), both $\alpha(\mathcal{J}_i)$ and $-\mathcal{J}_i$ solve $\cRes_p(L_{i,-},L_{i,+},G_m)$ where
    \[
    \alpha_i(\mathcal{J}_i)(x,y) := \mathcal{J}_i(\alpha(x),\alpha(y)).
    \]
    By the uniqueness of the energy minimizing unit flow, $\alpha_i(\mathcal{J}_i) = -\mathcal{J}_i$. This implies \ref{(FB4)}.

    Next, we construct the rest of the flows in the flow basis. Fix a distinct pair $i,j \in \mathcal{T}$ and consider the diagonal symmetry $\alpha_{i,j}^+$ defined on $G_m$. By the definition of symmetric hypercubic IGS, $\alpha_{i,j}^+$ is a graph isomorphism on $G_m$. Next we consider a graph $\widetilde{G}_m$ which is obtained by removing from $G_m^e$ vertices $x_{L_{k,\pm}}$ for all $k \notin \{i,j\}$ and the edges containing them. Note that $\mathcal{J}_i^*$ is supported on $\widetilde{G}_m$ so we understand $\mathcal{J}_i^*$ as a unit flow from $L_{i,+}^e$ to $ L_{i,-}^e$ on $\widetilde{G}_m$. We have a natural lift of $\alpha_{i,j}^+$ to $\widetilde{G}_m$ because $\alpha_{i,j}^+(L_{k,\pm}^e) \in \{ L_{i,+}^e,L_{i,-}^e,L_{j,+}^e,L_{j,-}^e \}$ for $k \in \{i,j\}$ by the same reasoning as in Remark \ref{rem:SymmHyper}. The lift is given by
    \[
        \alpha_{i,j}^+(x_{L_{k,\pm}}) = x_{\alpha(L_{k,\pm})}.
    \]
    
    We now claim that $\alpha_{i,j}^+$ on $\widetilde{G}_m$ is a reflectional symmetry as in Definition \ref{def:Reflection}. The partition $V(\widetilde{G}_m) = C \cup D \cup \Delta$ is given by
    \begin{align*}
        \Delta := & \{ x \in V(\widetilde{G}_m) : \alpha_{i,j}^+(x) = x \}\\
        = & \{ v_1\ldots v_m \in V^m \in c_{i}(v_k) = c_j(v_k) \text{ for all } i = 1,\dots,m \},
    \end{align*}
    \begin{align*}
        C := & \left\{ v_1\ldots v_m \notin \Delta : c_i(v_k) < c_j(v_k) \text{ where } k := \min_{l} c_i(v_l) \neq c_j(v_l) \right\} \\
        & \cup \{ x_{L_{i,+}} : x \in L_{i,-}  \} \cup \{ x_{L_{j,+}} : x \in L_{j,+}  \}
    \end{align*}
    and $D:=V(\widetilde{G}_m)\setminus(C\cup \Delta)$ contains the remaining vertices. It is immediate to verity $\alpha_{i,j}^+$ is a reflectional symmetry with respect to this partition.

    The flow $\mathcal{B}_{[L_{i,+} \to L_{j,-}]}$ is given by applying Proposition \ref{prop:TwistFlow} to $\mathcal{J}_{i}^*$ with $\alpha=\alpha_{i,j}^+$. Note that the source and target of $\mathcal{J}_{i}^*$ are $L_{i,+}^e$ and $L_{i,-}^e$, respectively, so the assumptions in Proposition \ref{prop:TwistFlow} are satisfied. Namely, $\mathcal{B}_{[L_{i,+} \to L_{j,-}]}$ is a well-defined unit flow from $L_{i,+} \to L_{j,-}$. Next we check that 
    the two pairs $\mathcal{B}_{[L_{i,+} \to L_{i,-}]}$ and $\mathcal{B}_{[L_{i,+} \to L_{j,-}]}$, and  $\mathcal{B}_{[L_{j,+} \to L_{j,-}]}$ and $\mathcal{B}_{[L_{i,+} \to L_{j,-}]}$ are compatible in the sense of \ref{(FB3)}. For the first pair, this follows immediately from the explicit form of the flow constructed in Proposition \ref{prop:TwistFlow} because for $x \in L_{i,+}$,
    \begin{align*}
        \mathcal{B}_{[L_{i,+} \to L_{j,-}]}(x_{L_{i,+}},x) & = \mathcal{J}_i^*(x_{L_{i,+}},x) + \mathcal{J}_i^*(\alpha(x)_{L_{j,+}},\alpha(x))\\
        & = \mathcal{J}_i^*(x_{L_{i,+}},x) = \mathcal{B}_{[L_{i,+} \to L_{i,-}]}(x_{L_{i,+}},x).
    \end{align*}
    For the second one, note that $\alpha_{i,j}^+(\mathcal{J}_i^*) = \mathcal{J}_j^*$ by the uniqueness of the energy minimizing unit flow, where $\alpha_{i,j}^+(\mathcal{J}_i^*)(x,y) := \mathcal{J}_i^*(\alpha_{i,j}^+(x),\alpha_{i,j}^+(y))$. Then we have
    \begin{align*}
        \mathcal{B}_{[L_{j,+} \to L_{j,-}]}(x_{L_{j,-}},x) & = \mathcal{J}_j^*(x_{L_{j,-}},x) = \mathcal{J}_i^*(\alpha(x_{L_{j,-}}),\alpha(x))\\
        & = \mathcal{J}_i^*(x_{L_{j,-}},x) +
        \mathcal{J}_i^*(\alpha(x_{L_{j,-}}),\alpha(x))\\
        & = \mathcal{B}_{[L_{i,+} \to L_{j,-}]}(x_{L_{j,-}},x).
    \end{align*}
    In the last equality we used the definition of the flow in Proposition \ref{prop:TwistFlow}.

    Finally, we construct the flows $\mathcal{B}_{[L_{i,+} \to L_{j,+}]}$ for distinct pairs $i,j \in \mathcal{T}$. We again use Proposition \ref{prop:TwistFlow} but this time we use the anti-diagonal symmetry $\alpha_{i,j}^{-}$. The verification of \ref{(FB3)} and \ref{(FB4)} are the same as above, and \ref{(FB1)} is trivially true. This concludes the construction of the flow basis $\mathcal{B}$.

    We still need to verify the energy estimate in the statement. First, note that the value $\cRes_p\left( I_{i,+}^m, I_{i,-}^m, G_m\right)$ is the same for all $i \in \mathcal{T}$. This follows by applying the diagonal symmetries to energy minimizing flows.
    By the construction in Proposition \ref{prop:TwistFlow} and the elementary inequality $(\abs{a} + \abs{b})^q \leq 2^{q-1}(\abs{a}^q + \abs{b}^q)$, therein constructed flow has at most $2^{q} = 2^{p/(p-1)}$ times larger energy than the original. Also by Lemma \ref{lemma: Resistance after joining},
    \[
    \mathcal{E}_q(\mathcal{J}_i^*) \leq 3C_{\deg}^{1/(p-1)}\mathcal{E}_q(\mathcal{J}_i) = 3C_{\deg}^{1/(p-1)} \cRes_p\left( I_{i,+}^m, I_{i,-}^m, G_m\right),
    \]
    where $C_{\deg}$ is the maximum degree of all replacement graphs $G_n$. 
    Note that $C_{\deg}$ is finite according to Lemma \ref{lemma: characterize bounded degree}. The desired energy estimate, and the continuity properties of the constant $C$, now follow by combining the previous estimates.
\end{proof}

\begin{theorem}\label{thm: Flow bases for pent}
    Let $\Xi$ be the Petagonal Sierpi\'nski carpet and $p \in (1,\infty)$.
    Then there is a constant $C \geq 1$ depending only on $p \in (1,\infty)$ and the geometry of $\Xi$ such that for all $m \geq 1$ there exists a flow basis $\mathcal{B}$ of $G_m^e$ satisfying
    \[
        \mathcal{E}_q(\mathcal{B}) \leq C \cRes_p\left( I_{t,+}^m, I_{t,-}^m, G_m\right)
    \]
    for all $t \in \mathcal{T}$.
    The constant $C$ can be chosen so that $C^{1-p}$ is continuous in $p \in (1,\infty)$ and has a limit as $p \to 1^+$.
\end{theorem}

\begin{proof}
    During the proof, the arithmetic on the symbols $\{ 0,1,2,3,4 \}$ is always modulo 5. 
    For simplicity we write $L_i := L_i^{m}$.
    We denote the energy minimizing unit flow from $L_i$ to $L_{i+2}$ in $G_m$, and from $L_i$ to $L_{i-2}$, as $\mathcal{J}_{i,+2}$ and $\mathcal{J}_{i,-2}$, respectively.
    As in the previous proofs, we extend these flows to $G_m^e$ and denote them $\mathcal{J}_{i,\pm 2}^*$.
    We lift symmetries $\alpha$ of the regular pentagon from the replacement graphs $G_m$ to $G_m^e$ so that $\alpha(x_{L_i}) = x_{\alpha(L_i)}$. These are clearly graph isomorphisms.

    Unfortunately, the  flows $\mathcal{J}_{i,\pm 2}$ can not be directly included in the flow basis $\mathcal{B}$, because they may fail to satisfy \ref{(FB3)}. This is because, for instance, $\mathcal{J}_{i,+2}$ and $\mathcal{J}_{i,-2}$ do not need to agree on the edges $\{ x_{L_i},x \}$. In fact, a fairly simple numerical calculation would verify this to be the case. Instead, the key idea will be to symmetrize the flows by taking an average $\mathcal{J}_i^* := 1/2( \mathcal{J}_{i,+2}^* + \mathcal{J}_{i,-2}^*)$, which is a unit flow $L_i^*$ to $L_{i+2}^* \cup L_{i-2}^*$. These flows are \emph{symmetric} in the following sense.
    Given a symmetry $\alpha$ of the regular pentagon and a flow $\mathcal{J}$ on $G_m^e$, we define the flow $\alpha(\mathcal{J})(x,y) := \mathcal{J}(\alpha^{-1}(x),\alpha^{-1}(y))$. We say that a flow $\mathcal{J}$ symmetric, if $\alpha(\mathcal{J})=\mathcal{J}$.
    Note that if $\alpha(L_i) = L_j$, then 
    \[
    \{ \alpha(L_{i+2}),\alpha(L_{i-2}) \} = \{ L_{j+2},L_{j-2} \}.
    \]
    By the uniqueness of the energy minimizing unit flow, we then have
    \[
        \{ \alpha(\mathcal{J}_{i,+2}), \alpha(\mathcal{J}_{i,-2}) \} = \{ \mathcal{J}_{j,+2}, \mathcal{J}_{j,-2} \}.
    \]
    This implies that $\alpha(\mathcal{J}_i^*) = \mathcal{J}_j^*$ and that the flows $\mathcal{J}_i^*$ are symmetric.

    For $i \in \{0,1,2,3,4\}$ let $\alpha_i$ be the non-trivial symmetry of the regular pentagon that fixes the vertex opposite to $L_i$. We observe that the isomorphisms $\alpha_i : G_m^e \to G_m^e$ satisfy the assumptions in Proposition \ref{prop:TwistFlow2}. To see this, first note that $\alpha_i^2 = \id_{V(G)}$ is clear.
    The partition $V(G_m^e) = C_{i,m} \cup D_{i,m} \cup \Delta_{i,m}$ is given as follows; see also Figure \ref{fig:PC}.
    Note that $\alpha_i$ only fixes one vertex, namely the constant sequence $v_1\ldots v_m \in V^m$ where $v_l = i$ for all $k$.
    For $G_1$ we have the natural partition
    \[
        V(G_1) = C_i \cup D_i \cup \Delta_i := \{ i-1,i-2 \} \cup \{ i+1,i+2 \} \cup \{ i \}.
    \]
    Then we define $\Delta_{i,m} = \{ i^m \}$, where $i^m$ is understood as the constant sequence,
    \[
        C_{i,m} \cap V^m := \{ v_1\ldots v_m \notin \Delta_{i,m} : v_k \in C_i \text{ where } k = \min\{l: v_l \neq i \}\}
    \]
    and
    \[
        D_{i,m} \cap V^m := \{ v_1\ldots v_m \notin \Delta_{i,m} : v_k \in D_i \text{ where } k = \min\{l: v_l \neq i \} \}.
    \]
    For $j \in \{0,1,2,3,4\}$ and $x \in L_j$, we define $x_{L_j} \in C_{i,m}$ if  $L_j \subseteq \Delta_{i,m} \cup (C_{i,m} \cap V^m)$. If $L_j \subseteq \Delta_{i,m} \cup (D_{i,m} \cap V^m)$ then $x_{L_j} \in D_{i,m}$.

    The flows in $\mathcal{B}$ are now constructed using Proposition \ref{prop:TwistFlow2} applied to $\mathcal{J}_i^*$ and $\alpha_j$ as follows. For $i \in \{0,1,2,3,4\}$, we use the symmetry $\alpha_{i+1}$ and the flow $\mathcal{J}_i^*$ to construct $\mathcal{B}_{[L_{i} \to L_{i+2}]}$. To construct $\mathcal{B}_{[L_i \to L_{i+1}]}$, we use $\alpha_{i-1}$ and apply it to the flow $\mathcal{J}_i^*$. The rest of the flows in the basis are constructed similarly. In order to verify the conditions of flow basis, we argue as follows. If $i \in \{0,1,2,3,4\}$, then the value $\abs{\mathcal{J}(x_{L_i},x)}$ is the same for all $\mathcal{J} \in \mathcal{B}$ such that $L_i$ is the source or target of $\mathcal{J}$. This follows from the explicit description of the flow constructed in Proposition \ref{prop:TwistFlow2} and the symmetry properties of $\mathcal{J}_i^*$ described above. The required conditions follow from this observation.

    Finally, we check the energy estimates. First we note that, by symmetry, the value $\cRes_p\left( I_{t,+}^m, I_{t,-}^m, G_m\right)$ is the same for all $t \in \mathcal{T}$. The desired energy estimate for the flow basis, and the continuity properties of the constant $C$, follows from a similar argument as in the proof of Theorem \ref{thm: Flow bases for cub}.
\end{proof}

\begin{proposition}\label{prop: Flow bases gask}
    Let $\Xi$ be the Sierpi\'nski gasket and $p \in (1,\infty)$.
    Then there is a constant $C \geq 1$ depending only on $p \in (1,\infty)$ and the geometry of $\Xi$ such that for all $m \geq 1$ there exists a flow basis $\mathcal{B}$ of $G_m^e$ satisfying
    \[
        \mathcal{E}_q(\mathcal{B}) \leq C \cRes_p\left( I_{t,+}^m, I_{t,-}^m, G_m\right)
    \]
    for all $t \in \mathcal{T}$.
    The constant $C$ can be chosen so that $C^{1-p}$ is continuous in $p \in (1,\infty)$ and has a limit as $p \to 1^+$.
\end{proposition}

\begin{proof}
    We set $\mathcal{B}$ to be the collection of flows which are obtained by extending the minimizers for
    $\cRes_p\left( I_{t,+}^m, I_{t,-}^m, G_m\right)$ to $G_m^e$ in an analogous manner we have done already multiple times. The collection of unit flows $\mathcal{B}$ already satisfies the conditions of flow basis because the sets $I_{t,\pm}^m$ only contain one vertex. Finally, because the Sierpi\'nski gasket has the symmetries of the equilateral triangle, the value
    $\cRes_p\left( I_{t,+}^m, I_{t,-}^m, G_m\right)$ is the same for all $t \in \mathcal{T}$. The desired energy estimate now follows from a similar argument as done in the previous two proofs.
\end{proof}

\section{Asymptotic geometric growth}\label{Section: Bounded geometry}
This section studies the geometry of replacement graphs.
The main goal of the section is the theorem below.
When $x \in V^n$ then $B(x,r) \subseteq V^n$ always denotes the open ball with respect to the path metric $d_{G_n}$.
We also recall the notation $x \cdot L$ from Subsection \ref{subsec:RR}.

\begin{theorem}\label{thm:Inclusions}
    Let $\Xi$ be a symmetric hypercubic IGS, the Sierpi\'nski gasket or the pentagonal carpet. Then there are constants $M_*,L_* \geq 1$ and $c_1,c_2 > 0$ depending on the given IGS so that the following inclusions hold. For all $n,m \in \N, \, x \in V^{n+m}$ and $r \leq c_1 L_*^m$,
    \[
        B(x,r) \subseteq \bigcup_{ z \in B(\pi_{n+m,n}(x),M_*)} z \cdot V^m,
    \]
    and for all $R \geq c_2L_*^m$,
    \[
        \bigcup_{ z \in B(\pi_{n+m,n}(x),M_*)} z \cdot V^m \subseteq B(x,R).
    \]
\end{theorem}

\subsection{Bounded geometry}
The first inclusion of Theorem \ref{thm:Inclusions} requires quite delicate arguments, and is another reason for us to restrict our focus to particular examples of IGS's instead of analyzing all cases simultaneously.
Our approach to this inclusion is based on studying the structure of a special class of paths that are defined below.

 Through the section, unless specified otherwise, $\Xi$ denotes a general IGS, and we refer to its data as stated in the paragraph below Definition \ref{def: Iterated graph system}. The replacement graphs are denoted $G_n$, and we assume that they have uniformly bounded degrees.

For a given path $\theta = [x_1,\ldots x_k]$ in $G_n$ and $1 \leq m \leq n$ the \emph{projected path}, denoted $\pi_{n,m}(\theta)$, is the path in $G_m$ obtained by removing the consecutive repeated tuples from the sequence $[\pi_{n,m}(x_1),\ldots \pi_{n,m}(x_k)]$.

\begin{definition}
    For $n > 1$ a path $\theta$ on $G_n$ is \emph{non-collapsing} if $\len(\theta) = \len(\pi_n(\theta))$. We say that an IGS is of \emph{bounded geometry}, if 
    \[
        \sup_{n \in \N} \max_\theta \diam(\theta,d_{G_n}) < \infty
    \]
    where the maximum is taken over all non-collapsing paths in $G_n$. Here $\diam$ refers to the diameter, with respect to the path metric, of the set of vertices contained in the path $\theta$.
\end{definition}

Verifying the bounded geometry condition seems to be unfortunately technical and tedious, and this work verifies it for the main examples only. We do not, however, exclude the possibility that there exists some relatively simple algorithm that could decide whether it holds or not based on the data of IGS only.

\begin{figure}[!ht]
\centering\includegraphics[width=350pt]{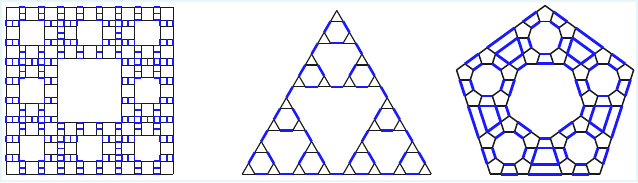} 
    \caption{The non-collapsing paths are exactly the paths which contain only blue edges.}
    \label{fig: Non-collapsing paths}
\end{figure}

\subsection{Bounded geometry}\label{Subsection: Verification of the bounded geometry condition}
We begin with the following technical lemma. For a given edge $\{ x,y \} \in E_n$ we denote $\Lambda(\{x,y\})$ the smallest number $m \geq 1$ so that $\pi_{n,m}(x) \neq \pi_{n,m}(y)$. For a given path $\theta = [x_1,\ldots x_k]$ we denote 
\[
    \Lambda(\theta) := \{ \Lambda(\{x_i,x_{i+1}\}) : 1 \leq i \leq k-1  \}.
\]

\begin{lemma}\label{lemma: Non-collapsing paths and bounded degree}
    Let $\theta = [x_1,\ldots x_k]$ be a path in $G_m$ and write $x_i := v_1^i\ldots v_m^i$. Assume that $J\subseteq \{1,\dots,m\}$ is a (discrete) sub-interval with $\Lambda(\theta) \cap J = \emptyset$. Then whenever
    $v_i^1 = v_j^1$ for some $i,j \in J$, we have $v_i^l = v_j^l$ for all $l = 1,\dots,k$.
\end{lemma}

\begin{proof}
    It is sufficient to prove the case where $k = 2$, so we write $\theta = [x_1,x_2]$. We may assume that $(x_1,x_2) \in E_m$ because otherwise we can consider the reversed path.

    Now, assume $v_i^1 = v_j^1$ for some $i,j \in J$. Since $J$ is an interval with $\Lambda(\{x_1,x_2\}) \notin J$, either
    $\Lambda(\{x_1,x_2\}) > i,j$ or $\Lambda(\{x_1,x_2\}) < i,j$.
    In the former case, by the definition of $\Lambda(\{x_1,x_2\})$, $v_i^2 = v_i^1 = v_j^1 = v_j^2$.
    In the latter case, by Lemma \ref{lemma:CharEdges}, we have $v_j^2 = \alpha_t(v_j^1) = \alpha_t(v_i^1) = v_1^2$.
\end{proof}

\begin{proposition}\label{prop: BIP for Gamecube}
    Hypercubic IGSs are of bounded geometry.
\end{proposition}

\begin{proof}
    Let $\Xi$ be a hypercubic IGS.
    Fix a non-collapsing path $\theta = [x_1,\ldots,x_k]$ in $G_m$ and write $x_i = v_1^i\ldots v_m^i$.
    We first show that $\abs{\Lambda(\theta)} \leq \abs{\mathcal{T}}$. More specifically, we show that every edge of the same type in $\theta$ have the same value in $\Lambda$.

    We fix an edge on $\theta$. Without loss of generality assume that we chose $\{x_1,x_2\}$ and let its type be $t$.
    Let $l$ be the smallest integer so that $1 < l \leq k$ and $\{ x_l,x_{l+1} \}$ is also of type $t$. We would be done if we can show that $\Lambda(\{ x_l,x_{l+1} \}) = \Lambda(\{ x_{1},x_{2} \})$.
    To do this, we do a counter-assumption. Because we can reverse the path, there is no harm in assuming
    $\Lambda(\{ x_1,x_{2} \}) < \Lambda(\{ x_{l},x_{l+1} \})$.
    
    Note that by the gluing rules of hypercubic IGSs and Lemma \ref{lemma:CharEdges}, if $x = v_1\ldots v_m, y = w_1\ldots w_m \in V^m$ are neighbors in $G_m$ and are connected by an edge of type $t' \neq t$, we have $c_t(v_i) = c_t(w_i)$ for all $i = 1,\ldots m$. Here $c_t$ are the coordinate functions. By the choice of $l$, it therefore holds that $c_t(v_i^2) = c_t(v_i^l)$ for all $i$.
    By Lemma \ref{lemma:CharEdges}, the gluing rules of hypercubic IGSs and $\Lambda(\{ x_1,x_{2} \}) < \Lambda(\{ x_{l},x_{l+1} \})$, $c_t(v_{j}^l) = c_t(v_{j}^2)$ are all the same value (either $1$ or $L_*$) for all $j \geq \Lambda(\{ x_{l},x_{l+1} \})$. This finally contradicts Lemma \ref{lemma:CharEdges}, the definition of edges and gluing rules of hypercubic IGSs and the fact that $\{ x_{l},x_{l+1} \}$ is an edge of type $t$. We now conclude that $\abs{\Lambda(\theta)} \leq \abs{\mathcal{T}}$.

    Now, fix $x = v_1\ldots v_m,\, y = w_1\ldots w_m \in V^m$ so that there is a non-collapsing path between $x$ and $y$.
    Write $\Lambda(\theta) = \{ l_1,\dots,l_k \}$ so that $l_1 < \dots < l_k$, and consider the discrete sub intervals 
    \[
        J_j :=
        \begin{cases}
            [1,l_1) & \text{ if } j = 0\\
            (l_k, m] & \text{ if } j = k\\
            (l_j,l_{j + 1}) & \text{ if } 1 \leq j \leq k - 1.
        \end{cases}
    \]
    Then write $U_j := \{ v_i : i \in J_j \}$ and for each $v \in U_j$ set $j_v$ to be the smallest index in $J_j$ so that $v_i = v$. By Lemma \ref{lemma: Non-collapsing paths and bounded degree} $y$ is uniquely determined by $x$ and $( w_{j_v} )_{ 1 \leq j \leq k,\, v \in U_j  }$. Since $k \leq \abs{\cT}$, there are at most $\abs{V}^{ \abs{\cT} (\abs{V} + 1) }$ many choices for $y \in V^m$. Because paths are connected, we get a bound for the diameters of non-collapsing paths.
\end{proof}

Next we deal with the Sierpi\'nski gasket and the pentagonal carpet. The argument for the former is direct, but the latter is considerably more tedious.

\begin{proposition}
    The Sierpi\'nski gasket and the pentagonal carpet are of bounded geometry.
\end{proposition}

\begin{proof}
    We first deal with the gasket. Let $\theta$ be a simple non-constant non-collapsing path in $G_m$. We will show that it contains exactly one edge. Let
    $x = v_1\ldots v_m,\, y = w_1\ldots w_m \in V^m$ be so that $(x,y) \in E_m$ is an edge contained in $\theta$. Using the symmetries of the equilateral triangle, we may assume that $(v_k,w_k) = (0,1)$, where $k = \Lambda(\{x,y\})$.
    By the gluing rules, $(v_i,w_i) = (1,0)$ for all $k < i \leq n$.
    From this, the fact that $\theta$ is non-collapsing and the gluing rules of the gasket we see that $\theta$ cannot contain other vertices than $x,y$.

    Next we move on to the pentagonal carpet.
    Let $\theta$ be a simple non-collapsing path in $G_m$ and $(x,y) \in E_m$ an edge contained in the path with
    \begin{equation}\label{eq: BIP for pentagon}
        k :=  \Lambda(\{x,y\}) = \min \Lambda(\theta).
    \end{equation}
    By applying the symmetries of the regular pentagon, we may assume that $(v_k,w_k) = (0,1)$.
    We will first prove that $\theta$ is contained in
    \[
        \left(\pi_{n,k}(x) \cdot I_{a,+}^{(n - k)}\right) \quad \bigcup \quad \left(\pi_{n,k}(y) \cdot I_{a,-}^{n-k}\right).
    \]
    Recall that $z \cdot A$ is the image of $A$ under the map $w \mapsto z \cdot w$.
    Note that $x$ and $y$ are contained in this set. Assume the converse, i.e., there is an edge $\{ z_1,z_2 \}$ in $\theta$ so that
    \[
        z_1 \in \left(\pi_{n,k}(x) \cdot I_{a,+}^{n-k}\right) \quad \bigcup \quad \left( \pi_{n,k}(y) \cdot I_{a,-}^{n-k}\right).
    \]
    and
    \[ 
        z_2 \notin \left(\pi_{n,k}(x) \cdot I_{a,+}^{n-k}\right) \quad \bigcup \quad \left( \pi_{n,k}(y) \cdot I_{a,-}^{n-k}\right).
    \]
    Write $z_j = u^{j}_1\dots u^{j}_n$ for $j = 1,2$ and, for simplicity, we assume that $z_1 \in \pi_{n,k}(x) \cdot I_{a,+}^{n-k}$. Observe that
    \[
        \pi_{n,k}(z_1) = \pi_{n,k}(x) \neq \pi_{n,k}(z_2).
    \]
    Hence $ \Lambda(\{z_1,z_2\}) \leq k$. But since $\{ z_1,z_2 \}$ is contained in $\theta$, by \eqref{eq: BIP for pentagon} we must have $\Lambda(\{z_1,z_2\}) = k$.
    Since $u_k^{1} = v_k = 0$, we must have $u_k^{2} \neq 1$. Otherwise, $(u_k^{1}, u_k^{2}) = (0,1)$ and 
    $u^{1} \in \pi_{n,k}(x) \cdot I_{a,+}^{(n - k)}$.
    In particular, we must have $u_k^{2} = 4$ and the edge is oriented as $(z_1,z_2)$.
    However, then $u^{1}_{k + 1} \in I_{a,+} \cap I_{e,+} = \emptyset$ which yields the contradiction. We thus have that claimed inclusion.

    We recall the notation $\Span{A}$ for the induced subgraph of $A \subseteq V^n$. We consider the subgraphs
    \[
        H_1 = \left\langle \left(\pi_{n,k}(x) \cdot I_{a,+}^{n-k}\right) \quad \bigcup \quad \left(\pi_{n,k}(y) \cdot I_{a,-}^{(n - k)}\right) \right\rangle \subseteq G_n
    \]
    and
    \[
        H_2 = \left\langle  \pi_{n,k}(x) \cdot I_{a,+}^{n-k} \right\rangle \subseteq G_n.
    \]
    Consider the simplicial map between graphs $f: H_1 \to H_2$,
    \[
        f(u \cdot v) =
        \begin{cases}
            \pi_{n,k}(x) \cdot v & \text{ if } u = \pi_{n,k}(x) \\
            \pi_{n,k}(x) \cdot \alpha_a(v) & \text{ if } u = \pi_{n,k}(y).
        \end{cases}
    \]
    Then the path $\theta'$ which is obtained by mapping $\theta$ by $f$ and removing the consequtive tuples, is a non-collapsing path contained in $\pi_{n,k}(x) \cdot I_{a,+}^{n-k}$. By proving that $\theta'$ contains at most 2 vertices, we can conclude that $\theta$ contains at most 4 vertices and we would be done.

    Let $(z_1,z_2)$ be an edge in $\theta'$ and $z_j = u^{j}_1\dots u^{j}_n$ and $l := \Lambda(\{z_1,z_2\})$.
    Note that
    \[
      \pi_{n,k}(z_1) = \pi_{n,k}(x) = \pi_{n,k}(z_2) 
    \]
    so we must have $l > k$. In particular, $u_l^1,u_l^2 \in I_{a,+}$, which gives $(u_l^1,u_l^2) = (1,2)$. Moreover, $u_i^1 \in I_{a,+} \cap I_{b,+} = \{2\}$ and $u^2_i \in I_{a,+} \cap I_{b,-} = \{ 1 \}$ for all $l < i \leq n$.
    Now suppose $\theta'$ contains an edge $\{ z_1,z_3 \}$. By replicating the previous argument, we have that $z_3 = z_2$. Similarly, if $\theta'$ contains $\{ z_2,z_3 \}$, we have $z_3 = z_1$. We may now conclude that $\theta'$ contains at most two vertices, so we are done.
    \end{proof}

\subsection{Asymptotic geometry}\label{Subsection: Asymptotic geometric growth of diameter}
We move on to proving Theorem \ref{thm:Inclusions}. We first perform some case-by-case analysis and derive suitable diameter upper bounds.

Let us first consider a symmetric hypercubic IGS $\Xi$.
For each $k=1,\dots, L_*$ denote with the bolded number $\bm{k}\in V$ the vertex, for which
    \[
        \begin{cases}
            c_1(\bm{k}) = k\\
            c_i(\bm{k}) = 1 & \text{ for all } i \neq 1\\
            \xi(\bm{k}) = 1.
        \end{cases}
    \]
    Such a vertex $\bm{k}$ is contained in $V$ by \ref{Def:SymmHypercubic(1)}.
    By \ref{Def:SymmHypercubic(2)}, there are edges $(\bm{k},\bm{k + 1}) \in E_1$ of type $t = 1$ for all $1 \leq k < L_*$.
    By \ref{Def:SymmHypercubic(3)}, $\bm{L_*} \in I_{1,+}$ and $\alpha(\bm{L_*}) = \bm{1}$.
    The graph $G_{m + 1}$ contains the edges $(\bm{k} \cdot (\bm{L_*})^m,(\bm{k + 1})\bm{1}^m)$.
    Here $\bm{k}^m \in V^m$ refers to the constant sequence whose all indices are equal to $\bm{k}$.
    
\begin{lemma}\label{lem: cornerdist}
    Let $\Xi$ be a symmetric hypercubic IGS.
    Then for all $m \in \N$,
    \begin{equation}\label{eq: Diam 1}
        d_{G_{m}}(\bm{1}^m,\bm{L_*}^m) \leq L_*^m - 1.
    \end{equation}
\end{lemma}

\begin{proof}
    Clear for $m = 1$. If it holds for $m \in \N$ (IH) then
    \begin{align*}
        & d_{G_{m + 1}}(\bm{1}^{m + 1},\bm{L_*}^{m + 1})\\
        \leq \quad & \sum_{k = 1}^{L_*} d_{G_{m + 1}}(\bm{k}\bm{1}^m,\bm{k}\bm{L_*}^m) + \sum_{k = 1}^{L_* - 1}d_{G_{m + 1}}(\bm{k}\bm{L_*}^{m},(\bm{k + 1})\bm{1}^m)\\
        \leq \quad & \sum_{k = 1}^{L_*} d_{G_{m}}(\bm{1}^m,\bm{L_*}^m) + L_* - 1\\
        \stackrel{\mathmakebox[\widthof{=}]{\text{(IH)}}}{\leq} \quad &  L_*\cdot (L_*^{m} - 1) + L_* - 1 = L_*^{m+1} - 1.
    \end{align*}
    In the second inequality we used $d_{G_{n+m}}(x \cdot y,x \cdot z) \leq d_{G_{n}}(y,z)$ which follows from the self-similarity of the graphs, Lemma \ref{lemma:CharEdges}.
\end{proof}

\begin{corollary}\label{cor: tilediam}
    Let $\Xi$ be a symmetric hypercubic IGS.
    Then for all $n,m \in \N$ and $x,y \in V^m$,
    \begin{equation}\label{eq: diam 2}
        d_{G_{n + m}}(x \cdot \bm{1}^{n},y \cdot \bm{1}^{n}) \leq d_{G_m}(x,y) \cdot L_*^{n}.
    \end{equation}

\end{corollary}
\begin{proof}
    It is sufficient to prove the case where $(w,v) \in E(G_m)$. Moreover, by \ref{Def:SymmHypercubic(4)}, we may assume that this edge is of type $t = 1$, since in the other cases one can use symmetries $\alpha_{1,i}^+$ which preserves the distances and fixes the point $\bm{1}$ in $V$. Then
    \begin{align*}
        d_{G_{n + m}}(x \cdot \bm{1}^{n},y \cdot \bm{1}^{n})
        & \leq d_{G_{n + m}}(x \cdot \bm{1}^{n},x \cdot \bm{L_*}^{n}) + d_{G_{n + m}}(x \cdot \bm{L_*}^{n},y \cdot \bm{1}^{n})\\
        & \leq d_{G_n}(\bm{1}^{n},\bm{L_*}^{n}) + 1\\
        & \leq L_*^{n} - 1 + 1 = L_*^n.
    \end{align*}
    The first inequality is triangle inequality, second is self-similarity and third is Lemma \ref{lem: cornerdist}.
\end{proof}

\begin{proposition}\label{prop:upperdiam}
Let $\Xi$ be a symmetric hypercubic IGS.
Then for all $m \in \N$ and $x,y \in V^m$,
    \begin{equation}\label{eq: diam 3}
        d_{G_m}(x,y) \leq 4\diam(G_1)\cdot L_*^m.
    \end{equation}   
\end{proposition}

\begin{proof}
   We first set $x = \bm{1}^m$ and also write $y = v_1\ldots v_m$. In the computation below we use the formal notations $\pi_{m,0}(y) \cdot z = z$ and $z \cdot \bm{1}^0 = z$. Now we estimate 
    \begin{align*}
        d_{G_m}(\bm{1}^m,y)
        \leq & \quad \sum_{j = 0}^{m - 1} d_{G_m}\left(\pi_{m,j}(y) \cdot \bm{1}^{m - j},\pi_{m,j+1}(y) \cdot \bm{1}^{m - (j + 1)}\right)\\
        = & \quad \sum_{j = 0}^{m - 1} d_{G_m}(\pi_{m,j}(y) \cdot \bm{1} \cdot \bm{1}^{m - (j + 1)},\pi_{m,j}(y) \cdot v_{j + 1} \cdot \bm{1}^{m - (j + 1)})\\
        \leq  & \quad \sum_{j = 0}^{m - 1} d_{G_{m - j}}(\bm{1}\cdot \bm{1}^{m - (j + 1)},v_{j + 1}\cdot \bm{1}^{m - (j + 1)})\\
        \leq  & \quad \diam(G_1) \cdot \sum_{j = 0}^{m - 1} L_*^{m - (j + 1)} \\
        \leq & \quad 2\diam(G_1)\cdot L_*^m. 
    \end{align*}
    In the second inequality, we used self-similarity and Corollary \ref{cor: tilediam} in the third.
    Thus, the aimed for inequality follows from
    \[
        d_{G_m}(x,y) \leq d_{G_m}(x,\textbf{1}^m) + d_{G_m}(\textbf{1}^m,y) \leq 4\diam(G_1)\cdot L_*^m. 
    \]
\end{proof}

Next we adjust the previous argument to get the desired upper bound for the Sierpi\'nski gasket and pentagonal carpet.

\begin{lemma}
    If $\Xi$ is the Sierpi\'nski gasket and $i,j \in \{ 0,1,2 \}$,
    \begin{equation*}
        d_{G_m}(i^m,j^m) \leq 2^m - 1.
    \end{equation*}
    If $\Xi$ is the pentagonal carpet and $i,j \in \{ 0,1,2,3,4 \}$ then
    \begin{equation*}
        d_{G_m}(i^m,j^m) \leq 2^{m+1} - 1.
    \end{equation*}
\end{lemma}

\begin{proof}
    First we give the argument for the Sierpi\'nski gasket. By symmetry we may assume $i = 0$ and $j = 1$. Since $(0,1) \in E$ is an edge of type $a$ and $(1,0) \in I_a$, the upper bound follows from an identical computation as in the proof of Lemma \ref{lem: cornerdist}.

    Next we consider the case of the pentagonal Sierpi\'nski carpet. First we assume that $i,j \in \{ 0,1,2,3,4 \}$ so that $j = i \pm 1 \mod 5$. By symmetry, we may assume that $i = 0$ and $j = 1$. Again, since $(0,1) \in E$ is an edge of type $a$ and $(1,0) \in I_a$, we can yet again conclude that
    \[
        d_{G_m}(0^m,1^m) \leq 2^m - 1.
    \]
    It now follows that if $j = i \pm 2 \mod 5$ and $k := i \pm 1 \mod 5$ we have
    \[
        d_{G_m}(i^m,j^m) \leq d_{G_m}(i^m,k^m) + d_{G_m}(k^m,j^m) \leq 2^m - 1 + 2^m - 1 \leq 2^{m+1} - 1.
    \]
\end{proof}

\begin{proposition}\label{prop: Diam UB Gasket and Pentagon}
    If $\Xi$ is the Sierpi\'nski gasket or the pentagonal carpet, then for all $m \in \N$ and $x,y \in V^m$,
    \begin{equation}\label{eq: Diam UP for both Gasket and Pentagon}
        d_{G_m}(x,y) \leq 32 \cdot L_*^m.
    \end{equation}
\end{proposition}

\begin{proof}
    We derive \eqref{eq: Diam UP for both Gasket and Pentagon} only using the estimate
    \[
        d_{G_m}(i^m,j^m) \leq 2^{m+1} - 1,
    \]
    and some properties of the IGSs satisfied by both examples.
    By the identical argument as given in Corollary \ref{cor: tilediam}, we have
    \[
        d_{G_{n + m}}(x \cdot 0^n,y \cdot 0^n) \leq d_{G_m}(x,y)\cdot 2^{m + 1}. 
    \]
    By following the argument given in Proposition \ref{prop:upperdiam}, we obtain
    \[
        d_{G_m}(0^m,y) \leq 4\diam(G_1) \cdot 2^{m+1} = 8 \cdot 2^{m+1}.
    \]
    The desired inequality then follows from the triangle inequality.
\end{proof}

We now have suitable diameter upper bounds. Next we deal with some lower bounds.

\begin{proposition}\label{prop: Distance between opposite faces}
Let $\Xi$ be a hypercubic IGS, Sierpi\'nski gasket or pentagonal carpet. Then for all $m \in \N$ and $t \in \mathcal{T}$,
    \begin{equation*}
        d_{G_m}\left(I_{t,+}^{m},I_{t,-}^{m}\right) \geq 1/2 \cdot L_*^m.
    \end{equation*}
\end{proposition}

\begin{proof}
    We deal with the hypercubic case first. Our method uses the following simple observation. Let $c_i$ be the coordinate functions for $i \in \mathcal{T}$ and consider the interval graph $\mathbb{I}_m$ whose vertex set is $\{0,1,\ldots L_*^m - 1\}$ and $\{k,l\} \in E(\mathbb{I}_m)$ if and only if $\abs{k-l} = 1$. Then the map $\varphi_{m,t} : V^m \to V(\mathbb{I}_m)$,
    \[
        \varphi_{m,t}(v_1\ldots v_m) := \sum_{i = 1}^{m} (c_t(v_i) - 1) \cdot L_*^{m-i}
    \]
    is simplicial $G_m \to \mathbb{I}_m$. These maps are counterparts of the coordinate projections $(x_1,\ldots,x_d) \mapsto x_i$ of $\R^d$. 
    The fact that $\varphi_{m,t}$ is simplicial is verified using the definition of the gluing rules and Lemma \ref{lemma:CharEdges}.
    
    Now, because $\varphi_m$ is simplicial, 
    \[
        d_{G_m}(x,y) \geq  d_{\mathbb{I}_m}(\varphi_m(x),\varphi_m(y)).
    \]
    Since $\varphi_m(I_{t,-}^m) = 0$ and $\varphi_m(I_{t,+}^m) = L_*^m-1$, this concludes the hypercubic case.
    
    For the gasket, we will show by induction on $m \in \N$ that
    \[
        d_{G_m}(0^m,1^m) \geq 2^m - 1
    \]
    for all $m \in \N$. By symmetry, this is sufficient.

    The claim is clear for $m = 1$. Assume that we have it for $m \in \N$ and let $\theta$ be a path connecting $0^{m + 1} $ to $ 1^{m + 1}$.
    Consider the non-trivial symmetry of the equilateral triangle $\alpha$ which fixes $0^{m+1}$. This induces a graph isomorphism $G_m \to G_m$ according to  Remark \ref{rem:gasketSYmm}. It also satisfies the conditions in Proposition \ref{prop:ReflFold} where the partition $V(G) = C \cup D \cup \Delta$ is $\Delta = \{0^{m+1}\}$,
    \[
        C = \{ v_1\ldots v_{m+1} \notin \Delta : v_k = 1 \text{ where } k = \min_{l} v_l \neq 0 \}
    \]
    and $D$ contains the rest.
    Now, let $f$ be as in Proposition \ref{prop:ReflFold} and $\theta'$ be the path obtained by mapping $\theta$ by $f$ and removing the consecutive tuples.
    We observe that $\theta'$ is contained in $0 \cdot V^m \cup 1 \cdot V^m$ and that $\{0 \cdot 1^m, 1 \cdot 0^m\}$ is the only edge connecting $0 \cdot V^m$ to $1 \cdot V^m$. So $\theta'$ must contain disjoint sub-paths $\theta_1$ connecting $0^{m+1}$ to $0 \cdot 1^m$ and contained in $0 \cdot V^m$, and $\theta_2$ connecting $1 \cdot 0^m$ to $1^{m+1}$ and contained in $1 \cdot V^m$. It also must contain the edge $\{0 \cdot 1^m,1\cdot 0^m \}$.
    By applying the induction hypothesis,
    \[
        \len(\theta) \geq \len(\theta') \geq \len(\theta_1) + \len(\theta_2) + 1 \geq 2^m -1 + 2^m - 1 + 1= 2^{m+1}.
    \]
    We conclude the argument for the gasket.

    Next, we deal with the pentagonal carpet; recall the notation from Example \ref{ex: SP}
    We show by induction that
    \[
        d_{G_m}(L_0^m, L_2^m \cup L_3^m) \geq 2^m - 1,
    \]
    which is sufficient by symmetry.
    The base case is yet again trivial.
    The inductive step is almost the same as for the gasket. We take the non-trivial symmetry of the regular pentagon that fixes the vertex $0^{m+1}$, and similarly fold the path onto a suitable smaller subset. Then we use induction hypothesis to find two disjoint paths with suitable length lower bound and add them together.
\end{proof}

\begin{proof}[Proof of Theorem \ref{thm:Inclusions}]
    We begin with the first inclusion, and we shall prove it for
    \[
        M_* := \sup_{n \in \N} \max_{\theta} \diam(\theta,G_n) + 1
    \]
    where the maximum is taken over all non-collapsing paths $\theta$ in $G_n$.
    Now, let $y \in V^m$ so that $d(\pi_{n+m,m}(x),y) \geq M_*$ and choose any $w \in y \cdot V^n$ with a path $\theta$ from $x$ to $w$. Let $\mathcal{C}$ and $\fold_m(\theta)$ be as in Lemma \ref{lemma:foldingFar}. First note that the condition for $\mathcal{C}$ in Lemma \ref{lemma:foldingFar} must fail, namely there is no common vertex for all the members in $\mathcal{C}$ because otherwise, the path $[h_1,\ldots,h_k]$ in Lemma \ref{lemma:foldingFar} would be non-collapsing. But this would contradict the choice of $M_*$.
    
    For these three classes of examples, the only possibility for the members in $\mathcal{C}$ to not have a common vertex, is that the folded path $\fold_m(\theta)$ connects the sets $\pi_{n+m,n}(x) \cdot I_{t,+}^m$ and $\pi_{n+m,n}(x) \cdot I_{t,-}^m$ for some $t \in \mathcal{T}$. Since the folded path is contained in $\pi_{n+m,n}(x) \cdot V^m$, it follows from Proposition \ref{prop: Distance between opposite faces} that
    \[
        \len(\theta) \geq \len(\fold_m(\theta)) \geq 1/2L_*^m.
    \]
    In conclusion, we have the first inclusion
    \[
    B(x,r) \subseteq \bigcup_{ z \in B(\pi_{n+m,n}(x),M_*)} z \cdot V^m
    \]
    when $r \leq 1/4L_*^m$.

    The latter inclusion, which reads
    \[
        \bigcup_{ z \in B(\pi_{n+n,m}(x),M_*)} z \cdot V^m \subseteq B(x,R),
    \]
    follows from the self-similarity of the graphs (Lemma \ref{lemma:CharEdges}) and the inequality
    \[
        \diam(V^m,G_m) \leq C L_*^{m},
    \]
    which was verified above for the examples.
\end{proof}

We need a couple more geometric properties. The first one regards the self-similarity of the path metric.

\begin{proposition}\label{prop:MetricSS}
    Let $\Xi$ be a symmetric hypercubic IGS, the Sierpi\'nski gasket or the pentagonal carpet. Then for all $n,m \in \N$ and $x \in V^n$, $y,z \in V^m$, the path metrics satisfy
    \[
        d_{G_{n+m}}(x \cdot y,x \cdot z) = d_{G_{m}}(y,z).
    \]
\end{proposition}

\begin{proof}
    The inequality $d_{G_{n+m}}(x \cdot y,x \cdot z) \leq d_{G_{m}}(y,z)$ follows from the self-similarity of the graphs, Lemma \ref{lemma:CharEdges}. We are left to show the opposite inequality. By induction, we may assume $n = 1$ and $x \in V$.

    We first check the symmetric hypercubic case. The dimension parameter is denoted $d$.
    We consider the following folding map. For every $w \in V$ we define the composition 
    \[
        \alpha_w := \prod_{i = 1}^d \alpha_{i}^{\abs{c_i(x) - c_i(w)}}
    \]
    of the of horizontal symmetries $\alpha_i$ for $i \in \mathcal{T}$.
    Note that the horizontal symmetries for different types commute with each other so the above definition is well-defined.
    
    Set $f(w \cdot y) := x \cdot \alpha_w(y)$. This is a simplicial map $G \to G$ onto $x \cdot V^m$ satisfying $f|_{x \cdot V^m} = \id_{x \cdot V^m}$. To see that it is simplicial, let $\{y,z\} \in E_{n+1}$. If $y,z \in w \cdot V^m$ for the same $w$, then $\{y,z\}$ is mapped to an edge because horizontal symmetries are graph isomorphisms. If $y = w \cdot y'$ and $y = u \cdot z'$ for distinct $w,u \in V$, then $\alpha_w = \alpha_u \circ \alpha_i$ where $i$ is the type of $\{x,w\}$. It then follows from Lemma \ref{lemma:CharEdges} that $\alpha_w(y') = \alpha_u(z')$, meaning $f(y) = f(z)$. Thus, the converse inequality follows because $f$ maps any path between points in $x \cdot V^m$ to path between the same vertices fully contained in $x \cdot V^m$ whose length is less than or equal to the original length. Here, we use that the induced graph of $x\cdot V^m$ is isomorphic to $G_m$, by Lemma \ref{lemma:CharEdges}.

    Next we deal with the pentagonal carpet.
    We again assume $n = 1$ and $x \in V$, and perform a similar folding argument, i.e. construct a simplicial map $f : G_{m+1} \to G_{m+1}$ onto $x \cdot V^m$ with $f|_{x \cdot V^m} = \id_{x \cdot V^m}$.
    By symmetry, we assume $x = 0$. Let $\alpha$ be the non-trivial symmetry of $G_{m+1}$ that fixes $2^{m+1}$. Then $\alpha$ satisfies the conditions in Proposition \ref{prop:ReflFold} and we can choose the partition so that $C \cup \Delta \subseteq 0 \cdot V^m \cup 1 \cdot V^m$. We first take the folding map $f_1$ onto the subset $C \cup \Delta$. Then if $\beta$ is the non-trivial symmetry of $G_{m+1}$ that fixes $3^{m+1}$, we see that $f(0 \cdot V^m \cup 1 \cdot V^m) = 0 \cdot V^m \cup 1 \cdot V^m$. By applying Proposition \ref{prop:ReflFold} to $\beta$ restricted to the subgraph induced by $0 \cdot V^m \cup 1 \cdot V^m$, we find a simplicial map $f_2$ onto $0 \cdot V^m$.
    We then take $f = f_2 \circ f_1$.

    The argument for the gasket is very similar to the one for the pentagon and employs twice applying Proposition \ref{prop:ReflFold} to symmetries fixing vertices $1^{m+1}$ and $2^{m+1}$.
\end{proof}

The following proposition yields the existence of an ``interior point'' to sets $G_m$ which is far from the ``boundary'' of $G_m$.
For $m \in \N$ we denote $\partial G_m \subseteq V^m$ the union of all the sets $I_{t,+}^m, \, I_{t,-}^m$ for all $t \in \mathcal{T}$. 

\begin{proposition}\label{prop:interior}
    Let $\Xi$ be a symmetric hypercubic IGS, the Sierpi\'nski gasket or the pentagonal carpet, and let $M_* \geq 1$ be the constant in Theorem \ref{thm:Inclusions}. Then there is $m \in \N$ and $v_* \in V^m$ so that $B(v_*,M_* + 1) \cap \partial G_m = \emptyset$.
\end{proposition}

\begin{proof}
    We consider the symmetric hypercubic case first. We use the functions $\varphi_{m,t}$ from the proof of Proposition \ref{prop: Distance between opposite faces}.
    First fix some $m \in \N$.
    Let $v,w \in V$ so that
    \[
        \begin{cases}
            c_t(v) = 1 & \text{ for all } t \in \mathcal{T}\\
            \xi(v) = 1
        \end{cases}
        \quad \text{ and } \quad 
        \begin{cases}
            c_t(w) = L_* & \text{ for all } t \in \mathcal{T}\\
            \xi(w) = 1.
        \end{cases}
    \]
    Then we choose the sequence $v_3,\dots,v_{m+2} \in V$ arbitrarily and set
    \[
      v_* := vwv_3\ldots v_{m+2} \in V^{m+2}.  
    \]
    Then for each $t \in \mathcal{T}$,
    \begin{align*}
        L_*^m & \leq (L_* - 1)\cdot L_*^{m} = (c_t(w) - 1)\cdot L_*^{m} \leq \varphi_{m+2,t}(v_*).
    \end{align*}
    Because $c_t(v) = 1$, we have
    \[
        \varphi_{m+2,t}(v_*) \leq \sum_{i = 2}^{m+2} (L_* - 1)L_*^{m+2-i} = L_*^{m+1} - 1.
    \]
    Because the maps $\varphi_{m,t}$ are simplicial and $\varphi_{m,t}(I_{t,\pm}^{m+2}) \in \{0,L_*^{m+2} - 1\}$,
    \[
        d_{G_{m+2}}(v_*,\partial G_{m+2}) \geq L_*^m,
    \]
    and we are therefore done by choosing $m \in \N$ so that $L_*^m > M_*$.

    Next we consider the gasket and the pentagonal carpet.
    A similar argument works for both so we omit the gasket case. First fix $m \geq 1$ and set $v_* := 02^m$. Consider the folding map $f$ onto $0 \cdot V^m$ that was defined during the proof of Proposition \ref{prop:MetricSS}. Note that $f$ satisfies
    \[
        f(\partial G_{m + 1}) = \partial G_{m + 1} \cap 0 \cdot V^m = 0 \cdot L_2^m \cup 0 \cdot L_3^m.
    \]
    Because $f|_{0 \cdot V^m} = \id_{0 \cdot V^m}$ and $f$ simplicial onto the subgraph induced by the subset $0 \cdot V^m$,
    \[
        d_{G_{m+1}}(v_*,\partial G_{m+1}) = d_{G_{m+1}}(v_*,0 \cdot L_2^m \cup 0 \cdot L_3^m).
    \]
    By combining Propositions \ref{prop: Distance between opposite faces} and \ref{prop:MetricSS},
    \[
        d_{G_{m+1}}(v_*,0 \cdot L_2^m \cup 0 \cdot L_3^m) = d_{G_m}(2^m,L_2^m  \cup L_3^m) \geq 1/2L_*^m.
    \]
    We are done by combining the previous two equations.
\end{proof}

Finally, we need one more technical condition for symmetric hypercubic IGSs and the pentagonal carpet.

\begin{proposition}\label{prop:ManyPaths}
    Let $\Xi$ be a symmetric hypercubic IGSs or the pentagonal carpet. Then for all $m \in \N$ and $t \in \mathcal{T}$ the graph $G_m$ contains $2^m$ pair-wise disjoint paths from $I_{t,+}^m$ to $I_{t,-}^m$.
\end{proposition}

\begin{proof}
    We argue the hypercubic case first. By symmetry we assume $ t = 1 \in \mathcal{T}$. Note that we essentially constructed one path during the proof of Lemma \ref{lem: cornerdist} by constructing a ``line segment'' between ``corner points''.    
    The construction is quite intuitive by examining Figure \ref{fig:SC}. Let $\theta_1$ be that path on $G_1$.
    Then we take another similar line segment $\theta_2 \in G_1$ between two corner points. For instance it can be obtained by applying to $\theta_2$ the horizontal symmetry $\alpha_2$ for the type $2 \in \mathcal{T}$. These are clearly disjoint paths. This gives the base case of the induction. Supposing now that we have constructed $2^m$ paths in $G_m$, we can obtain $2^m$ paths in $G_{m+1}$ by gluing the self-similar copies of these in the subgraphs $v\cdot G_m$ for $v\in \theta_i$ for $i=1,2$. This way one can construct $2^{m+1}$ disjoint paths between $I_{1,+}^{m+1}$ to $I_{1,-}^{m+1}$. 

    Next we deal with the pentagonal carpet. By symmetry we  assume $\{ I_{t,+},I_{t,-} \} = \{L_2,L_4\}$.
    The idea is essentially the same as above but we make use of a symmetry trick. 
    We warmly recommend the reader to make heavy use of Figure \ref{fig:PC}.

    We prove by induction on $m$ that for all $x \in L_2^m$ there is a path connecting $x$ to $L_4^m$ so that these paths are disjoint for different $x,y \in L_2^m$.
    The base case $m = 1$ is obvious so we make an inductive hypothesis (IH) for $m \in \N$.
    By self-similarity, for every vertex in $x \in 0 \cdot L_2^m$ there is a path that connects it to some $y \in 1 \cdot L_2^m$ and the paths are pair-wise disjoint. By applying self-similarity and IH once again, we can find a path from $y$ to $z \in 1 \cdot L_4^m \subseteq L_4^{m+1}$. By connecting these two paths together, we obtain a disjoint family of $2^m$ paths from $0 \cdot L_2^m$ to $1 \cdot L_4^m$.
    It is clear by construction that these paths are contained in $0 \cdot V^m \cup 1 \cdot V^m$.

    Next we construct another $2^m$ pair-wise disjoint paths so that they are contained in $2 \cdot V^m \cup 3 \cdot V^m \cup 4 \cdot V^m$. Then we would be done because these do not intersect with the previously constructed ones.
    By applying a suitable symmetry of the regular pentagon and the IH, we have $2^m$ disjoint paths from $L_2^m$ to $L_0^m$.
    Thus, we have a path connecting $x \in 4 \cdot L_2^m$ to $y \in 3 \cdot L_2^m$. Then by applying the IH again, we can connect $y \in 3 \cdot L_2^m$ to some $z \in 2 \cdot L_1^m$. We use the IH one last time, and a suitable symmetry in the same way as in the beginning of the paragraph, to get a path from $z \in 2 \cdot L_1^m$ to $w \in 2 \cdot L_4^m$. We have found additional $2^m$ pair-wise disjoint paths. Because these are contained in $2 \cdot V^m \cup 3 \cdot V^m \cup 4 \cdot V^m$, and the previously constructed once are contained in $0 \cdot V^m \cup 1 \cdot V^m$, the union of these two families is a family of paths we were looking for.
\end{proof}

\begin{remark}
    In the previous proof, one may notice an interesting feature of the pentagonal carpet. The argument shows that the constructed family of paths is a partition of $V^m$, meaning not only that it is a disjoint family of paths but also that it covers the vertex set.
    We do not make any use of this here but it may be useful in future works.
\end{remark}

\section{Construction of the limit space}\label{sec:LS}
The rest of the work studies limit spaces produced by IGSs/replacement graphs. These are constructed by taking limits of rescaled path-metrics.
This choice is largely driven by a desire to give a simpler analysis and to avoid many technical issues. 
For instance, it simplifies the proof of approximate self-similarity.
Other metrics could also be considered, such as visual metrics \cite{BM} that correspond to an associated iterative subdivision of the constructed fractal; see also \cite{Kigamiweighted}.
Conceptually, the main requirements for the construction is for it to produce a quasivisual approximation in the sense of \cite{bonkmeyertrees}.
We do not need or use this language in the present paper since the rescaled path metrics offer a very convenient and simple framework for our results. More complex metrics may be needed in non-symmetric settings; compare for example with the well known example of the so called Rickman rug \cite{KleinerSchioppa}.

\subsection{Assumptions}
Let us begin the discussion on the limit space by introducing lists of assumptions.
We always explicitly state which conditions are assumed to hold at each stage.

The first one is needed to construct a self-similar metric space suitable for our goals.
We remark that the notation $B(v,r)$ where $r \geq 0$ and $v \in V^n$ for $n \in \N$ always refers to the open ball with respect to the path metric $d_{G_n}$.
We also recall the notation $z \cdot L$ from Subsection \ref{subsec:RR}.

\begin{assumption}\label{Assumption: For limit good limit space}
    Let $\Xi$ be an IGS whose replacement graphs have uniformly bounded degree.
    We assume that there are constants $C,L_*,M_* \geq 2$ and $c_1,c_2 > 0$ so that the following hold.
    \begin{enumerate}[label={\textup{(\textcolor{blue}{M\arabic*})}}, widest=a, leftmargin=*]
        \vspace{3pt}
        \item \label{item:A1}
        For all $n,m \in \N, \, v \in V^{n+m}$,
        \[
            B(v,c_1L_*^m) \subseteq \bigcup_{z \in B(\pi_{n+m,n}(v),M_*)} z \cdot V^m.
        \]
        \vspace{3pt}
        \item \label{item:A2}
        For all $n,m \in \N, \, v \in V^{n+m}$,
        \[
            \bigcup_{z \in B(\pi_{n+m,n}(v),M_*)} z \cdot V^m \subseteq B(v,c_2 L_*^m).
        \]
        \vspace{3pt}
        \item \label{item:InteriorPoint}
        There is $n_0 \in \N$ and $v_* \in V^n$ so that $B(v_*,M_* + 1) \cap \partial G_n = \emptyset$ for all $t \in \mathcal{T}$.
        Here $\partial G_n$ is the union of the
        sets $I_{t,+}^n,\,I_{t,-}^n$ for $t \in \mathcal{T}$.
        \vspace{3pt}
    \end{enumerate}
    \end{assumption}

The rest of the assumptions are needed in Section \ref{Section: CLP}. The following is used to ensure properties of discrete moduli that are necessary for our proofs.
Their validity, in general, is quite non-trivial and the use of symmetries of the construction seems to be the best general approach to establish them.
We also expect that the slit carpet in Figure \ref{fig: Slit carpet} fails to satisfy these assumptions.

\begin{assumption}\label{assumption:Symm}
    Let $\Xi$ be an IGS, and $\mathcal{A}$ be some collection of unordered pairs $\{L_1,L_2\}$ so that $L_i \in \{ I_{t_i,+}, I_{t_i,-} \}$ for some $t_i \in \mathcal{T}$ and $i = 1,2$. We assume $\Xi$ and $\mathcal{A}$ to satisfy the following conditions.
    \begin{enumerate}[label={\textup{(\textcolor{blue}{S\arabic*})}}, widest=a, leftmargin=*]
    \vspace{3pt}
        \item \label{item:S1}
        $L_1 \cap L_2 = \emptyset$ for all $\{ L_1,L_2 \} \in \mathcal{A}$.
        \vspace{3pt}
        \item \label{item:FlowAss}
        For all $p \in (1,\infty)$ there is a constant $C \geq 1$ so that for all $m \in \N$ there is a flow basis $\mathcal{B}$ on $G_m^e$ satisfying
        \[
            \mathcal{E}_q(\mathcal{B}) \leq C \cRes_p(L_1^m,L_2^m,G_m)
        \]
        for all $\{L_1,L_2\} \in \mathcal{A}$. 
        Here $q = p/(p-1)$.
        The constant $C$ is assumed to be depending on $p$ in a way that $C^{1-p}$ is continuous in $p \in (1,\infty)$ and has a limit when $p \to 1^+$.
        \vspace{3pt}
        \item \label{item:Folding}
        For all $n,m \in \N$ and a path $\theta = [v_1,\ldots v_k]$ in $G_{n+m}$ so that $v_k \notin B(v_1,M_*)$, there is $\{L_1,L_2\} \in \mathcal{A}$ so that the folded path $\fold_m(\theta)$ (as defined in Subsection \ref{subsec:Symmetries}) intersects with both $\pi_{n,m}(v_k) \cdot L_1^{n}$ and $\pi_{n,m}(v_k) \cdot L_2^{n}$.
    \end{enumerate}
\end{assumption}

A good choice for $\mathcal{A}$ in many cases seems to be to simply take all pairs $\{L_1,L_2\}$ so that $L_1 \cap L_2 = \emptyset$.

The next condition is essential for the combinatorial Loewner property. The Sierpi\'nski gasket fails to satisfy it.

\begin{assumption}\label{Assumption: No local cutpoints}
Let $\Xi$ be an IGS. We assume that for all $M > 0$ there are $m \in \N$ and $t \in \mathcal{T}$ so that the family of paths connecting $I_{t,+}^m$ to $I_{t,-}^m$ contains at least $M$ pair-wise disjoint paths.
\end{assumption}

\begin{proposition}\label{prop: Collecting assumptions}
    Assumptions \ref{Assumption: For limit good limit space} and \ref{assumption:Symm} are satisfied by the Sierpi\'nski gasket.
    Assumptions \ref{Assumption: For limit good limit space}, \ref{assumption:Symm} and \ref{Assumption: No local cutpoints} are satisfied by symmetric hypercubic IGSs and the pentagonal carpet.
\end{proposition}

\begin{proof}
    The conditions in Assumption \ref{Assumption: For limit good limit space} were verified in Theorem \ref{thm:Inclusions} and Proposition \ref{prop:interior}.

    To verify the conditions in Assumption \ref{assumption:Symm}, we define $\mathcal{A}$ to consists of all pairs $L_1,L_2$ so that $L_i \in \{ I_{t_i,+}, I_{t_i,-} \}$ for some $t_i \in \mathcal{T}$ and $i = 1,2$. Note that, for the considered examples, the pairs in $\mathcal{A}$ are exactly the pairs $\{ I_{t,+},I_{t,-} \}$ for some $t \in \mathcal{T}$.
    Thus, \ref{item:S1} is clear.
    The condition \ref{item:FlowAss} follows from Theorems \ref{thm: Flow bases for cub} and \ref{thm: Flow bases for pent} and Proposition \ref{prop: Flow bases gask}. The remaining condition \ref{item:Folding} in Assumption \ref{assumption:Symm} follows from the argument in the proof of Theorem \ref{thm:Inclusions}.

    Finally, Assumption \ref{Assumption: No local cutpoints} was verified in Proposition \ref{prop:ManyPaths}.
\end{proof}

\subsection{Geometry of the limit space}
We are now well-prepared to define the limit space.
Throughout the rest of the section, we consider a fixed IGS satisfying Assumption \ref{Assumption: For limit good limit space}.

We first define the \emph{symbol space} $\Sigma := V^{\N}$, and equip it with a metric $\delta$ and a Radon measure $\mathfrak{m}$. The points in $\Sigma$ are written with the sequence notation $\omega = \omega_1\omega_2\ldots \in \Sigma$.
We define the natural projections $\pi_n(\omega_1\omega_2\ldots) := \omega_1\omega_2\ldots \omega_n \in V^n$. For a given $v \in V^n$ we denote $\Sigma_v := \pi_n^{-1}(v) \subseteq \Sigma$. 

The metric $\delta$ is the \emph{word metric}
\[
    \delta( \omega,\tau ) :=
    \begin{cases}
    0 & \text{ if } \omega = \tau\\
    2^{-k} & \text{ if $\omega \neq \tau$ where } k = \min \{l : \pi_l(\omega) \neq \pi_l(\tau) \}.
    \end{cases}
\]
The Radon measure $\mathfrak{m}$ is defined as the measure which satisfies
\[
    \mathfrak{m}(\Sigma_v) := \abs{V}^{-n} \text{ for all } v \in V^n.
\]
In other words, $\mathfrak{m}$ is a Bernoulli measure with uniform weights. It follows from Caratheodory's extension theorem that the condition above uniquely determines a Radon measure on the metric space $(\Sigma,\delta)$.

\begin{definition}\label{def:LS}
    Let $d$ be the semi-metric on $\Sigma$ given by
    \[
        d(\omega,\tau) := \limsup_{n \to \infty} L_*^{-n}d_{G_n}(\pi_n(\omega),\pi_n(\tau)).
    \]
    Here $L_*$ is the same as in Assumption \ref{Assumption: For limit good limit space}.
    We define the  set $X$ as the quotient $X := \Sigma / \sim$ where $\omega \sim \tau$ if and only if $d(\omega,\tau)=0$.
    Note that $d(x,y)$ for $x,y \in X$ does not depend on the choices of the representatives in $\Sigma$ by the triangle inequality of path metrics.
    We denote the canonical projection $\chi : \Sigma \to X$.
    Lastly, we define the push-forward measure $\mu := \chi_*(\mathfrak{m})$ given by $\mu(A) := \mathfrak{m}(\chi^{-1}(A))$. The metric space $(X,d)$ is the \emph{limit space} of the IGS $\Xi$.
\end{definition}

The following proposition shows that our limit spaces are self-similar. 
For $v \in V^n$ and $\omega \in \Sigma$ we define the symbolic product $v \cdot \omega \in \Sigma$ by appending the sequence $\omega$ by $v$. This satisfies $\pi_{n+k}(v \cdot \omega) = v \cdot \pi_{k}(\omega)$ for all $k \in \N$.

\begin{proposition}\label{prop:Similarity}
    Let $n \in \N$ and $v \in V^n$.
    The mapping $\chi(\omega) \mapsto \chi(v \cdot \omega)$ is a well-defined $L_*^{-n}$-Lipschitz function $X \to X$.
\end{proposition}

\begin{proof}
    It follows from Lemma \ref{lemma:CharEdges},
    \begin{align*}
        d_{G_{n+m}}( \pi_{n+m}(v \cdot \omega), \pi_{n+m}(v \cdot \tau)) & = d_{G_{n+m}}( v \cdot \pi_{m}(\omega), v \cdot \pi_{m}(\tau))\\
        & \leq d_{G_m}(\pi_{m}(\omega),\pi_{m}(\tau)).
    \end{align*}
    The claim now follows from the definition of the metric.
\end{proof}

\begin{remark}
    If the mappings defined in Proposition \ref{prop:Similarity} are injective, then our construction would produce a self-similar structure in the sense of Kigami \cite{AnalOnFractals}. This for instance would hold if the conclusion of Proposition \ref{prop:MetricSS} is valid.
\end{remark}

We state the main goal of the section. Recall the terminology from Section \ref{Section: Preliminary}.

\begin{theorem}\label{thm:LS}
    Let $(X,d)$ be a metric space which can be realized as a limit space of an IGS $\Xi$ satisfying Assumption \ref{Assumption: For limit good limit space}, and let $c_1,c_2,L_*$ be as stated therein.
    \begin{enumerate}
        \vspace{3pt}
        \item \label{item:LS(1))} $(X,d)$ is compact, quasiconvex and $Q$-Ahlfors regular for 
        \[
            Q := \frac{\log(\abs{V})}{\log(L_*)}.
        \]
        \vspace{3pt}
        \item \label{item:LS(2))} $(X,d)$ is approximately self-similar.
        \vspace{3pt}
        \item \label{item:LS(3))} There exists $\alpha > 1$ so that for all $n \in \N$ and $v \in V^n$ the set $X_v := \chi(\Sigma_v)$ contains a point $z_v \in X_v$ satisfying the following properties.
        \begin{enumerate}
            \item $B(z_v,\alpha^{-1} L_*^{-n}) \subseteq X_v \subseteq B(z_v,\alpha L_*^{-n})$.
            \item If $v,w \in V^n$ are distinct,
            \[
                B(z_v,\alpha^{-1} L_*^{-n}) \cap B(z_w,\alpha^{-1} L_*^{-n}) = \emptyset.
            \]
        \end{enumerate}
        In particular, the coverings $\{X_v\}_{v \in V^n}$ form an $\alpha$-approximation of $X$.
    \end{enumerate}
\end{theorem}

Before proving this theorem, we establish first some auxiliary lemmas.

\begin{lemma}\label{lemma:WhenContained}
    Let $M_*$ be as in Assumption \ref{Assumption: For limit good limit space} and $x = \chi(\omega) \in X$.
    If $x \in X_v$ for $v \in V^n$, then $v \in B(\pi_n(\omega),M_*)$.
\end{lemma}

\begin{proof}
    This follows from \ref{item:A1}. Indeed, if $v \notin B(\pi_n(\omega),M_*)$, then by the definition of the metric $d$,
    \[
        d(x,y) \geq c_1 L_*^{-n} \text{ for all } y \in X_v.
    \]
\end{proof}

\begin{lemma}\label{lemma: Ball inside rooms}
    Let $c_1,c_2,M_*,L_*$ be as in Assumption \ref{Assumption: For limit good limit space} and $x := \chi(\omega) \in X$.
    Then for all $n \in \N$,
    \begin{equation*}
        B(x,c_1L_*^{-n}) \subseteq \bigcup_{w \in B(\pi_{n}(\omega) ,M_*)} X_w,
    \end{equation*}
    and also
    \begin{equation*}
         \bigcup_{w \in B(\pi_n(\omega),M_*)} X_w \subseteq B(x,c_2L_*^{-n}).
    \end{equation*}
\end{lemma}

\begin{proof}
    Let $w \in V^n \setminus B(\pi_n(\omega,M_*))$ and $y = \chi(\tau)$ so that $\tau \in \Sigma_w$. To prove the first inclusion, it is sufficient to show $y \notin B(x,c_1L_*^n)$.

    Now, by \ref{item:A1}, $d_{G_{n+m}}(\pi_{n+m}(\omega),\pi_{n+m}(\tau)) \geq c_1L_*^{m}$. By letting $m \to \infty$,
    \[
        d(x,y) = \limsup_{m \to \infty} L_*^{-(n+m)}d_{G_{n+m}}(\pi_{n+m}(\omega),\pi_{n+m}(\tau)) \geq c_1L_*^{-n}.
    \]

    The second inclusion follows from a similar argument using \ref{item:A2}.
\end{proof}

\begin{proof}[Proof of Theorem \ref{thm:LS}-\eqref{item:LS(1))}]
    Because $(\Sigma,\delta)$ is compact, the compactess of $(X,d)$ follows once we show that $\chi$ is continuous. First, by applying \ref{item:A2} for $n = 1$, we see that there is a constant $C$ satisfying $\diam(V^m,d_{G_m}) \leq C L_*^m$ for all $m \in \N$. By the self-similarity of the graphs, $\diam(v \cdot V^m,d_{G_{n+m}}) \leq CL_*^{m}$ for all $v \in V^n$ and $n \in \N$. This then implies
    \[
        \diam(\chi(\Sigma_v)) \leq CL_*^{-n},
    \]
    so the continuity of $\chi$ follows from the definition of the word metric $\delta$.

    We next show that $(X,d,\mu)$ is $Q$-Ahlfors regular. Note that Lemma \ref{lemma:WhenContained} implies for all $v \in V^n$,
    \[
        \chi^{-1}(X_v) \subseteq \bigcup_{w \in B(v,M_*)} \Sigma_w.
    \]
    We combine this with the first inclusion in Lemma \ref{lemma: Ball inside rooms} and the definition of $\mu,\, \mathfrak{m}$,
    \[
        \mu(B(x,c_1 L_*^{-n})) \leq \bigcup_{w \in B(\pi_n(\omega),2M_*)} \mathfrak{m}(\Sigma_w) \leq C \abs{V}^{-n} = C(L_*^{-n})^Q,
    \]
    where $C$ depends on $M_*$ and the maximum degree of the replacement graphs. Now for a given $r \leq 2\diam(X)$, we get the upper inequality in the definition of $Q$-Ahlfors regularity by taking $n \in \N$ with $CL_*^{-n-1} \leq r \leq C L_*^{-n}$ for suitable $C$. The lower inequality is obtained similarily by using the second inclusion in Lemma \ref{lemma: Ball inside rooms}. This gives $Q$-Ahlfors regularity.

    We now verify quasiconvexity. Let $x,y \in X$ be distinct points and write $x = \chi(\omega), \, y = \chi(\tau)$. We also fix a constant $C \geq 1$ which we determine during the proof. Let $n \in \N$ so that $CL_*^{-n - 1} \leq d(x,y) < CL_*^{-n}$. By a suitable choice of $C$, there is, according to Lemma \ref{lemma: Ball inside rooms}, another constant $K \geq 1$ depending on $C$ so that for all $m \in \N$,
    \[
        \pi_{n+m}(\omega),\, \pi_{n+m}(\tau) \in \bigcup_{w \in B(\pi_n(\omega),K)} w \cdot V^m.
    \]
    By the argument at the beginning of the proof, $d_{G_{n+m}}(\pi_{n+m}(\omega), \pi_{n+m}(\tau)) \leq C_1 L_*^m$ for all $m \in \N$.

    Next for all $m \in \N$ we take discrete paths $\theta_m = [u_{1,m},\ldots, u_{k_m + 1,m}]$ in $G_{n+m}$ connecting $\pi_{n+m}(\omega) $ to $ \pi_{n+m}(\tau)$ and having length at most $k_m \leq C_1 L_*^m$. We also fix arbitrary points $z_{i,m} \in X_{u_{i,m}}$. According to the second inclusion in Lemma \ref{lemma: Ball inside rooms}, $d(z_{i,m},z_{i+1,m}) \leq c_2 L_*^{-(n+m)}$. Thus we have
    \[
        \sum_{i = 1}^{k_m} d(z_{i,m},z_{i+1,m}) \leq C_1 c_2 L_*^{-n}.
    \]
    By letting $m \to \infty$, we obtain by Arzela--Ascoli theorem a continuous curve from $x$ to $y$ of length at most $C_1 c_2 L_*^{-n}$; see for instance \cite[Proof of Theorem 4.32]{bjorn2011nonlinear} for an analogous argument. By recalling how $n$ was defined, quasiconvexity follows.
\end{proof}

Next, we deal with approximate self-similarity. To this end, we need a new definition. First, we fix a sufficiently large constant $N \geq 1$ which is determined during the proof below. We let $\mathcal{N}$ to be a certain finite collection of balls $B(v,N)$ where $v \in V^n$ for some $n \in \N$.
The number $n$ can be different for different members.

The members of $\mathcal{N}$ are chosen so that the following condition holds. For every $n \in \N$ and $v \in V^n$ there is $B \in \mathcal{N}$ and a graph isomorphism between the subgraphs induced by $B$ and $B(v,N)$ that preserves types and orientations of the edges. Such isomorphism lifts to higher level graphs by Lemma \ref{lemma:CharEdges}; the explicit description is given in the proof below.
The set $\mathcal{N}$ can be chosen to be finite because the degrees of replacement graphs are assumed to be uniformly bounded.

\begin{proof}[Proof of Theorem \ref{thm:LS}-\eqref{item:LS(3))}]
    Let $x = \chi(\omega)$ and $r > 0$. Recall that the goal is to find a suitable biLipschitz embedding $F : B(x,r) \to X$.
    We may assume that $r$ is suitably small because otherwise we may simply choose $F$ to be the inclusion $B(x,r) \hookrightarrow X$.
    Fix a constant $C \geq 1$, whose value is determined later in the proof, and let $n \in \N$ with $CL_*^{-n-1} < r \leq CL_*^{-n}$.

    Now, let $B(w,N) \in \mathcal{N}$ for $w \in V^m$ so that there exists an isomorphism $f : B(\pi_n(\omega),N) \to B(w,N)$ as we discussed right before the proof. Then $f$ has natural lifts $f_k : \pi_{n+k,n}^{-1}(B(\pi_n(\omega),N)) \to \pi_{m+k,m}^{-1}(B(w,N))$ given by
    \[
        f_k(v \cdot z) := f(v) \cdot z \text{ for all } v \in B(\pi_n(\omega),N) \text{ and } z \in V^k.
    \]
    We claim that the restriction of $f_k$ onto ${\pi_{n+k,n}^{-1}(B(\pi_n(\omega),2M_*))}$ is an isometry onto the image with respect to the path metrics if we choose $N$ to be suitably large. To see this, we combine \ref{item:A1} and \ref{item:A2},
    \[
        \pi_{n+k,n}^{-1}(B(\pi_n(\omega),2M_*)) \subseteq B(z,C_1 L_*^k) \subseteq B(z,3C_1 L_*^k) \subseteq \pi_{n+k,n}^{-1}(B(\pi_n(\omega),N))
    \]
    where $z=\pi_{n+k}(\omega)$, $C_1$ depends only on the fixed constant $M_*$ and $N$ is sufficiently large. In particular, the shortest paths between vertices in $\pi_{n+k,n}^{-1}(B(\pi_n(\omega),2M_*))$ are contained in $\pi_{n+k,n}^{-1}(B(\pi_n(\omega),N))$. These paths are thus mapped by $f_k$ to the $ \pi_{n+k,n}^{-1}(B(w,N))$. This shows that $f_k$ does not increase distances. By performing the same analysis in reverse, we conclude that the restriction of $f_k$ to $\pi_{n+k,n}^{-1}(B(\pi_n(\omega),M_*))$ is an isometry.

    The final step is to take a limit of $f_k$. By the observation above,
    \[
    F : \bigcup_{v \in B(\pi_n(\omega),2M_*)} X_v \to \bigcup_{u \in B(w,2M_*)} X_u, \quad  F(\chi(v \cdot \tau)) := \chi(f(v) \cdot \tau),
    \]
    is a well-defined homemorphims.
    Since the discrete functions were isometric embeddings, $F$ satisfies
    \[
        \abs{F(y) - F(z)} = L_*^{n-m}d(y,z).
    \]
    Recall that $n \in \N$ was chosen so that $r \approx L_*^{-n}$. Also $m$ has an upper bound because the set $\mathcal{N}$ is finite. 
    We get
    \[
        L^{-1}\frac{d(y,z)}{r}\leq \abs{F(y) - F(z)} \leq L\frac{d(y,z)}{r}
    \]
    where $L$ depends on $m$ and the other constants used during the proof.
    
    To conclude the approximate self-similarity, what is left is to show that $F(B(x,r)) \subseteq X$ is an open set. To see this, notice that $F(B(x,r))=B(F(x),rL_*^{n-m})$. Indeed, this equality follows since 
    \[
        B(x,r) \subseteq \bigcup_{v \in B(\pi_{n}(\omega),M_*)} X_v,
    \]
    and
    \[
     B(F(x),r) \subset \bigcup_{u \in B(w,M_*)} X_u
    \]
    and $F$ is an isometry on a set containing $\bigcup_{v \in B(\pi_{n}(\omega),M_*)} X_v$ onto a set that contains $ \bigcup_{u \in B(w,M_*)} X_u$ (with the metric scaled by $L_*^{n-m}$). 
\end{proof}

We are left to check the properties of an $\alpha$-approximation.

\begin{proof}[Proof of Theorem \ref{thm:LS}-\eqref{item:LS(3))}]
    Let $n_0,\, v_*$ be as in \ref{item:InteriorPoint}. Fix $n \in \N$ and choose $v \in V^n$. Then we choose $z_v = \chi(\omega_v) \in X_{v \cdot v_*}$ for any point $\omega_v\in \Sigma_{v\cdot v^*}$.

    First, we claim that
    \begin{equation}\label{eq:noboundary}
        B(v \cdot v_*,M_* + 1) \subseteq v \cdot  (V^{n_0} \setminus \partial  G_{n_0}).
    \end{equation}
    This follows from the following observation. According to Lemma \ref{lemma:CharEdges}, any path in $G_{n + n_0}$ which exits $v \cdot V^{n_0}$ goes through $v \cdot \partial  G_{n_0}$.
    Then \ref{eq:noboundary} follows from \ref{item:InteriorPoint}.
    
    Then we combine \eqref{eq:noboundary} with the first inclusion of Lemma \ref{lemma: Ball inside rooms},
    \[
        B(z_v,c_1L_*^{-n}) \subseteq \bigcup_{w \in B(\pi_n(\omega_v),M_*)} X_w \subseteq X_v.
    \]
    Also the combination of \eqref{eq:noboundary} and Lemma \ref{lemma:WhenContained} yield
    \[
        B(z_v,c_1L_*^{-n}/2) \cap X_w = \emptyset
    \]
    for all $w \in V^n \setminus \{v\}$.
    Thus we have the desired result.
\end{proof}

We finish the section by establishing some properties of the natural incidence graphs associated to the closed coverings $\{ X_v \}_{v \in V^n}$ for $n \in \N$.
We denote these $\mathbb{G}_n := (V^n,\mathbb{E}_n)$ where we have identified the vertices $V^n$ with their associated sets $X_v$.
Our covering is not open, but if one wishes to have an open covering, we could replace $X_v$ with the interiors of
\[
    \bigcup_{w \in B(v,M_*)} X_w.
\]
This would result in adding some edges, and is only a matter of the preference. Both should work equally well because of Proposition \ref{prop: Modulus change levels}.

\begin{corollary}\label{lemma:DicsreteModuliSame}
    If $\mathbb{G}_n$ are as above then the identity map $\id_{V^n} : (V^n,d_{G_n}) \to (V^n,d_{\mathbb{G}_n})$ is a $2M_*$-biLipschitz simplicial map $G_n \to \mathbb{G}_n$. Moreover, there is a constant $C \geq 1$ so that for all $p \geq 1$, $n \in \N$ and any disjoint pair of non-empty subsets  $A,B \subseteq V^n$,
    \[
         C^{-1} \Mod_p(A,B,\mathbb{G}_n) \leq \Mod_p(A,B,G_n) \leq C \Mod_p(A,B,\mathbb{G}_n).
    \]
\end{corollary}

\begin{proof}
    To see that the identity is simplicial $G_n \to \mathbb{G}_n$, or in the other words $E_n \subseteq \mathbb{E}_n$, fix an edge $(v,w) \in E_n$. By Lemma \ref{lemma:CharEdges} there are $\omega,\tau \in \Sigma$ so that
    \[
    (\pi_m(v \cdot \omega),\pi_m(w \cdot \tau)) \in E_m \text{ for all } m \geq n.
    \]
    Thus $\chi(v \cdot \omega) = \chi(w \cdot \tau)$, meaning that $X_v \cap X_w \neq \emptyset$. By definition, $ \{v,w\} \in \mathbb{E}_n$, and $E_n \subseteq \mathbb{E}_n$ is now clear.

    Note that $d_{\mathbb{G}_n} \leq d_{G_n}$ is obvious from $E_n \subseteq \mathbb{E}_n$. On the other hand, if $\{v,w\} \in \mathbb{E}_n$, then $d_{G_n}(v,w) < 2M_*$. Indeed, according to Lemma \ref{lemma:WhenContained}, $\chi(\omega) \in X_v \cap X_w$ for some $\omega \in \Sigma$ implies
    \[
        d_{G_n}(v,w) \leq d_{G_n}(v,\pi_n(\omega)) + d_{G_n}(\pi_n(\omega),w) < 2M_*.
    \]
    We conclude $d_{G_n} \leq 2M_*d_{\mathbb{G}_n}$.

    Lastly, we deal with the modulus estimate.
    Note that the upper inequality holds for $C = 1$ because $E_n \subseteq \mathbb{E}_n$. For the lower one, we argue as follows.
    To prove the converse, let $\rho : V^n \to [0,\infty)$ be $\Theta(A,B,G_n)$-admissible. Define
    \[
        \rho'(v) := 2M_* \max_{d_{G_n}(v,w) \leq 2M_*} \rho(w).
    \]
    Because $d_{G_{m}} \leq 2M_* d_{\mathbb{G}_m}$, $\rho'$ is $\Theta(A,B,\mathbb{G}_n)$-admissible, and it has the $p$-mass estimate
    \[
        \cM_p(\rho') \leq 2M_* C^{2M_*} \cM_p(\rho)
    \]
    where $C \geq 1$ is a constant bounding the degrees of $G_n$.
\end{proof}

\section{Combinatorial Loewner property}\label{Section: CLP}
The final section of the paper regards the combinatorial Loewner property of the limit space, and the super-multiplicativity inequality of discrete modulus.

We recall the limit space $X := (X,d)$ of an IGS from Definition \ref{def:LS}, and also the notation related to discrete modulus from Subsection \ref{Subsec: Discrete modulus and conformal dimension}.
In particular, we recall the discrete moduli \eqref{eq:moddeltadef}, which are defined
\[
    \cM_{p,\delta}^{(n)} := \Mod_p(\Gamma_\delta,\mathbb{G}_n)
\]
where the $\alpha$-approximations $\mathbb{G}_n := (V^n,\mathbb{E}_n)$ are given by the incidence graphs of $\{X_v\}_{v \in V^n}$ and $\Gamma_\delta = \{\gamma : \diam(\gamma) \geq \delta\}$.

The goal of the rest of the paper is to prove the following two theorems. The first one regards the \emph{super-multiplicativity inequality}. For expository reason we include the sub-multiplicativity as well but note that it already follows from Bourdon--Kleiner (Proposition \ref{prop: Sub-mult}).

\begin{theorem}\label{thm:supermult}
    Let $(X,d)$ be a metric space which can be realized as a limit space of an IGS $\Xi$ satisfying Assumptions \ref{Assumption: For limit good limit space} and \ref{assumption:Symm}.
    There is a $\delta_0 > 0$ so that for any given $p \in [1,\infty)$ and $\delta \in (0,\delta_0)$ there exists a constant $C \geq 1$ so that the super-multiplicativity inequality,
    \[
       C^{-1}\cM_{p,\delta}^{(n)} \cdot \cM_{p,\delta}^{(m)}\leq  \cM_{p,\delta}^{(n+m)} \leq C \cM_{p,\delta}^{(n)} \cdot \cM_{p,\delta}^{(m)},
    \]
    holds for all $n,m \in \N$.
    The constant $C$ can be chosen to be continuous in $p\in [1,\infty)$.
    In particular, there is a constant $\cM_p \in (0,\infty)$,
    \[
       C^{-1}\cM_p^n \leq \cM_{p,\delta}^{(n)} \leq C\cM_p^n.
    \]
    The constant $\cM_p$ is continuous and non-increasing in $p\in [1,\infty)$.
\end{theorem}

\begin{remark}
    We expect in many cases that $p \mapsto \cM_p$ is strictly decreasing. However, this seems to be very difficult to prove in most cases.
\end{remark}

The second result is the combinatorial Loewner property.

\begin{theorem}\label{thm: CLP}
    Let $(X,d)$ be a metric space which can be realized as a limit space of an IGS $\Xi$ satisfying Assumptions \ref{Assumption: For limit good limit space}, \ref{assumption:Symm} and \ref{Assumption: No local cutpoints}. For a given $p \in [1,\infty)$ let $\cM_p \in (0,\infty)$ be as in Theorem \ref{thm:supermult}.
    Then there is a unique $Q > 1$ which solves $\cM_Q = 1$ and for which $(X,d)$ satisfies the Combinatorial $Q$-Loewner property. This value $Q$ is the (Ahlfors regular) conformal dimension of $(X,d)$.
\end{theorem}

\subsection{Super-multiplicativity}
Throughout the section, we consider a fixed IGS $\Xi$ satisfying Assumptions \ref{Assumption: For limit good limit space} and \ref{assumption:Symm}. The limit space is denoted $(X,d)$.
We recall $\mathcal{A}$ from Assumption \ref{assumption:Symm}.

\begin{lemma}\label{lemma:ModuliFaces}
    For all $p \in [1,\infty)$ there is a constant $C \geq 1$ so that for any pair $\{ L_1,L_2 \},\{L_3,L_4\} \in \mathcal{A}$ and $n \in \N$,
    \[
         C^{-1} \Mod_p\left( L_3^n,L_4^n,G_n \right) \leq \Mod_p\left( L_1^n,L_2^n,G_n \right) \leq C\Mod_p\left( L_3^n,L_4^n,G_n \right).
    \]
    The constant $C$ can be chosen so that it is continuous in $p \in [1,\infty)$.
\end{lemma}

\begin{proof}
    The case $p \in (1,\infty)$ follows from duality (Proposition \ref{prop: Duality}), \ref{item:FlowAss} and the trivial inequality
    \[
        \cRes_p(L_1^n,L_2^n,G_n) \leq \mathcal{E}_q(\mathcal{B})
    \]
    which follows from the definition of flow basis. Note that the constant $C$ in the claim can be chosen to have a limit when $p \to 1^+$. The claim for $p = 1$ then follows from the continuity of the optimization problem $\Mod_p$ because it has only finitely many variables and the target function is continuous both in the variables and the parameter $p$.
\end{proof}

By the previous lemma, we adopt the simpler notation
\[
    \cM_{p,n} := \max_{ \{L_1,L_2\} \in \mathcal{A}} \Mod_p\left( L_1^n,L_2^n,G_n \right). 
\]
While this is dangerously similar to $\cM_{p,\delta}^{(n)}$ we justify the choice in Proposition \ref{prop:preliModuli}.

\begin{lemma}\label{lemma:helpful}
    For all $p \in [1,\infty)$ there is a constant $C \geq 1$ so that for any $n,k \in \N$ and $v,w \in V^n$ with $d_{G_n}(v,w) = M_*$,
    \[
        C^{-1} \cM_{p,k} \leq \Mod_p(v \cdot V^k,w \cdot V^k,G_{n+k}) \leq  C \cM_{p,k}. 
    \]
    The constant $C$ can be chosen to be continuous in $p\in [1,\infty)$.
\end{lemma}

\begin{proof}
    We first prove the lower estimate using replacement flow. Take a shortest path $\theta = [v_1,\ldots,v_{M_* + 1}]$ from $v$ to $w$ in $G_n$ and take $\mathcal{J}$ to be the constant flow along this path, i.e., $\mathcal{J}(v_i,v_{i+1}) = 1$ and zero on all other edges. Its $p$-energy is equal to $M_*$. Then we use Theorem \ref{thm: Replacement flow} to construct a flow $\mathcal{J}'$ from $v \cdot V^k $ to $ w \cdot V^k$ with the energy estimate
    \[
        \mathcal{E}_q(\mathcal{J}') \leq C_1\mathcal{E}_q(\mathcal{J}) \mathcal{E}_q(\mathcal{B}) = C_1M_* \mathcal{E}_q(\mathcal{B}).
    \]
    The lower bound in the claim follows from the duality (Proposition \ref{prop: Duality}), \ref{item:FlowAss} and Lemma \ref{lemma:helpful}.

    For the upper bound, we use a suitable folding which will be provided by \ref{item:Folding}.
    First we take a density $\rho : V^k \to [0,\infty)$ which is $\Theta(L_1,L_2,G_k)$-admissible for every $\{L_1,L_2\} \in \mathcal{A}$. We can choose it so that its $p$-mass is comparable to $\cM_{p,k}$. To see why, since $\mathcal{A}$ is finite we can use Lemma \ref{lemma:helpful} to define $\rho$ as a vertex-wise maximum over all $\{L_1,L_2\} \in \mathcal{A}$ of $\Theta(L_1,L_2,G_k)$-admissible densities with $p$-mass comparable to $\cM_{p,k}$.
    We can also assume that $\rho \circ \alpha_t = \rho$ where $\{\alpha_t\}_{t \in \mathcal{T}}$ are the symmetries of the IGS as in Definition \ref{def:SymmetricIGS}.
    To see this, note that the size of the subgroup of bijections $V^k \to V^k$ generated by $\{ \alpha_t \}_{t \in \mathcal{T}}$ does not depend on $k$ because of the way we defined symmetries of IGSs (Definition \ref{def:symmetry}). The property $\rho \circ \alpha_t = \rho$ can be achieved by a symmetrization argument - at the possible expense of a constant.

    Now, we define a density $\widetilde{\rho}$ on $V^{n+k}$ by setting
    \[
        \widetilde{\rho}(u) :=
        \begin{cases}
            \rho(\sigma_{n+k,n}(u)) & \text{ if } d_{G_n}(\pi_{n+k,n}(u),v) \leq M_*\\
            0 & \text{ otherwise}.
        \end{cases}
    \]
    We claim that $\widetilde{\rho}$ is $\Mod_p(v \cdot V^k,w \cdot V^k,G_{n+k})$-admissible. Because the degrees of replacement graphs are uniformly bounded, the $p$-masses of $\rho$ and $\widetilde{\rho}$ are comparable. This would yield the desired claim.

    Let $\theta$ be any path from $v \cdot V^k $ to $ w \cdot V^k$ and consider the folding $\fold_k(\theta)$.
    By the construction of $\fold_k(\theta)$ (Subsection \ref{subsec:Symmetries}) and
    \ref{item:Folding}, $\fold_k(\theta)$ is a path contained in $w \cdot V^k$ and connects $w \cdot L_1^k$ to $w \cdot L_2^k$ for some pair $\{ L_1,L_2 \} \in \mathcal{A}$. Because $\rho$ is $\Theta(L_1^k,L_2^k,G_k)$-admissible, it follows from the definition of $\widetilde{\rho}$ and self-similarity (Lemma \ref{lemma:CharEdges}),
    \[
    1 \leq \sum_{v \in \fold_k(\theta)} \widetilde{\rho}(v).
    \]
    On the other hand, by the symmetry $\rho \circ \alpha_t = \rho$ and the construction of the folding,
    \[
        \sum_{v \in \fold_k(\theta)} \widetilde{\rho}(v) \leq \sum_{v \in \theta} \widetilde{\rho}(v) .
    \]
    The admissibility follows by combining the previous inequalities.
\end{proof}

The following proposition is a key ingredient in the proof of super-multiplicativity.
A similar method was used in \cite[Lemma 4.4]{BourK}.

\begin{proposition}\label{prop:preliModuli}
There is $\delta_0 > 0$ so that for any given $p \in [1,\infty)$ and $\delta \in (0,\delta_0)$, there is a constant $C \geq 1$,
\[
    C^{-1}\cM_{p,\delta}^{(n)} \leq \cM_{p,n} \leq C\cM_{p,\delta}^{(n)}.
\]
The constant $C$ can be chosen to be continuous in $p\in [1,\infty)$.
\end{proposition}

We remark that in the above statement we compare two different notions of discrete moduli; one is defined for curves in the limit space and the other on the paths on the replacement graphs.

\begin{proof}[Proof of Proposition \ref{prop:preliModuli}]
    We first perform some simple estimates to $\cM_{p,\delta}^{(n)}$.
    Let $\alpha > 1$ and the points $z_v \in X$ be as in Theorem \ref{thm:LS}-\eqref{item:LS(3))}.
    We also fix a suitably small constant $\delta>0$ determined during the proof.

    Let $l_1 \in \N$ so that $\delta \leq \alpha L_*^{-l_1} < 4^{-1}\diam(X)$. This choice is possible when $\delta$ is sufficiently small. Then every continuous curve connecting $B(z_v,\alpha L_*^{-l_1})$ to the complement of $B(z_v,2\alpha L_*^{-l_1})$ is contained in $\Gamma_\delta$.
    By the monotonicity of modulus,
    \begin{equation}\label{eq:Annuli1}
        \Mod_p(B(z_v,\alpha L_*^{-l_1}),X\setminus B(z_v,2\alpha L_*^{-l_1}),\mathbb{G}_n) \leq \cM_{p,\delta}^{(n)}.
    \end{equation}
    We also have a converse estimate. Let $l_2 \in \N$ so that $2^{-1}\alpha L_*^{-l_2} < \delta$. Because the balls $B(z_v,\alpha L_*^{-l_2})$ cover $X$, every curve in $\Gamma_\delta$ connects $B(z_v,\alpha L_*^{-l_2})$ to the complement of $B(z_v,2\alpha L_*^{-l_2})$ for some $v \in V^n$. This implies
    \begin{equation}\label{eq:Annuli2}
        \cM_{p,\delta}^{(n)} \leq \abs{V}^{l_2} \max_{v \in V^{l_2}} \Mod_p(B(z_v,\alpha L_*^{-l_2}),X\setminus B(z_v,2\alpha L_*^{-l_2}),\mathbb{G}_n).
    \end{equation}
    To see why, if $\rho_v$ is admissible for the modulus problem in the right hand side for $v \in V^k$, then we get a $\Gamma_\delta$-admissble density by taking a maximum of $\rho_v$ over all $v \in V^k$.
    In conclusion, we derive the desired estimate by finding suitable moduli estimates over annuli.

    To this end, fix a radius $r > 0$ and a point $x \in X$. We assume that $\alpha L_*^{-j-1} \leq r < \alpha L_*^{-j}$ for $j \in \N$.
    Our goal is to find certain estimates for
    \[
    \Mod_p(B(x,r),X\setminus B(x,2r),\mathbb{G}_n)
    \]
    which are allowed to depend on $j$.
    By Lemmas \ref{lemma:WhenContained} and \ref{lemma: Ball inside rooms} there are $i \in \N$ and subsets $U_1 \subseteq U_2 \subseteq V^{i+j}$ such that
    \[
        B(x,r) \subseteq \bigcup_{u \in U_1} X_u \subseteq \bigcup_{u \in U_2} X_u \subseteq B(x,2r).
    \]
    By choosing $i$ to be suitably large, it always holds that a continuous curve $\gamma$ connecting $B(x,r)$ to the complement of $B(x,2r)$,
    $\gamma$ intersects with $X_u$ and $X_v$ for some pair $u \in U_1$ and $v \in U_2$ with $d_{G_{i+j}}(u,v) = M_*$.
    Because $r \approx L_*^{-j}$, the size of $U_1,U_2$ can be bounded by constants depending on $i$, and $i$ can be bounded by the geometry of $X$. 
    By the same maximum argument as above, we thus have
    \[
        \Mod_p(B(x,r),X\setminus B(x,2r),\mathbb{G}_{i+j+n}) \leq \abs{U_1}\abs{U_2} \max_{u \in U_1} \max_{ v  } \Mod_p( u \cdot V^n, v \cdot V^n, \mathbb{G}_{i+j+n} )
    \]
    where the second maximum is taken over $v \in U_2$ with $d_{G_{i+j}}(u,v) = M_*$.
    By combining the previous inequality with Lemma \ref{lemma:helpful} and Lemma \ref{lemma:DicsreteModuliSame},
    \[
        \Mod_p(B(x,r),X\setminus B(x,2r),\mathbb{G}_{i+j+n}) \leq C_1 \cM_p^{(n)}.
    \]
    We then apply Proposition \ref{prop: Modulus change levels} to obtain an estimate, which depends on $i$ and $j$,
    \begin{equation}\label{eq:ForCLP2}
         \Mod_p(B(x,r),X\setminus B(x,2r),\mathbb{G}_{n}) \leq C_2\Mod_p(B(x,r),X\setminus B(x,2r),\mathbb{G}_{i+j+n}).
    \end{equation}
    The combination of the previous two inequalities with \eqref{eq:Annuli2} for $l_2 = j$ concludes the lower bound of the desired inequality.

    We move on to the converse inequality. Again, 
    we fix a radius $r > 0$ and a point $x \in X$ so that $\alpha L_*^{-j-1} \leq r < \alpha L_*^{-j}$ for $j \in \N$.
    This time we are looking for lower bounds of
    \[
    \Mod_p(B(x,r),X\setminus B(x,2r),\mathbb{G}_n)
    \]
    that are allowed to depend on $j$.
    Fix $i \in \N$ and $u,v \in V^{i+j}$ so that
    \[
        X_u \subseteq B(x,r) \text{ and } X_v \subseteq X \setminus B(x,2r) 
    \]
    and $d_{G_{i+j}}(u,v) \leq K$ where $K$ depends on $i$. Such choices are possible as long as $i,j$ are sufficiently large according to Lemmas \ref{lemma:WhenContained}, \ref{lemma: Ball inside rooms}.

    Now, let $\rho : V^{i+j+n} \to [0,\infty)$ be $\Gamma(B(x,r),X\setminus B(x,2r))$-admissible.
    By quasiconvexity of $X$ (Theorem \ref{thm:LS}-\eqref{item:LS(2))}) and self-similarity (Proposition \ref{prop:Similarity}) every discrete path $\theta$ in $\mathbb{G}_{i+j+n}$ can be used to construct a continuous curve $\gamma$ in the limit space $X$.
    By using Lemma \ref{lemma:WhenContained}, we can choose $\gamma$ so that
    \[
        V^{i+j+n}[\gamma] \subseteq \bigcup_{v \in \theta} B(v,M_*)
    \]
    where the balls are with respect to the path metric of $G_{i+j+n}$.
    By a similar argument as in Lemma \ref{lemma:DicsreteModuliSame}, there is a constant $C_3$ such that the density
    \[
        \widetilde{\rho}(w) := C_3\max_{w'} \rho(w'),
    \]
    where the maximum is taken over all $w \in V^{i+j+n}$ with $d_{G_{i+j+n}}(w,w') \leq M_*$, is  
    $\Theta(u \cdot V^{n}, v \cdot V^n,\mathbb{G}_{i+j+n})$-admissible density with $p$-mass is comparable to $\rho$.
    This shows
    \begin{align*}
         \Mod_p(u \cdot V^{n}, v \cdot V^n,\mathbb{G}_{i+j+n}) & \leq  C_4\Mod_p(B(x,r),X\setminus B(x,2r),\mathbb{G}_{i+j+n})\\
         & \leq C_5\Mod_p(B(x,r),X\setminus B(x,2r),\mathbb{G}_{n})
    \end{align*}
    where the last inequality follows from Proposition \ref{prop: Modulus change levels}.

    By using the same replacement flow argument and Lemma \ref{lemma:DicsreteModuliSame}, we get a lower bound
    \[
        \Mod_p(u \cdot V^{n}, v \cdot V^n,\mathbb{G}_{i+j+n}) \geq C_6^{-1}\cM_{p,n}
    \]
    where $C_6$ only depends on the distance upper bound $d_{G_{i+j}}(u,v) \leq K$.
    By combining the previous two inequalities with \eqref{eq:Annuli1} for $j = l_1$, we obtain the second inequality that we aimed to show.
\end{proof}

\begin{proof}[Proof of Theorem \ref{thm:supermult}]
    We already have the sub-multiplicativity inequality by the result of Bourdon--Kleiner (Proposition \ref{prop: Sub-mult}), and Theorem \ref{thm:LS}. By the replacement flow argument (Theorem \ref{thm:ReplacementflowSuper}) combined with Proposition \ref{prop:preliModuli}, Lemma \ref{lemma:ModuliFaces} and \ref{item:FlowAss}, we get the super-multiplicativity inequality. 
    In conclusion, for all $p \in [1,\infty)$ and $n,m \in \N$,
    \[
       C^{-1}\cM_{p,\delta}^{(n)} \cdot \cM_{p,\delta}^{(m)}\leq  \cM_{p,\delta}^{(n+m)} \leq C \cM_{p,\delta}^{(n)} \cdot \cM_{p,\delta}^{(m)}.
    \]
    The continuity of $C$ in $p$ also follows.
    We remark that the case $p=1$ is obtained by a similar continuity argument as in Lemma \ref{lemma:ModuliFaces}.

    It is well-known that a sequence satisfying both sub- and super-multiplicativity is comparable to an exponential; see e.g. \cite[Theorem 5.1]{BarlowBassResistanceOnSC}. We nevertheless sketch the argument.
    Let $a_n := \log(C\cM_{p,\delta}^{(n)})$. Then $a_n$ is sub-additive, meaning $a_{n+m} \leq a_n + a_m$. It follows from Fekete's subadditive lemma that
    \[
        \rho_1 := \lim_{n \to \infty} \frac{a_n}{n} = \inf_{n > 0} \frac{a_n}{n}
    \]
    exists.
    Also $b_n := \log(C^{-1}\cM_p^{(n)})$ is super-additive, so
    \[
       \rho_2 := \lim_{n \to \infty} \frac{b_n}{n} = \sup_{n > 0} \frac{b_n}{n}
    \]
    exists as-well, and we clearly have $\rho_1 = \rho = \rho_2$.
    But then
    \[
        \rho - \frac{\log(C)}{n} \leq \frac{\log(\cM_p^{(n)})}{n} \leq \rho - \frac{\log(C^{-1})}{n}.
    \]
    By taking exponentials,
    \begin{equation}\label{eq:EXP}
        C^{-1}(e^\rho)^n \leq \cM_{p,\delta}^{(n)} \leq C(e^\rho)^n,
    \end{equation}
    and are done by choosing $\cM_p := e^\rho$.

    We show that $p \mapsto \cM_p$ is continuous and non-increasing.
    First, to see that $\cM_p$ is non-increasing, assume it is not. Namely, $\cM_{p} < \cM_{p'}$ for some $p < p'$. By \eqref{eq:EXP}, for large enough $n$ we have $\cM_{p,\delta}^{(n)} < \cM_{p',\delta}^{(n)}$. But this contradicts the fact that the discrete modulus is non-increasing.
    Second, to see that $p \mapsto \cM_p$ is continuous, again assume not. This means that there is $p \in [1,\infty)$ and $\varepsilon > 0$ so that $\cM_{p} > \cM_{p'} + \varepsilon$ for all $p' > p$ or $\cM_{p} < \cM_{p'} - \varepsilon$ for all $1 \leq p' < p$.
    Because the constant $C$ in \eqref{eq:EXP} is continuous in $p$, we can choose it to be uniform on any bounded interval. If either of the discontinuities would hold, by \eqref{eq:EXP}, for large enough $n \in \N$ the discrete modulus $\cM_{p,\delta}^{(n)}$ would not be continuous.
\end{proof}

\subsection{Combinatorial Loewner Property}
The remainder of the work establishes the combinatorial Loewner property of the limit space. In addition to Assumption \ref{Assumption: For limit good limit space}, we also assume Assumption \ref{Assumption: No local cutpoints}. This is necessary because, for instance, Sierpi\'nski gasket is not Combinatorially Loewner due to its local cut points; see \cite[Proposition 2.5]{BourK}.

\begin{lemma}
    For a given $p \in [1,\infty)$ let $\cM_p$ be as in Theorem \ref{thm:supermult}. Then $\cM_1  > 1$, and there is a $Q > 1$ such that $\cM_Q = 1$. 
\end{lemma}

\begin{proof}
    It follows from Assumption \ref{Assumption: No local cutpoints} and Lemma \ref{lemma:ModuliFaces} that for any $t \in \mathcal{T}$,
    \[
        \lim_{n \to \infty} \Mod_1( I_{t,+}^m,I_{t,-}^m,G_m) = \infty.
    \]
    Then $\cM_1 > 1$ follows from Proposition \ref{prop:preliModuli}, Lemma \ref{lemma:helpful} and Theorem \ref{thm:supermult}.
    To show that $\cM_Q = 1$ for some $Q > 1$, since $\cM_p$ is continuous in $p$ by Theorem \ref{thm:supermult}, it is sufficient to show that $\cM_p < 1$ for some $p > 1$.
    We use Theorem \ref{thm:LS} where we showed that $(X,d)$ is $Q'$-Ahlfors regular for $Q' := \log(\abs{V})/\log(\abs{L_*})$.
    This already implies that $\cM_p < 1$ for all $p > Q'$. We sketch the argument here. Consider any $\alpha$-approximation given by the coverings $\{X_v\}_{v \in V^n}$ of any $Q'$-Ahlfors regular compact metric space $(X,d)$.
    Then $\rho(v) := \diam(X_v)/\delta$ is $\Gamma_\delta$ admissible. Because $X_v$ are roughly balls with radii comparable to $L_*^{-n}$ and the centers are $L_*^{-n}$ separated, a volume counting argument shows that $\abs{V^n} \approx (L_*^n)^Q$.
    Then the $p$-mass of $\rho$ can be estimated
    \[
        \cM_p(\rho) \leq C (L_*^n)^Q \cdot (L_*^{-n})^p = C (L_*^n)^{Q-p}.
    \]
    Thus $\cM_{p,\delta}^{(n)} \to 0$ as $n \to \infty$ when $p > Q$.
\end{proof}

 Throughout the rest of the paper, we let $Q > 1$ be any value s.t. $\cM_Q = 1$. We have yet to show uniqueness of such $Q$ and that $Q$ is equal to the Ahlfors regular conformal dimension of $X$, but these will follow from the $Q$-CLP.
The remaining moduli computations are performed for $p = Q$. The dual exponent of $Q$ is denoted $q := Q/(Q-1)$.
We included some discussion for general $p$ in the last part of the paper.

Before verifying CLP we first establish certain moduli lower bounds of curve families that connect balls, namely the \emph{ball-Loewner estimate}.
It is well-known among the experts of Quasiconformal geometry that ball-Loewner estimates imply Loewner estimates bewteen any pair of continua. This follows from an iteration technique of Bonk--Kleiner \cite{BK05}. The variant of this estimate we are interested in comes from Bourdon--Kleiner \cite[Section 2.3]{BourK}.
The analysis on fractals community seems to refer to the ball-Loewner estimate as the \emph{Knight move condition};
see \cite{kigami,shimizu,KusuokaZhou,barlow1989construction}.
See also Murugan--Shimizu \cite{murugan2023first} for an analogical purely graph theoretical notion.

\begin{proposition}\label{prop: Combinatorial Ball-Loewner}
    Let $x,y \in X$, $r > 0$ and $A > 0$. Take two disjoint open balls $B_1 = B(x,r), B_2 = B(y,r)$ with $\dist(B_1,B_2) \leq A \cdot r$.
    Then there are constants $K,L > 0$ depending on $A$ satisfying the following condition.
    If $\Gamma$ is the family of continuous curves
    \begin{equation*}
        \Gamma := \{ \gamma \in \Gamma(B_1,B_2) : \diam(\gamma) \leq L \cdot r \},
    \end{equation*}
    and $n \in \Z$ so that $L_*^{-n} \leq r$, then for all $m \in \N$,
    \begin{equation*}\label{eq: Combinatorial Ball-Loewner}
        \Mod_Q(\Gamma, \mathbb{G}_{n+m}) \geq K.
    \end{equation*}
    The constants can be explicitly chosen $L := C(A+2)$ and $K := C(A+2)^{1-Q}$ for some constant $C$ depending on the geometry of $(X,d)$.
\end{proposition}

\begin{proof}
    Let $j \in \Z$ so that $L_*^{-j} \leq r < L_*^{-j + 1}$, $n \geq j$ and $m \in \N$.
    Note that $j$ has a lower bound because $(X,d)$ is compact.
    By the second inclusion in Lemma \ref{lemma: Ball inside rooms} there is $i \in \N$ with $i>-j$, depending only on the constants in Assumption \ref{Assumption: For limit good limit space}, and $u,v \in V^{i+j}$ so that
    \[
        x\in X_u \subseteq B_1 \text{ and } y\in X_v \subseteq B_2.
    \]
    Since the balls are disjoint, we obviously have $u \neq v$.
    
   We show that there exists a constant $C\geq 1$ s.t  $d_{G_{i+j}}(u,v) \leq C(A + 2)$. First iterate Lemma \ref{lemma: Ball inside rooms} to get
    \[
        B(\chi(w),(l \cdot c_1) L_*^{-i-j}) \subseteq \bigcup_{z \in B(\pi_{i+j}(w),l \cdot M_*)} X_z
    \]
    for all $w \in \Sigma$ and $l \in \N$. 
    Take $l \in \N$ to be the smallest integer satisfying $(l \cdot c_1) > (A+2){L_*^{i+2}}$. If we had $d_{G_{i+j}}(u,v) \leq l \cdot M_*$, then $d(X_u, X_v)> (A+2){L_*^{-j+1}}$.  This contradicts $d(x,y) \leq (A + 2)r \leq (A+2)L_*^{-j+1}$. Thus, with $C\geq l \cdot M_*$  we have $d_{G_{i+j}}(u,v) \leq C(A + 2)$.

    Next, take a shortest path $\theta := [v_1,\ldots v_l]$ from $u$ to $v$ in $G_{i+j}$ and take the unit flow along this path, meaning $\mathcal{J}(v_i,v_{i+1}) = 1$ and the flow vanishes on all other edges.
    Note that the $Q$-energy of $\mathcal{J}$ is equal to $\mathcal{E}_q(\mathcal{J}) = d_{G_{i+j}}(u,v)$.
    Then we apply replacement flow (Theorem \ref{thm: Replacement flow}) and the flow basis estimate \ref{item:FlowAss}, to construct a flow $\widetilde{\mathcal{J}}$ on $G_{n+m+i}$ satisfying the energy estimate
    \[
        \mathcal{E}_q(\widetilde{\mathcal{J}}) \leq C_2 \mathcal{E}_q(\mathcal{J}) \mathcal{E}_q(\mathcal{B}) \leq C_3 \mathcal{E}_q(\mathcal{J}) = C_3 d_{G_{i+j}}(u,v) \leq C_4 (A+2).
    \]
    In the second inequality we used the fact that $\mathcal{E}_q(\mathcal{B})$ has a constant upper bound because $\cM_Q = 1$.

    By construction of the replacement flow, the flow $\widetilde{\mathcal{J}}$ vanishes on edges that are not adjacent to the set
    \[
        Z_n := \bigcup_{w \in \theta} w \cdot V^{n-j + m}.
    \]
    In particular, if $\Theta$ is the family of paths
    \[
        \Theta := \{ \theta' \in \Theta(u \cdot V^{n-j+m},v \cdot V^{n-j+m},G_{n+m+i}) : \theta' \subseteq Z_n \}
    \]
    we have the modulus lower bound
    \[
        \Mod_Q(\Theta,G_{n+m+i}) \geq C_5 \mathcal{E}_q(\widetilde{\mathcal{J}})^{1-Q} \geq C_6(A+2)^{1-Q}.
    \]
    This follows by applying the duality (Proposition \ref{prop: Duality}) to the subgraph induced by the subset $Z_n \subseteq V^{n+m+i}$. By a similar argument as in the proof of Lemma \ref{lemma:DicsreteModuliSame}, we may replace $G_{n+m+i}$ by $\mathbb{G}_{n+m+i}$ in the previous inequality.

    The final step is to transfer the previous modulus lower bound for discrete paths to a modulus lower bound of continuous curves.
    Let $Z \subseteq X$ be the subset
    \[
        Z := \bigcup_{w \in \theta} X_w.
    \]
     The second inclusion in Lemma \ref{lemma: Ball inside rooms} and $\len(\theta) \leq C_1(A+2)$ imply that $\diam(Z) \leq C_1(A+2)L_*^{-j} \leq C_7(A+2)r$.
    Now, we choose $L := C(A+2)$ and let $\Gamma' \subseteq \Gamma$ be the family of continuous curves connecting $X_u$ to $X_v$ and contained in $Z$.
    It follows from the combination of the previous moduli estimate and a similar argument using quasiconvexity of $X$ as in the proof of Proposition \ref{prop:preliModuli},
    \[
        \Mod_Q(\Gamma',\mathbb{G}_{n+m+i}) \geq C_8 \Mod_Q(\Theta,\mathbb{G}_{n+m+i}).
    \]
    By combining this with the monotonicity of discrete modulus and Proposition \ref{prop: Modulus change levels}, we obtain the inequality
    \begin{align*}
        \Mod_Q(\Gamma,\mathbb{G}_{n+m+i}) & \geq \Mod_Q(\Gamma',\mathbb{G}_{n+m+i}) \geq C_{10} \Mod_Q(\Gamma',\mathbb{G}_{n+m})\\
        & \geq C_{10}C_8 \Mod_Q(\Theta,\mathbb{G}_{n+m})
        \geq C_{10}C_8 C_6 (A+2)^{1-Q}.
    \end{align*}
    The claim follows if we set $K := C(A+2)^{1-Q}$.
\end{proof}

We are now prepared to establish the Combinatorial Loewner property of the limit space $(X,d)$; recall definition \ref{def:CLP}.
The method is the same as in \cite[Theorem 5.2]{anttila2024constructions}

\begin{proof}[Proof of Theorem \ref{thm: CLP}]
    We need to verify the two conditions of CLP. The first one, \ref{CLP1}, follows from the ball-Loewner estimates (Proposition \ref{prop: Combinatorial Ball-Loewner}) and \cite[Proposition 2.9]{BourK}.
    We are left to deal with the second one \ref{CLP2}.

    Fix $x \in X$ and $L_*^{-n} \leq r$ for $n \in \N$.
    We first show the existence of a constant $K \geq 1$, independent of $n,m \in \N$ and $x \in X, \, r > 0$,
    \begin{equation}\label{eq:2Annuli}
        \Mod_Q(B(x,r),X \setminus B(x,2r),\mathbb{G}_{n+m}) \leq K.
    \end{equation}
    The idea is essentially the same as in the proof of \eqref{eq:ForCLP2} in Proposition \ref{prop:preliModuli}. The details are as follows.

    We may assume $r \leq 1$. If not, since the ambient space is compact, the cases $r > 1$ can be deduced by a covering argument from the cases $r \leq 1$.
    Let $k \in \N \cup \{0\}$ so that $L_*^{-n + k} \leq r < L_*^{-n + k + 1}$.
    By Lemmas \ref{lemma:WhenContained} and \ref{lemma: Ball inside rooms}, there are $i \in \N$ and subsets $U_1 \subseteq U_2 \subseteq V^{n-k+i}$ such that
    \[
        B(x,r) \subseteq \bigcup_{u \in U_1} X_u \subseteq \bigcup_{u \in U_2} X_u \subseteq B(x,2r).
    \]
    By choosing $i$ to be suitably large, the choices can be performed so that every continuous curve $\gamma$ connecting $B(x,r)$ to the complement of $B(x,2r)$ intersects with $X_u$ and $X_v$ for some pair $u \in U_1$ and $v \in U_2$ with $d_{G_{n-k+i}}(u,v) = M_*$.
    Thus, by the argument in the proof of Proposition \ref{prop:preliModuli},
    \begin{align*}
        & \quad \Mod_Q(B(x,r),X \setminus B(x,2r),\mathbb{G}_{n+m+i})\\
        \leq & \quad \abs{U_1}\abs{U_2} \max_{u \in U_1} \max_{v} \Mod_Q(u \cdot V^{m + k},v \cdot V^{m + k},\mathbb{G}_{n+m+i})\\
        \leq & \quad C\abs{U_1}\abs{U_2} \cM_Q^{m+k} \leq K,
    \end{align*}
    where the last inequality follows from $\cM_Q = 1$. We arrive to \eqref{eq:2Annuli} by applying Proposition \ref{prop: Modulus change levels}.

    Now we fix $C > 0$ and our objective is to estimate
    \[
        \Mod_Q(B(x,r),X \setminus B(x,(C+1)r),\mathbb{G}_{n+m})
    \]
    in terms of $C$.
    We consider the cases $C \leq 1$ and $C \geq 1$ separately. Let us assume the former. We take a covering $\{ B(y,Cr/4) \}_{y \in I}$ of $B(x,r)$ so that each $y \in B(x,r)$ and $d(y_1,y_2) \geq Cr/4$ for all distinct $y_1,y_2 \in I$. By using the $Q'$-Ahlfors regularity of $(X,d)$ (Theorem \ref{thm:LS}), we have a bound $\abs{I} \leq C_1C^{-Q'}$. This follows from a volume counting argument. Now for $y \in I$ we take
    \[
        \rho_y : V^{n+j} \to [0,\infty), \, \text{$\Gamma(B(y,Cr/4),B(y,Cr/2)$)-admissible},
    \]
    whose $Q$-mass can be bounded by $K$ according to \eqref{eq:2Annuli}.
    Then the vertex-wise maximum $\rho := \max_{y \in I} \rho_y$ is $\Gamma(B(x,r),X \setminus B(x,(C+1)r))$ admissible whose $Q$-mass can be bounded by $K \cdot\abs{I} \leq C_2 C^{-Q'}$. We can set $\psi(t) := C_2t^{Q'}$ when $t \geq 1$.

    Next we consider the case $C > 1$.
    Our method here is a quite standard technique in quasiconformal geometry; see e.g. \cite[Main Lemma 4.12]{HK} and \cite[Proposition 3.4]{BourK}.
    Let $N := \ceil{\log_2(C)}-1$, and for each $i = 1,\dots, N-1$ we take $\Gamma( B(x,2^{i}r), X \setminus B(x,2^{i+1}r))$-admissible density $\rho_i : V^{n+j} \to [0,\infty)$ with $Q$-mass bounded by $K$.
    Then we take the average $\rho := A/(N-1)\sum_{i = 1}^{N-1} \rho_i$ where $A$ is chosen so that every $X_v$ for $v \in V^{n+i}$ intersects at most $A$ different balls $B(x,2^i r)$. Such choice is possible by Theorem \ref{thm:LS}-\eqref{item:LS(3))}.

    Because every curve connecting $B(x,2^{i}r)$ to $B(x,2^{i+1}r)$ contains a sub-curve that is fully contained in the annuli $B(x,2^{i+1}r) \setminus B(x,2^{i}r)$, we can always choose $\rho_i$ so that $\rho_i(v) = 0$ unless $X_v \cap (B(x,2^{i+1}r) \setminus B(x,2^{i}r)) \neq \emptyset$. Because the members of the covering $\{X_v\}_{v \in V^{n+i}}$ have bounded overlapping (by \eqref{item:LS(1))} and \eqref{item:LS(3))} of Theorem \ref{thm:LS}), there is a constant $M$ depending only on the geometry of $(X,d)$ so that for every $v \in V^{n+i}$ the value $\rho_i(v)$ is non-zero for at most $M$ different $i = 1,\ldots l$. Using this, we get the $Q$-mass bound
    \[
        \cM_Q(\rho) \leq \frac{A M^{Q-1}}{(N-1)^Q} \sum_{i = 1}^{N-1}\cM_Q(\rho_i) \leq \frac{C_3}{N^{Q-1}} \leq C_4 \log(C)^{1-Q}
    \]
    Thus, when $t \leq 1$, we set $\psi(t) := C_4 \log(1/t)^{1-Q}$.

    We have now shown the $Q$-CLP  for any $Q$ s.t. $\cM_Q=1$. By \cite[Lemma 4.2]{EBConf} any such $Q$ must be equal to the conformal dimension of $X$ and is thus also unique.
\end{proof}

\section{Some further results}\label{sec: Some further results}
We finish this paper by discussing some further applications of our techniques to constructing Sobolev spaces, $p$-energies and diffusions on the limit spaces using the Kusuoka--Zhou method; we refer to \cite{murugan2023first,KusuokaZhou,kigami} for further discussion.

Let us assume that $\Xi$ is a fixed IGS satisfying Assumptions \ref{Assumption: For limit good limit space} and \ref{assumption:Symm}.
We denote the Hausdorff dimension
\[
    Q := \frac{\log(\abs{V})}{\log(L_*)}.
\]
We also define the \emph{$p$-walk dimension} by
\[
    d_{w,p} := \frac{\log(\abs{V} \cdot \cM_p^{-1})}{\log(L_*)},
\]
where $\cM_p$ is as in Theorem \ref{thm:supermult}.
We first verify a suitable volume growth condition for the graphs.

\begin{proposition}\label{prop: AR of replacement graphs}
    The replacement graphs $G_n$ are uniformly $Q$-Ahlfors regular, i.e. there is constant $C \geq 1$ so that for all $x \in V^n$ and $r \in [1,2\diam(G_n)]$,
    \begin{equation*}
       C^{-1}r^Q  \leq \abs{B(w,r)} \leq Cr^Q.
    \end{equation*}
\end{proposition}

\begin{proof}
    Follows from an identical argument as in the proof of $Q$-Ahlfors regularity of the limit space in Theorem \ref{thm:LS}-\eqref{item:LS(1))}.
\end{proof}

Next we have a natural upper estimate for the discrete moduli over annuli.

\begin{proposition}\label{prop: Discrete cap over annuli}
Let $p \in [1,\infty)$. There is a constant $C \geq 1$ so that for all $x \in V^n$ and $r > 0$
	\begin{equation}\label{eq: RG (Cap)}
		\Mod_p(B(x,r), V^n \setminus B(x,2r), \mathbb{G}_n) \leq C \frac{\abs{B(x,r)}}{r^{d_{w,p}}}.
	\end{equation}
\end{proposition}

\begin{proof}
    Similar arguments as in the proof of Theorem \ref{thm: CLP}, combined with Theorem \ref{thm:supermult}, show that
    \[
        \Mod_p(B(x,r), V^n \setminus B(x,2r),\mathbb{G}_n) \leq C \cM_p^m,
    \]
    where $r \approx L_*^{m}$. 
    The desired upper bound now follows from the definition of $d_{w,p}$ and Proposition \ref{prop: AR of replacement graphs}.
\end{proof}

The following proposition verifies the $p$-combinatorial Ball Loewner condition from \cite{murugan2023first}.

\begin{proposition}\label{prop: Discrete Ball-Loewner}
Let $p \in [1,\infty)$ and $A > 0$. There are constants $L,K > 0$ depending on $A$ so that for all $n \in \N$, $x,y \in V^n$ and $r > 0$ the following implication holds.
If $B_1 := B(x,r),\, B_2 := B(y,r) \subseteq V^n$ with $\dist(B_1,B_2) \leq A \cdot r$ and
	\[
		\Theta := \{ \theta \in \Theta(B_1,B_2) : \diam(\theta) \leq L \cdot r \}
	\]
	then
	\begin{equation}\label{eq: RG (BCL)}
		\Mod_p(\Theta,G_n) \geq K \cdot r^{Q - d_{w,p}}.
	\end{equation}
\end{proposition}

\begin{proof}
The idea is the same as in the proof of Proposition \ref{prop: Combinatorial Ball-Loewner}. First, we note that $r$ can be assumed to be suitably large. This is because otherwise we have only finitely many combinatorially distinct cases, depending only on $A$ and the maximum degree of the replacement graphs.
Then we choose $m \leq n$ according to the given radius $r > 0$. It is chosen as the smallest positive integer so that
\[
    \pi_{n,m}(x) \cdot V^{n-m} \subseteq B(x,r) \text{ and } \pi_{n,m}(y) \cdot V^{n-m} \subseteq B(y,r).
\]
We choose a shortest path $\theta$ from $\pi_{n,m}(x)$ to $\pi_{n,m}(y)$. Then take the unit flow along this path and apply the replacement flow (Theorem \ref{thm: Replacement flow}).
By the same argument as in the proof of Proposition \ref{prop: Combinatorial Ball-Loewner},
\[
    \Mod_p(\Theta,G_n) \geq C \cdot \cM_p^m \geq K r^{Q-d_{w,p}}.
\]
\end{proof}

Let $p \in (1,\infty)$.
According to Murugan--Shimizu \cite{murugan2023first}, the previous three propositions and the dimension condition $Q - d_{w,p} < 1$ are sufficient for the uniform elliptic Harnack inequality of $p$-harmonic functions and the uniform $(p,p)$-Poincar\'e-inequality on the replacement graphs $G_n$. By combining these estimates with the methods of Murugan--Shimizu, it should be possible to extend the constructions of Sobolev spaces to higher dimensional Menger sponges, and also to the broader class of generalized Sierpi\'nski carpets; see \cite[Section 7]{murugan2023first}.
In particular, this construction works when $p$ is equal to the conformal dimension of the ambient limit space. We also get a non-probabilistic construction of the self-similar Dirichlet form on the Menger sponge.
It is a difficult open problem whether the condition $Q - d_{w,p} < 1$ can be relaxed.

We expect the method of Murugan--Shimizu also applies on many other fractal spaces that we considered it the work, including the pillow space and the pentagonal carpet. A detailed study of this is left for future work.

Under the more restrictive dimensionality condition $d_{w,p} > Q$, it should be possible to verify the conductive homogeneity condition of Kigami \cite{kigami}; note that $d_{w,p} > Q$ is equivalent to $p$ being strictly larger than the conformal dimension of the ambient limit space. This is because the statement in Proposition \ref{prop: Discrete Ball-Loewner} is equivalent to the \emph{$p$-Knight move condition} of Kigami (see \cite{kigami,shimizu} and also \cite{barlow1989construction,barlow1999brownian,KusuokaZhou}) which is equivalent to the conductive homogeneity condition. We note that, for this equivalence, the condition $d_{w,p} > Q$ seems to be quite crucial.

\bibliographystyle{acm}
\bibliography{clp}

\begin{thebibliography}{10}

\bibitem{ACFPC}
{\sc Albin, N., Clemens, J., Fernando, N., and Poggi-Corradini, P.}
\newblock Blocking duality for {$p$}-modulus on networks and applications.
\newblock {\em Ann. Mat. Pura Appl. (4) 198}, 3 (2019), 973--999.

\bibitem{anttila2024constructions}
{\sc Anttila, R., and Eriksson-Bique, S.}
\newblock On constructions of fractal spaces using replacement and the
  combinatorial {L}oewner property.
\newblock {\em arXiv preprint arXiv:2406.08062\/} (2024).

\bibitem{barlow}
{\sc Barlow, M.~T.}
\newblock Diffusions on fractals.
\newblock In {\em Lectures on probability theory and statistics
  ({S}aint-{F}lour, 1995)}, vol.~1690 of {\em Lecture Notes in Math.} Springer,
  Berlin, 1998, pp.~1--121.

\bibitem{barlow1989construction}
{\sc Barlow, M.~T., and Bass, R.~F.}
\newblock The construction of brownian motion on the sierpinski carpet.
\newblock In {\em Annales de l'IHP Probabilit{\'e}s et Statistiques\/} (1989),
  vol.~25, pp.~225--257.

\bibitem{BarlowBassResistanceOnSC}
{\sc Barlow, M.~T., and Bass, R.~F.}
\newblock On the {R}esistance of the {S}ierpinski {C}arpet.
\newblock {\em Proceedings: Mathematical and Physical Sciences 431}, 1882
  (1990), 345--360.

\bibitem{barlow1999brownian}
{\sc Barlow, M.~T., and Bass, R.~F.}
\newblock Brownian motion and harmonic analysis on sierpinski carpets.
\newblock {\em Canadian Journal of Mathematics 51}, 4 (1999), 673--744.

\bibitem{UniquenessofBrownianMotionOnGSC}
{\sc Barlow, M.~T., Bass, R.~F., Kumagai, T., and Teplyaev, A.}
\newblock Uniqueness of {B}rownian motion on {S}ierpiński carpets.
\newblock {\em Journal of the European Mathematical Society 12}, 3 (2010),
  655--701.

\bibitem{bjorn2011nonlinear}
{\sc Bj{\"o}rn, A., and Bj{\"o}rn, J.}
\newblock {\em Nonlinear potential theory on metric spaces}, vol.~17.
\newblock European Mathematical Society, 2011.

\bibitem{BK01}
{\sc Bonk, M., and Kleiner, B.}
\newblock Quasisymmetric parametrizations of two-dimensional metric spheres.
\newblock {\em Inventiones mathematicae 150\/} (08 2001).

\bibitem{BK05}
{\sc Bonk, M., and Kleiner, B.}
\newblock Conformal dimension and gromov hyperbolic groups with 2-sphere
  boundary.
\newblock {\em Geometry \& Topology 9\/} (2005), 219--246.

\bibitem{BonkMerenkov}
{\sc Bonk, M., and Merenkov, S.}
\newblock Quasisymmetric rigidity of square {S}ierpi\'{n}ski carpets.
\newblock {\em Ann. of Math. (2) 177}, 2 (2013), 591--643.

\bibitem{BM}
{\sc Bonk, M., and Meyer, D.}
\newblock {\em Expanding Thurston Maps}, vol.~225 of {\em Mathematical Surveys
  and Monographs}.
\newblock American Mathematical Society, Providence, RI, 2017.

\bibitem{bonkmeyertrees}
{\sc Bonk, M., and Meyer, D.}
\newblock Quasiconformal and geodesic trees.
\newblock {\em Fund. Math. 250}, 3 (2020), 253--299.

\bibitem{BourK}
{\sc Bourdon, M., and Kleiner, B.}
\newblock Combinatorial modulus, the combinatorial {Loewner} property, and
  {Coxeter} groups.
\newblock {\em Groups Geom. Dyn. 7}, 1 (2013), 39--107.

\bibitem{BourdonPajot}
{\sc Bourdon, M., and Pajot, H.}
\newblock Cohomologie {$l_p$} et espaces de {B}esov.
\newblock {\em J. Reine Angew. Math. 558\/} (2003), 85--108.

\bibitem{bowers}
{\sc Bowers, P.~L., and Stephenson, K.}
\newblock A ``regular'' pentagonal tiling of the plane.
\newblock {\em Conform. Geom. Dyn. 1\/} (1997), 58--68.

\bibitem{Resistance4Ncarpets}
{\sc Canner, C., Hayes, C., Huang, R., Orwin, M., and Rogers, L.~G.}
\newblock Resistance scaling on {$4N$}-carpets.
\newblock {\em Forum Mathematicum 34}, 1 (2022), 61--75.

\bibitem{CFKP}
{\sc Cannon, J.~W., Floyd, W.~J., Kenyon, R., and Parry, W.~R.}
\newblock Constructing rational maps from subdivision rules.
\newblock {\em Conform. Geom. Dyn. 7\/} (2003), 76--102.

\bibitem{CFP}
{\sc Cannon, J.~W., Floyd, W.~J., and Parry, W.~R.}
\newblock Finite subdivision rules.
\newblock {\em Conform. Geom. Dyn. 5\/} (2001), 153--196.

\bibitem{cao2024whether}
{\sc Cao, S., and Chen, Z.-Q.}
\newblock Whether $ p $-conductive homogeneity holds depends on $ p$.
\newblock {\em Journal of Fractal Geometry\/} (2024).

\bibitem{Carrasco}
{\sc Carrasco~Piaggio, M.}
\newblock On the conformal gauge of a compact metric space.
\newblock {\em Ann. Sci. \'{E}c. Norm. Sup\'{e}r. (4) 46}, 3 (2013), 495--548.

\bibitem{chilakamarri2013self}
{\sc Chilakamarri, K.~B., Khan, M., Larson, C., and Tymczak, C.}
\newblock Self-similar graphs.
\newblock {\em arXiv preprint arXiv:1310.2268\/} (2013).

\bibitem{clais}
{\sc Clais, A.}
\newblock Combinatorial modulus on boundary of right-angled hyperbolic
  buildings.
\newblock {\em Anal. Geom. Metr. Spaces 4}, 1 (2016), 1--531.

\bibitem{david2024analytically}
{\sc David, G.~C., and Eriksson-Bique, S.}
\newblock Analytically one-dimensional planes and the combinatorial loewner
  property.
\newblock {\em arXiv preprint arXiv:2408.17279\/} (2024).

\bibitem{David_Schul_2017}
{\sc David, G.~C., and Schul, R.}
\newblock The analyst’s traveling salesman theorem in graph inverse limits.
\newblock {\em Annales Fennici Mathematici 42}, 2 (2017), 649–692.

\bibitem{EBConf}
{\sc Eriksson-Bique, S.}
\newblock Equality of different definitions of conformal dimension for
  quasiself-similar and {CLP} spaces.
\newblock {\em Ann. Fenn. Math. 49}, 2 (2024), 405--436.

\bibitem{godsil}
{\sc Godsil, C.~D., and McKay, B.~D.}
\newblock A new graph product and its spectrum.
\newblock {\em Bull. Austral. Math. Soc. 18}, 1 (1978), 21--28.

\bibitem{He}
{\sc Heinonen, J.}
\newblock {\em \normalfont ``{L}ectures on analysis on metric spaces''}.
\newblock Universitext. Springer-Verlag, New York, 2001.

\bibitem{HK}
{\sc Heinonen, J., and Koskela, P.}
\newblock Quasiconformal maps in metric spaces with controlled geometry.
\newblock {\em Acta Math. 181}, 1 (1998), 1--61.

\bibitem{MN}
{\sc Kajino, N., and Murugan, M.}
\newblock On the conformal walk dimension: quasisymmetric uniformization for
  symmetric diffusions.
\newblock {\em Invent. Math. 231}, 1 (2023), 263--405.

\bibitem{AnalOnFractals}
{\sc Kigami, J.}
\newblock {\em Analysis on Fractals}.
\newblock Cambridge Tracts in Mathematics. Cambridge University Press,
  Cambridge, 2001.

\bibitem{Kigamiweighted}
{\sc Kigami, J.}
\newblock {\em Geometry and analysis of metric spaces via weighted partitions},
  vol.~2265 of {\em Lecture Notes in Mathematics}.
\newblock Springer, Cham, [2020] \copyright 2020.

\bibitem{kigami}
{\sc Kigami, J.}
\newblock {\em Conductive homogeneity of compact metric spaces and construction
  of {$p$}-energy}, vol.~5 of {\em Memoirs of the European Mathematical
  Society}.
\newblock European Mathematical Society (EMS), Berlin, [2023] \copyright 2023.

\bibitem{KleinerICM}
{\sc Kleiner, B.}
\newblock The asymptotic geometry of negatively curved spaces: Uniformization,
  geometrization and rigidity.
\newblock In {\em International Congress of Mathematicians. Vol. II}. European
  Mathematical Society, Z{\"u}rich, 2006, pp.~743--768.

\bibitem{KleinerSchioppa}
{\sc Kleiner, B., and Schioppa, A.}
\newblock P{I} spaces with analytic dimension 1 and arbitrary topological
  dimension.
\newblock {\em Indiana Univ. Math. J. 66}, 2 (2017), 495--546.

\bibitem{kron2002green}
{\sc Kr{\"o}n, B.}
\newblock Green functions on self-similar graphs and bounds for the spectrum of
  the laplacian.
\newblock {\em Annales de l'institut Fourier 52}, 6 (2002), 1875--1900.

\bibitem{KusuokaZhou}
{\sc Kusuoka, S., and Zhou, X.~Y.}
\newblock Dirichlet forms on fractals: {P}oincar\'{e} constant and resistance.
\newblock {\em Probab. Theory Related Fields 93}, 2 (1992), 169--196.

\bibitem{Kwapisz}
{\sc Kwapisz, J.}
\newblock Conformal dimension via {\tt p}-resistance: {S}ierpi\'{n}ski carpet.
\newblock {\em Ann. Acad. Sci. Fenn. Math. 45}, 1 (2020), 3--51.

\bibitem{Laakso}
{\sc Laakso, T.~J.}
\newblock Plane with {$A_\infty$}-weighted metric not bi-{L}ipschitz embeddable
  to {${\Bbb R}^N$}.
\newblock {\em Bull. London Math. Soc. 34}, 6 (2002), 667--676.

\bibitem{LangPlaut}
{\sc Lang, U., and Plaut, C.}
\newblock Bilipschitz embeddings of metric spaces into space forms.
\newblock {\em Geom. Dedicata 87}, 1-3 (2001), 285--307.

\bibitem{Leeslash}
{\sc Lee, J.~R., and Raghavendra, P.}
\newblock Coarse differentiation and multi-flows in planar graphs.
\newblock {\em Discrete Comput. Geom. 43}, 2 (2010), 346--362.

\bibitem{MT}
{\sc Mackay, J.~M., and Tyson, J.~T.}
\newblock {\em Conformal dimension}, vol.~54 of {\em University Lecture
  Series}.
\newblock American Mathematical Society, Providence, RI, 2010.
\newblock Theory and application.

\bibitem{malo2015discrete}
{\sc Malo, R.~J.}
\newblock {\em Discrete Extremal Lengths of Graph Approximations of Sierpiński
  Carpets}.
\newblock Ph.d. thesis, Montana State University, Ann Arbor, MI, 2015.

\bibitem{Mcgillivray}
{\sc Mcgillivray, I.}
\newblock Resistance in higher-dimensional {S}ierpiński carpets.
\newblock {\em Potential Analysis 16\/} (2002), 289--303.

\bibitem{merenkov}
{\sc Merenkov, S.}
\newblock A {S}ierpi\'nski carpet with the co-{H}opfian property.
\newblock {\em Invent. Math. 180}, 2 (2010), 361--388.

\bibitem{meyerquasi}
{\sc Meyer, D.}
\newblock Quasisymmetric embedding of self similar surfaces and origami with
  rational maps.
\newblock {\em Ann. Acad. Sci. Fenn. Math. 27}, 2 (2002), 461--484.

\bibitem{MurCA}
{\sc Murugan, M.}
\newblock Conformal {A}ssouad dimension as the critical exponent for
  combinatorial modulus.
\newblock {\em Ann. Fenn. Math. 48}, 2 (2023), 453--491.

\bibitem{murugan2023first}
{\sc Murugan, M., and Shimizu, R.}
\newblock First-order {S}obolev spaces, self-similar energies and energy
  measures on the {S}ierpi{\'n}ski carpet.
\newblock {\em Comm. Pure Appl. Math. (in press)\/} (2025).

\bibitem{NakamuraYamasakiDuality}
{\sc Nakamura, T., and Yamasaki, M.}
\newblock {Generalized extremal length of an infinite network}.
\newblock {\em Hiroshima Mathematical Journal 6}, 1 (1976), 95 -- 111.

\bibitem{neroli2024fractal}
{\sc Neroli, Z.}
\newblock Fractal dimensions for iterated graph systems.
\newblock {\em Proceedings of the Royal Society A 480}, 2300 (2024), 20240406.

\bibitem{P89}
{\sc Pansu, P.}
\newblock Dimension conforme et sph\`ere \`a l'infini des vari\'{e}t\'{e}s \`a
  courbure n\'{e}gative.
\newblock {\em Ann. Acad. Sci. Fenn. Ser. A I Math. 14}, 2 (1989), 177--212.

\bibitem{Carthesis}
{\sc Piaggio, M.~C.}
\newblock {\em Jauge conforme des espaces m{\'e}triques compacts}.
\newblock PhD thesis, Universit{\'e} de Provence-Aix-Marseille I, 2011.

\bibitem{Qi2017ConsensusIS}
{\sc Qi, Y., Zhang, Z., Yi, Y., and Li, H.}
\newblock Consensus in self-similar hierarchical graphs and {S}ierpiński
  graphs: {C}onvergence speed, delay robustness, and coherence.
\newblock {\em IEEE Transactions on Cybernetics 49\/} (2017), 592--603.

\bibitem{AndreaSchioppa2015}
{\sc Schioppa, A.}
\newblock Poincaré inequalities for mutually singular measures.
\newblock {\em Analysis and Geometry in Metric Spaces 3}, 1 (2015), 40--45,
  electronic only.

\bibitem{Shaconf}
{\sc Shanmugalingam, N.}
\newblock On {C}arrasco {P}iaggio's theorem characterizing quasisymmetric maps
  from compact doubling spaces to {A}hlfors regular spaces.
\newblock In {\em Potentials and partial differential equations---the legacy of
  {D}avid {R}. {A}dams}, vol.~8 of {\em Adv. Anal. Geom.} De Gruyter, Berlin,
  [2023] \copyright 2023, pp.~23--48.

\bibitem{ShimizuParabolicIndex}
{\sc Shimizu, R.}
\newblock {\em Parabolic index of an infinite graph and Ahlfors regular
  conformal dimension of a self-similar set}.
\newblock De Gruyter, Berlin, Boston, 2021, pp.~201--274.

\bibitem{shimizu}
{\sc Shimizu, R.}
\newblock Construction of {$p$}-energy and associated energy measures on
  {S}ierpi\'{n}ski carpets.
\newblock {\em Trans. Amer. Math. Soc. 377}, 2 (2024), 951--1032.

\bibitem{XI2017ScaleFree}
{\sc Xi, L., Wang, L., Wang, S., Yu, Z., and Wang, Q.}
\newblock Fractality and scale-free effect of a class of self-similar networks.
\newblock {\em Physica A: Statistical Mechanics and its Applications 478\/}
  (2017), 31--40.

\end{thebibliography}

\end{document}